\documentclass[a4paper]{amsart}
\usepackage{amscd, amssymb, color, booktabs}
\usepackage{myown}
\title[Lusztig Correspondence and Howe Correspondence]
{Lusztig Correspondence and Howe Correspondence for Finite Reductive Dual Pairs}
\author{Shu-Yen Pan}
\address{Department of Mathematics,
National Tsing Hua University, Hsinchu 300, Taiwan}
\email{sypan@math.nthu.edu.tw}


\date{\today}

\begin{document}

\begin{abstract}
Let $(G,G')$ be a reductive dual pair of a symplectic group and an orthogonal group over a finite field
of odd characteristic.
The Howe correspondence establishes a correspondence between a subset of irreducible characters of $G$
and a subset of irreducible characters of $G'$.
The Lusztig correspondence is a bijection between the Lusztig series indexed by the conjugacy class
of a semisimple element $s$ in the connected component $(G^*)^0$ of the dual group of $G$
and the set of irreducible unipotent characters of the centralizer of $s$ in $G^*$.
In this paper, we prove the commutativity (up to a twist of the sign character)
between these two correspondences.
As a consequence, the Howe correspondence can be explicitly described in terms of Lusztig's
parametrization for classical groups.
\end{abstract}

\maketitle
\tableofcontents


\section{Introduction}

\subsection{}
Let $(\bfG,\bfG')$ be a reductive dual pair over a finite field $\bfF_q$ of odd characteristic.
By restricting the Weil character (\cf.~\cite{gerardin}) to $(\bfG,\bfG')$ with respect to
a nontrivial additive character $\psi$ of $\bfF_q$,
we have the \emph{Weil character} $\omega^\psi_{\bfG,\bfG'}$,
which has non-negative integral decomposition
\[
\omega^\psi_{\bfG,\bfG'}=\sum_{\rho\in\cale(G),\ \rho'\in\cale(G')}m_{\rho,\rho'}\rho\otimes\rho'
\]
where $\cale(G)$ denotes the set of irreducible characters of group of rational points $G$ of $\bfG$.
We say that $(\rho,\rho')$ \emph{occurs} in the \emph{Howe correspondence} if $m_{\rho,\rho'}\neq 0$,
i.e., there is a relation
\[
\Theta^\psi_{\bfG,\bfG'}=\{\,(\rho,\rho')\in\cale(G)\times\cale(G')\mid m_{\rho,\rho'}\neq 0\,\}
\]
between $\cale(G)$ and $\cale(G')$.
The main task of this article is to describe the correspondence explicitly in terms of
Lusztig parametrization.

\subsection{}
Let $\cale(G)_s$ denote the \emph{Lusztig series} associated to the conjugacy class $(s)$
of a semisimple element $s$ in the connected component $(G^*)^0$ of the dual group $G^*$ of $G$.
In particular, $\cale(G)_1$ is the set of irreducible unipotent characters of $G$.
It is known that these Lusztig series partition the set of irreducible characters:
\[
\cale(G)=\bigcup_{(s)\subset (G^*)^0}\cale(G)_s.
\]
Irreducible characters in $\cale(G)_1$ are called \emph{unipotent}.
In \cite{lg} Lusztig shows that there exists a bijection:
\[
\grL_s\colon\cale(G)_s\longrightarrow\cale(C_{G^*}(s))_1
\]
where $C_{G^*}(s)$ denotes the centralizer of $s\in(G^*)^0$.
Such a bijection $\grL_s$ will be called the \emph{Lusztig correspondence}.
Note that $\grL_s$ is in general not uniquely determined.
A discussion on the ambiguity of $\grL_s$ is in \cite{pan-ambiguity}.

If $(\bfG,\bfG')=(\Sp_{2n},\rmO^\epsilon_{2n'})$ for some $n,n'$,
it is known that the unipotent characters are preserved by the Howe correspondence (\cf.~\cite{adams-moy} theorem 3.5).
For this case, the restriction of $\Theta^\psi_{\bfG,\bfG'}$ to unipotent characters is denoted by
$\Theta_{\bfG,\bfG',1}$.
Moreover, it is known in \cite{lg} that irreducible unipotent characters of
a symplectic group or an orthogonal group can be parametrized by
\emph{symbols} satisfying certain conditions (\cf.~Proposition~\ref{0232}).
We can define an explicit relation $\calb_{\bfG,\bfG'}$ on those symbols (\cf.~(\ref{0227})).
The following proposition is \cite{pan-finite-unipotent} theorem 3.4:

\begin{prop*}
Let $(\bfG,\bfG')=(\Sp_{2n},\rmO^\epsilon_{2n'})$.
The unipotent part $\omega_{\bfG,\bfG',1}$ of the Weil character $\omega^\psi_{\bfG,\bfG'}$ has the decomposition:
\begin{align*}
\omega_{\bfG,\bfG',1}
&=\sum_{(\Lambda,\Lambda')\in\calb_{\bfG,\bfG'}}\rho_\Lambda\otimes\rho_{\Lambda'}
\end{align*}
where $\rho_\Lambda$ denotes the irreducible unipotent character parametrized by the
symbol $\Lambda$.
\end{prop*}

The above proposition means that we can give a very explicit description of
the finite Howe correspondence on unipotent characters.
Our next step is to reduce the general Howe correspondence to the correspondence of
unipotent characters via the Lusztig correspondence.

\subsection{}
For a semisimple element $s$ in a symplectic group or an orthogonal group,
we can define three groups $G^{(0)}$, $G^{(1)}$ and $G^{(2)}$ depending on $(s)$,
a decomposition $s=s^{(0)}\times s^{(1)}\times s^{(2)}$ (\cf.~(\ref{0613})),
and a \emph{modified Lusztig correspondence}
\[
\Xi_s\colon\cale(G)_s\rightarrow
\cale(G^{(0)}\times G^{(1)}\times G^{(2)})_1
\]
(\cf.~Subsection~\ref{0506} and Subsection~\ref{0607}).

Then we have the following result on the (almost) commutativity between
the modified Lusztig correspondences $\Xi_s,\Xi_{s'}$ and Howe correspondence
$\Theta^\psi_{\bfG,\bfG'}$ (\cf.~Theorem~\ref{0606}):

\begin{thm*}
Let $(\bfG,\bfG')=(\Sp_{2n},\rmO^\epsilon_{2n'})$,
and let $\rho\in\cale(G)_s$ and $\rho'\in\cale(G')_{s'}$ for some semisimple elements $s\in G^*$
and $s'\in (G'^*)^0$.
Write $\Xi_s(\rho)=\rho^{(0)}\otimes\rho^{(1)}\otimes\rho^{(2)}$ and\/
$\Xi_{s'}(\rho')=\rho'^{(0)}\otimes\rho'^{(1)}\otimes\rho'^{(2)}$.
Let $\{\rho'_i\}$ denote the image under $\Xi_{s'}^{-1}$ of
\begin{multline*}
\{\rho'^{(0)}\otimes\rho'^{(1)}\otimes\rho'^{(2)},\rho'^{(0)}\otimes(\rho'^{(1)}\cdot\sgn)\otimes\rho'^{(2)}, \\
\rho'^{(0)}\otimes\rho'^{(1)}\otimes(\rho'^{(2)}\cdot\sgn),
\rho'^{(0)}\otimes(\rho'^{(1)}\cdot\sgn)\otimes(\rho'^{(2)}\cdot\sgn)\}.
\end{multline*}
Then one of the $\{\rho\otimes\rho'_i\}$ occurs in $\omega_{\bfG,\bfG'}^\psi$ if and only if
the following conditions hold:
\begin{itemize}
\item $s^{(0)}=s'^{(0)}$ (up to conjugation), and $\rho^{(0)}=\rho'^{(0)}$;

\item $G^{(1)}\simeq G'^{(1)}$, and $\rho^{(1)}$ is equal to $\rho'^{(1)}$ or $\rho'^{(1)}\cdot\sgn$;

\item $\rho^{(2)}\otimes\rho'^{(2)}$ or $\rho^{(2)}\otimes(\rho'^{(2)}\cdot\sgn)$ occurs in
$\omega_{\bfG^{(2)},\bfG'^{(2)},1}$.
\end{itemize}
\end{thm*}
The theorem is a little ambiguous due to the non-uniqueness of the modified Lusztig correspondences
$\Xi_s,\Xi_{s'}$.
If $\Xi_s,\Xi_{s'}$ are properly chosen as in \cite{pan-ambiguity},
we can have a more precise result (\cf.~Proposition~\ref{0808}):

\begin{prop*}
Let $(\bfG,\bfG')=(\Sp_{2n},\rmO^\epsilon_{2n'})$,
and let $\rho\in\cale(G)_s$ and $\rho'\in\cale(G')_{s'}$ for some semisimple elements $s\in G^*$
and $s'\in (G'^*)^0$.
Let $\Xi_s,\Xi_{s'}$ be the unique Lusztig correspondences given in \cite{pan-ambiguity} \S 9, \S 8 respectively.
Write $\Xi_s(\rho)=\rho^{(0)}\otimes\rho^{(1)}\otimes\rho^{(2)}$ and\/
$\Xi_{s'}(\rho')=\rho'^{(0)}\otimes\rho'^{(1)}\otimes\rho'^{(2)}$.
Then $\rho\otimes\rho'$ occurs in $\omega_{\bfG,\bfG'}^\psi$ if and only if the following conditions hold:
\begin{itemize}
\item $s^{(0)}=s'^{(0)}$ (up to conjugation), and $\rho^{(0)}=\rho'^{(0)}$;

\item $G^{(1)}\simeq G'^{(1)}$, and $\rho^{(1)}=\rho'^{(1)}$;

\item $\rho^{(2)}\otimes\rho'^{(2)}$ occurs in $\omega_{\bfG^{(2)},\bfG'^{(2)},1}$.
\end{itemize}
Therefore, we have the commutative diagram
\[
\begin{CD}
\cale(G)_s @> \Theta^\psi_{\bfG,\bfG'} >> \cale(G')_{s'} \\
@V \Xi_s VV @VV \Xi_{s'} V \\
\cale(G^{(0)}\times G^{(1)}\times G^{(2)})_1
@> \id\otimes\id\otimes\Theta_{\bfG^{(2)},\bfG'^{(2)},1} >> \cale(G'^{(0)}\times G'^{(1)}\times G'^{(2)})_1
\end{CD}
\]
\end{prop*}
We also have analogous results for dual pair $(\Sp_{2n},\SO_{2n'+1})$
(\cf.~Theorem~\ref{0502} and Proposition~\ref{0803}).
Therefore, the description of general Howe correspondence
for dual pair of a symplectic group and an orthogonal group is now completely obtained.

For a dual pair $(\rmU_n,\rmU_{n'})$ of two unitary groups the above results also hold and are known
by Aubert-Michel-Rouquier in \cite{amr}.
The main difference between the correspondence for $(\rmU_n,\rmU_{n'})$ and the correspondences
for $(\Sp_{2n},\rmO^\epsilon_{2n'})$ or $(\Sp_{2n},\SO_{2n'+1})$ is 
that all irreducible unipotent characters of unitary groups are uniform
(i.e., can be written as linear combinations of Deligne-Lusztig virtual characters)
but $\omega^\psi_{\Sp_{2n},\rmO^\epsilon_{2n'}}$, $\omega^\psi_{\Sp_{2n},\SO_{2n'+1}}$ 
and most of irreducible characters of a symplectic group or an orthogonal group are not.
Therefore we need to make more efforts to characterize carefully an irreducible character from
its uniform projection.

\subsection{}
The contents of this article are as follows.
In Section 2 we first give the definition and notation which are used in this paper.
Then we recall some basic properties of Deligne-Lusztig virtual character and the characterization
of unipotent characters by Lusztig.
In Section 3 we provide some results on uniform projection of an irreducible character of a symplectic group
or an orthogonal group.
In Section 4 we focus on the Lusztig correspondence and the decomposition of the uniform projection of the Weil character.
In Section~5 for completeness we rephrase the result by Aubert-Michel-Rouquier on the Howe correspondence of irreducible
unipotent characters for the dual pair of two unitary groups,
and we also recall the commutativity between the Howe correspondence and Lusztig correspondence for this case.
In Section 6 we prove our main result: Theorem~\ref{0606} on the almost commutativity between Howe correspondence
and Lusztig correspondence on irreducible characters for the dual pair of a symplectic group and an
even orthogonal group.
In Section 7 we consider the analogous result (Theorem~\ref{0502}) for the dual pair of
a symplectic group and a special odd orthogonal group.
In the final section, we recall some results from \cite{pan-ambiguity} on the choice
of the modified Lusztig correspondences.
And then we obtain the full commutativity between Howe correspondence and modified Lusztig correspondence.


\section{Unipotent Characters}

\subsection{Basic notations}\label{0205}
Let $\bfF_q$ be a finite field of $q$ elements where $q$ is a power of an odd prime.
Let $\bfG$ be a classical group defined over $\bfF_q$,
$F$ the corresponding Frobenius endomorphism,
$G=\bfG^F$ the group of rational points.

Let $\cale(G)$ denote the set of the characters of irreducible representations of $G$,
and let $\calv(G)$ denote the space of complex-valued class functions on $G$.
It is known that $\calv(G)$ is an inner product space
(with the inner product denoted by $\langle,\rangle_\bfG$)
with $\cale(G)$ as an orthonormal basis.
Let ${\bf1}={\bf1}_\bfG$ denote the trivial character of $G$,
and let $\sgn=\sgn_\bfG$ denote the sign character if $G$ is an orthogonal group.

Let $\bfT_0$ be a fixed maximally split rational maximal torus of $\bfG$,
$\bfN(\bfT_0)$ the normalizer of $\bfT_0$,
$\bfW=\bfW_\bfG=\bfN(\bfT_0)/\bfT_0$ the \emph{Weyl group} of $\bfG$.
For $w\in\bfW$, choose an element $g\in\bfG$ such that $g^{-1}F(g)\in\bfN(\bfT_0)$ and
its image in $\bfW$ is $w$,
and define $\bfT_w=g\bfT_0g^{-1}$,
a rational maximal torus of $\bfG$.
Every rational maximal torus of $\bfG$ is $G$-conjugate to $\bfT_w$ for some $w\in\bfW$.

The Weyl groups $\bfW_{\Sp_{2n}}$ and $\bfW_{\SO_{2n+1}}$ is identified with the group $W_n$
of permutations on the set $\{1,2,\ldots,n,n^*,(n-1)^*,\ldots,1^*\}$ which
commutes with the involution
\[
(1,1^*)(2,2^*)\cdots(n,n^*)
\]
where $(i,j)$ denote the transposition of $i,j$.
For $i=1,\ldots,n-1$, let $s_i=(i,i+1)(i^*,(i+1)^*)$ and let $\sigma_n=(n,n^*)$.
It is known that $W_n$ is generated by $\{s_1,\ldots,s_{n-1},\sigma_n\}$.
The kernel $W_n^+$ of the homomorphism $\varepsilon\colon W_n\rightarrow\{\pm1\}$
given by $s_i\mapsto 1$ and $\sigma_n\mapsto -1$ is subgroup of index two and is generated
by $\{s_1,\ldots,s_{n-1},\sigma_n s_{n-1}\sigma_n\}$.
Let $W_n^-=W_n\smallsetminus W_n^+$.
The mapping $W^+_n\rightarrow W_n^-$ given by $x\mapsto x\sigma_n$ is a bijection.
The Weyl group $\bfW_{\SO_{2n}^\epsilon}$ where $\epsilon=\pm$ is identified with $W_n^+$.

For a finite set $X$,
let $|X|$ denote the number of elements of $X$.

\subsection{Centralizer of a semisimple element}\label{0501}
From \cite{amr} subsection~1.B, we know that the centralizer in $\bfG$ of a rational semisimple element $s\in G^0$
can be described as follows.
Suppose that $\bfG$ is a classical group of rank $l$,
and $\bfT_l\simeq \overline\bfF_q\times\cdots\times\overline\bfF_q$
is a rational maximal torus where $\overline\bfF_q$ denotes a fixed algebraic closure of $\bfF_q$.
For $s=(\lambda_1,\ldots,\lambda_l)\in T_l$,
let $\nu_\lambda(s)$ denote the number of the $\lambda_i$'s which are equal to $\lambda$,
and let $\langle\lambda\rangle$ denote the set of all roots in $\overline\bfF_q$ of the irreducible polynomial of $\lambda$ over $\bfF_q$.
Now the group $C_\bfG(s)$ decomposed as a product
\[
C_{\bfG}(s)=\prod_{\langle\lambda\rangle\subset\{\lambda_1,\ldots,\lambda_l\}}\bfG_{[\lambda]}(s)
\]
where $\bfG_{[\lambda]}(s)$ is a reductive quasi-simple group
of rank equal to $|\langle\lambda\rangle|\nu_\lambda(s)$.
\begin{itemize}
\item If $\bfG=\rmU_l$ or if $\lambda\neq \pm 1$,
then $\bfG_{[\lambda]}(s)$ is either a general linear group or a unitary group.

\item If $\bfG=\SO_{2l+1}$,
then $\bfG_{[-1]}(s)\simeq\rmO_{2\nu_{-1}(s)}^\pm$ and $\bfG_{[1]}(s)\simeq\SO_{2\nu_1(s)+1}$

\item If $\bfG=\rmO_{2l}^\pm$,
then $\bfG_{[-1]}(s)\simeq\rmO_{2\nu_{-1}(s)}^\pm$ and $\bfG_{[1]}(s)\simeq\rmO_{2\nu_1(s)}^\pm$.

\item If $\bfG=\Sp_{2l}$,
then $\bfG_{[-1]}(s)\simeq\Sp_{2\nu_{-1}(s)}$ and $\bfG_{[1]}(s)\simeq\Sp_{2\nu_1(s)}$.
\end{itemize}

\subsection{Deligne-Lusztig virtual characters}\label{0206}
For a rational maximal torus $\bfT$ in a connected group $\bfG$ and $\theta\in\cale(T)$ where $T=\bfT^F$,
let $R_{\bfT,\theta}=R^\bfG_{\bfT,\theta}$ denote the \emph{Deligne-Lusztig (virtual) character} of $G$
defined in \cite{dl}.
For an orthogonal group, we define
\begin{equation}\label{0213}
R_{\bfT,\theta}^{\rmO_n^\epsilon}
=\Ind_{\SO_n^\epsilon}^{\rmO_n^\epsilon}(R_{\bfT,\theta}^{\SO_n^\epsilon})
\end{equation}
An irreducible character $\eta$ of $G$ is \emph{unipotent} if
$\langle\eta,R_{\bfT,\bf1}\rangle_\bfG\neq 0$ for some $\bfT$.
Let $\cale(G)_1$ denote the set of irreducible unipotent characters of $G$.

Now we recall some definitions from \cite{lg-CBMS} 3.17, 3.19.
Suppose that $\bfG$ is connected.
The Frobenius map $F$ acts on the Weyl group $\bfW=\bfW_\bfG$ and on the set $\cale(\bfW)$ of irreducible characters of $\bfW$.
The set of fixed points is denoted by $\cale(\bfW)^F$.
An irreducible representation $(\pi,V)$ of $\bfW$ is in $\cale(\bfW)^F$ if and only if there exists
a vector space isomorphism $\phi\colon V \rightarrow V$ such that the following diagram
\begin{equation}\label{0907}
\begin{CD}
V @>\pi(w)>> V\\
@V \phi VV @VV \phi V\\
V @>\pi(F(w))>> V
\end{CD}
\end{equation}
commutes for all $w\in\bfW$.
The isomorphism $\phi$ is unique up to a nonzero constant multiple.
For $(\pi,V)\in\cale(\bfW)^F$,
we define
\begin{equation}\label{0908}
R_\pi^\bfG=\frac{1}{|\bfW|}\sum_{w\in\bfW}\Tr(\phi\circ\pi(w),V)R_{\bfT_w,\bf1}^\bfG.
\end{equation}
It is known that these $R_\pi^\bfG$'s, for $(\pi,V)\in\cale(\bfW)^F$,
form an orthonormal basis of the space $\calv(G)_1^\sharp$ of unipotent uniform class functions on $G$.

If $\bfG$ is $\Sp_{2n}$ or $\SO_{2n}^+$,
the action of $F$ is trivial,
then $\phi$ is trivial and for an irreducible character $\chi$ of $\bfW=W_n$ or $W_n^+$,
we have
\begin{align*}
R^{\Sp_{2n}}_\chi
&=\frac{1}{|W_n|}\sum_{w\in W_n}\chi(w) R_{\bfT_w,\bf1}^{\Sp_{2n}};\\
R_\chi^{\SO_{2n}^+}
&=\frac{1}{|W_n^+|}\sum_{w\in W_n^+}\chi(w) R_{\bfT_w,\bf1}^{\SO_{2n}^+}.
\end{align*}

Now we consider that case $\bfG=\SO_{2n}^-$.
Identify $\bfW_\bfG=W_n^+$.
The action $F\colon W_n^+\rightarrow W_n^+$ is given by $x\mapsto\sigma_nx\sigma_n$ where $\sigma_n$ is given in
Subsection~\ref{0205}.
If $\chi\in\cale(W_n)$,
then the restriction of $\chi$ to $W_n^+$ is fixed by $F$ (but might not be irreducible)
and the isomorphism $\phi$ in (\ref{0907}) can be chosen to be $\chi(\sigma_n)$,
therefore (\ref{0908}) becomes
\[
R^{\SO_{2n}^-}_\chi=\frac{1}{|W_n^+|}\sum_{w\in W_n^+}\chi(w\sigma_n)R_{\bfT_w,\bf1}^{\SO_{2n}^-}.
\]
Now the mapping $w\mapsto w\sigma_n$ gives a bijection between $W_n^+$ and $W_n^-$,
and we define $R_{\bfT_{w\sigma_n},1}^{\SO_{2n}^-}=R_{\bfT_w,1}^{\SO_{2n}^-}$ for $w\in W_n^+$ by convention.
So finally we have
\[
R^{\SO_{2n}^-}_\chi=\frac{1}{|W_n^-|}\sum_{w\in W_n^-}\chi(w)R_{\bfT_w,\bf1}^{\SO_{2n}^-}
\]
for any $\chi\in\cale(W_n)$.

\subsection{Symbols and unipotent characters}\label{0234}
In this subsection,
we recall some basic definitions and properties of ``symbols'' introduced by Lusztig \cite{lg}.

A \emph{$\beta$-set} $A=\{a_1,\ldots,a_m\}$ is a finite subset (possibly empty) of non-negative integers
written in (strictly) decreasing order, i.e., $a_1>a_2>\cdots>a_m$.
A \emph{symbol}
\[
\Lambda=\binom{A}{B}=\binom{a_1,a_2,\ldots,a_{m_1}}{b_1,b_2,\ldots,b_{m_2}}
\]
is an ordered pair of two $\beta$-sets.
The first row (resp.~second row) of $\Lambda$ is also denoted by $\Lambda^*$ (resp.~$\Lambda_*$).
A symbol $\Lambda$ is called \emph{reduced} if $0\not\in\Lambda^*\cap\Lambda_*$;
and $\Lambda$ is called \emph{degenerate} if $\Lambda^*=\Lambda_*$, and \emph{nondegenerate} otherwise.
The \emph{rank} and the \emph{defect} of $\Lambda$ are defined by
\begin{align*}
{\rm rk}(\Lambda) &=\sum_{a\in\Lambda^*} a+\sum_{b\in\Lambda_*}b-\left\lfloor\left(\frac{|\Lambda^*|+|\Lambda_*|-1}{2}\right)^2\right\rfloor, \\
{\rm def}(\Lambda) &=|\Lambda^*|-|\Lambda_*|.
\end{align*}
Moreover, we denote $|\Lambda|=|\Lambda^*|+|\Lambda_*|$.
Note that the definition of ${\rm def}(\Lambda)$ here is slightly different from the original
definition in \cite{lg}.
It is easy to check that
\begin{equation}\label{0210}
{\rm rk}(\Lambda)\geq\left\lfloor\left(\frac{{\rm def}(\Lambda)}{2}\right)^2\right\rfloor.
\end{equation}
A symbol $\Lambda$ is called \emph{cuspidal} if
${\rm rk}(\Lambda)=\bigl\lfloor\bigl(\frac{{\rm def}(\Lambda)}{2}\bigr)^2\bigr\rfloor$.
It is not difficult to check that $\Lambda$ is cuspidal if and only if
$\Lambda^*=\{m_1,m_1-1,\ldots,0\}$ and $\Lambda_*=\{m_2,m_2-1,\ldots,0\}$ for some $m_1,m_2$.

Define $\binom{A}{B}^\rmt=\binom{B}{A}$ which is called the \emph{transpose} of $\binom{A}{B}$.
It is clear that ${\rm rk}(\Lambda^\rmt)={\rm rk}(\Lambda)$ and ${\rm def}(\Lambda^\rmt)=-{\rm def}(\Lambda)$.
On the set of symbols we define an equivalence relation generated by
\[
\binom{A}{B}\sim
\binom{\{\,a+1\mid a\in A\,\}\cup\{0\}}{\{\,b+1\mid b\in B\,\}\cup\{0\}}.
\]
It is not difficult to see that two equivalent symbols have the same rank and the same defect.
Moreover, each symbol is equivalent to a unique reduced symbol.

Let $\cals$ denote the set of all (reduced) symbols,
and let $\cals_{n,\delta}\subset\cals$ denote the subset of symbols of rank $n$ and defect $\delta$.
And we define
\begin{align*}
\cals_\Sp &=\{\,\Lambda\in\cals\mid{\rm def}(\Lambda)\equiv 1\pmod 4\,\}, &
\cals_{\Sp_{2n}} &=\{\,\Lambda\in\cals_\Sp\mid{\rm rk}(\Lambda)=n\,\}; \\
\cals_{\rmO^+} &=\{\,\Lambda\in\cals\mid{\rm def}(\Lambda)\equiv 0\pmod 4\,\}, &
\cals_{\rmO_{2n}^+} &=\{\,\Lambda\in\cals_{\rmO^+}\mid{\rm rk}(\Lambda)=n\,\}; \\
\cals_{\rmO^-} &=\{\,\Lambda\in\cals\mid{\rm def}(\Lambda)\equiv 2\pmod 4\,\}, &
\cals_{\rmO_{2n}^-} &=\{\,\Lambda\in\cals_{\rmO^-}\mid{\rm rk}(\Lambda)=n\,\}.
\end{align*}
Then we have the following parametrization of irreducible unipotent characters from
\cite{lg} theorem 8.2:

\begin{prop}[Lusztig]\label{0232}
Let\/ $\bfG$ be a symplectic group or an even orthogonal group.
Then there exists a bijective correspondence between $\cals_\bfG$ and $\cale(G)_1$.
\end{prop}
The irreducible unipotent character associated to the symbol $\Lambda\in\cals_\bfG$ is denoted by $\rho_\Lambda$.
It is known that $\rho_{\Lambda^\rmt}=\rho_\Lambda\cdot\sgn$ when $\bfG$ is an even orthogonal group.

\subsection{Partitions and bi-partitions}\label{0229}
For a partition $\lambda=[\lambda_1,\ldots,\lambda_r]$ with $\lambda_1\geq\lambda_2\geq\cdots\geq\lambda_r>0$,
we define $\|\lambda\|=\lambda_1+\cdots+\lambda_r$ and $\ell(\lambda)=r$.

For a $\beta$-set $A=\{a_1,\ldots,a_m\}$, we define
\[
\Upsilon\colon\{a_1,\ldots,a_m\}\mapsto[a_1-(m-1),a_2-(m-2),\ldots,a_{m-1}-1,a_m].
\]
Then it is easy to check that the map
\begin{equation}\label{0219}
\Upsilon(\Lambda):=\sqbinom{\Upsilon(\Lambda^*)}{\Upsilon(\Lambda_*)}
\end{equation}
induces a bijection
\[
\Upsilon\colon\cals_{n,\delta}\longrightarrow
\textstyle\calp_2\left(n-\bigl\lfloor\bigl(\frac{\delta}{2}\bigr)^2\bigr\rfloor\right)
\]
where $\calp_2(n)$ denotes the set of \emph{bi-partitions} of $n$,
i.e., the set of ordered pair $\sqbinom{\lambda}{\mu}$ of two partitions $\lambda,\mu$ such that
$\left\|\sqbinom{\lambda}{\mu}\right\|:=\|\lambda\|+\|\mu\|=n$.
In particular, $\Upsilon$ induces a bijection from $\cals_{n,1}$ onto $\calp_2(n)$ and
a bijection from $\cals_{n,0}$ onto $\calp_2(n)$.

It is known that $\cale(W_n)$ is parametrized by $\calp_2(n)$ (\cf.~\cite{Geck-Pfeiffer} theorem 5.5.6)
and so the irreducible character associated to $\sqbinom{\lambda}{\mu}$ is denoted by $\chi_{\sqbinom{\lambda}{\mu}}$.
If $\Sigma$ is in $\cals_{n,1}$ or $\cals_{n,0}$ so that $\Upsilon(\Sigma)=\sqbinom{\lambda}{\mu}$,
then $R^\bfG_{\chi_{\sqbinom{\lambda}{\mu}}}$ (\cf.~Subsection~\ref{0206}) is also denoted by $R^\bfG_\Sigma$.

\subsection{Unipotent characters of\/ $\Sp_{2n}(q)$}\label{0910}
In this subsection, let $\bfG=\Sp_{2n}$.
A symbol $Z=\binom{a_1,a_2,\ldots,a_{m+1}}{b_1,b_2\ldots,b_m}$ of defect $1$
is called \emph{special} if
\[
a_1\geq b_1\geq a_2\geq b_2\geq\cdots\geq a_m\geq b_m\geq a_{m+1}.
\]
Let $Z$ be a special symbol of defect $1$.
Let $Z_\rmI=Z\smallsetminus\binom{Z^*\cap Z_*}{Z^*\cap Z_*}$ be the subsymbol of ``\emph{singles}'' in $Z$.
The \emph{degree} $\deg(Z)$ of a special symbol $Z$ of defect $1$ is defined to be the non-negative integer
$\frac{|Z_\rmI|-1}{2}$.
For a subsymbol $M\subset Z_\rmI$, we define
\begin{align*}
\Lambda_M &= (Z\smallsetminus M)\cup M^\rmt, \\
\cals_Z &=\{\,\Lambda_M\mid M\subset Z_\rmI,\ |M|\text{ even}\,\}, \\
\cals_{Z,1} &= \{\,\Lambda_M\mid M\subset Z_\rmI,\ |M^*|=|M_*|\,\}.
\end{align*}
It is clear that $\cals_Z\subset\cals_\Sp$.

For a special symbol $Z=\binom{a_1,a_2,\ldots,a_{m+1}}{b_1,b_2,\ldots,b_m}$ of defect $1$,
we define (\cf.~\cite{lg-symplectic} (2.5.2))
\begin{multline}\label{0913}
a_Z=\sum_{1\leq i<j\leq m+1}\min(a_i,a_j)+\sum_{1\leq i<j\leq m}\min(b_i,b_j)\\
+\sum_{1\leq i\leq m+1,\ 1\leq j\leq m}\min(a_i,b_j)-\tfrac{1}{6}m(m-1)(4m+1).
\end{multline}
For a symbol $\Lambda\in\cals_\Sp$,
let $a_\Lambda=a_Z$ if $\Lambda\in\cals_Z$ for some special symbol $Z$ of defect $1$.

\begin{lem}\label{0912}
Let $Z=\binom{a_1,a_2,\ldots,a_{m+1}}{b_1,b_2,\ldots,b_m}$ be a special symbol of defect $1$.
If $m\geq 1$,
then $a_Z>0$.
\end{lem}
\begin{proof}
Let
\[
Z_0=\binom{m,m-1,\ldots,1,0}{m,m-1,\ldots,1}\quad\text{or}\quad\binom{m+1,m,\ldots,1}{m-1,m-2,\ldots,0}.
\]
It is easy to check that $a_{Z_0}=m^2$ from (\ref{0913}).
Now for a special symbol $Z=\binom{a_1,a_2,\ldots,a_{m+1}}{b_1,b_2,\ldots,b_m}$,
we know that either
\begin{itemize}
\item $a_i\geq i-1$ and $b_j\geq j$ for $i=1,\ldots,m+1$ and $j=1,\ldots, m$; or

\item $a_i\geq i$ and $b_j\geq j-1$ for $i=1,\ldots,m+1$ and $j=1,\ldots, m$.
\end{itemize}
Hence it is easy to check that $a_Z\geq a_{Z_0}=m^2\geq1$ form (\ref{0913}).
\end{proof}

The following result is from \cite{lg-symplectic} theorem 5.8.

\begin{prop}[Lusztig]\label{0901}
Let $\bfG=\Sp_{2n}$, $Z$ a special symbol of rank $n$ and defect $1$.
For $\Sigma\in\cals_{Z,1}$ and $\Lambda\in\cals_{\Sp_{2n}}$,
we have
\[
\langle R^\bfG_\Sigma,\rho_{\Lambda}\rangle_\bfG
=\begin{cases}
\frac{1}{2^{\deg(Z)}}(-1)^{\langle\Sigma,\Lambda\rangle}, & \text{if $\Lambda\in\cals_Z$};\\
0, & \text{otherwise}
\end{cases}
\]
where
$\langle,\rangle\colon\cals_{Z,1}\times\cals_Z\rightarrow\bfF_2$ is given by
$\langle\Lambda_N,\Lambda_M\rangle=|N\cap M|\pmod 2$.
\end{prop}
Proposition~\ref{0901} means that we have decompositions
\[
\calv(G)_1=\bigoplus_Z\calv(G)_Z,\qquad
\calv(G)^\sharp_1=\bigoplus_Z\calv(G)^\sharp_Z
\]
where $Z$ runs over all special symbols of rank $n$ and defect $1$,
$\calv(G)_Z$ is the span of $\{\,\rho_\Lambda\mid\Lambda\in\cals_Z\,\}$
and $\calv(G)^\sharp_Z$ is the span of $\{\,R^\bfG_\Sigma\mid\Sigma\in\cals_{Z,1}\,\}$.

Now we recall the concept of ``\emph{cells}'' by Lusztig in \cite{lg-symplectic}, \cite{lg-orthogonal}
(see also \cite{pan-finite-unipotent} \S 4.3, \S 4.4).
Let $Z$ be a special symbol of defect $1$, and let $d=\deg(Z)$.
An \emph{arrangement} of $Z_\rmI$ is a partition $\Phi$
of the $2d+1$ singles in $Z_\rmI$ into $d$ (disjoint) pairs and
one isolated element such that each pair contains one element in the first row and one element in the
second row of $Z_\rmI$.
A set $\Psi$ of some pairs (possibly empty) in $\Phi$ is called a \emph{subset of pairs} of $\Phi$
and is denoted by $\Psi\leq\Phi$.
For a subset of pairs $\Psi$ of an arrangement $\Phi$,
we define
\begin{equation}\label{0302}
C_{\Phi,\Psi}=\{\,\Lambda_M\in\cals_Z
\mid |M\cap\Psi'|\equiv|(\Phi\smallsetminus\Psi)\cap\Psi'^*|\pmod 2\text{ for each }\Psi'\leq\Phi\,\}
\end{equation}
where $\Psi'^*$ denotes the first row of $\Psi'$.
Such a subset $C_{\Phi,\Psi}$ of $\cals_Z$ is called a \emph{cell}.
From (\ref{0302}) we can see that a symbol $\Lambda_M\in\cals_Z$ is in $C_{\Phi,\Psi}$ if and only if
$M$ satisfies the following two conditions:
\begin{itemize}
\item $M$ contains either none or two elements of each pair in $\Psi$; and

\item $M$ contains exactly one element of each pair in $\Phi\smallsetminus\Psi$.
\end{itemize}

The following lemma is from \cite{pan-finite-unipotent} \S 4.3:
\begin{lem}\label{0226}
Let $Z$ be a special symbol of rank $n$ and defect $1$.
\begin{enumerate}
\item[(i)] The class function $\sum_{\Lambda\in C_{\Phi,\Psi}}\rho_\Lambda$ of\/ $\Sp_{2n}(q)$ is uniform.

\item[(ii)] Let $\Lambda_1,\Lambda_2$ be two distinct symbols in $\cals_Z$.
There exists an arrangement $\Phi$ of $Z_\rmI$ with two subsets of pairs $\Psi_1,\Psi_2$
such that $\Lambda_i\in C_{\Phi,\Psi_i}$ for $i=1,2$ and
$C_{\Phi,\Psi_1}\cap C_{\Phi,\Psi_2}=\emptyset$.

\item[(iii)] For any given $\Lambda\in\cals_Z$,
there exist two arrangements $\Phi_1,\Phi_2$ of $Z_\rmI$ with subsets of pairs $\Psi_1,\Psi_2$ respectively such that
$C_{\Phi_1,\Psi_1}\cap C_{\Phi_2,\Psi_2}=\{\Lambda\}$.
\end{enumerate}
\end{lem}

\subsection{Unipotent characters of $\rmO_{2n}^\epsilon(q)$}
In this subsection, let $\bfG=\rmO^\epsilon_{2n}$ where $\epsilon=+$ or $-$.
A symbol $Z=\binom{a_1,a_2,\ldots,a_m}{b_1,b_2\ldots,b_m}$
of defect $0$ is called \emph{special} if
\[
a_1\geq b_1\geq a_2\geq b_2\geq\cdots\geq a_m\geq b_m.
\]
Let $Z$ be a special symbol of defect $0$ and let $Z_\rmI$ be the
subsymbol of singles.
The the \emph{degree} $\deg(Z)$ of $Z$ is defined to be $\frac{|Z_\rmI|}{2}$.
Define
\begin{align*}
\cals^+_Z &=\{\,\Lambda_M\mid M\subset Z_\rmI,\ |M|\text{ even}\,\}, \\
\cals^-_Z &=\{\,\Lambda_M\mid M\subset Z_\rmI,\ |M|\text{ odd}\,\}, \\
\cals_{Z,0} &= \{\,\Lambda_M\mid M\subset Z_\rmI,\ |M^*|=|M_*|\,\}.
\end{align*}
It is clear that $\cals_Z^\epsilon\subset\cals_{\rmO^\epsilon}$.
The following proposition is a modification for $\rmO^\epsilon_{2n}$
from \cite{lg-orthogonal} theorem 3.15 (\cf.~\cite{pan-uniform} proposition 3.6):

\begin{prop}[Lusztig]\label{0315}
Let $\bfG=\rmO^\epsilon_{2n}$, and let $Z$ be a non-degenerate special symbol of rank $n$ and defect $0$.
For any $\Sigma\in\cals_{Z,0}$ and $\Lambda\in\cals_\bfG$,
we have
\[
\langle R^\bfG_\Sigma,\rho_\Lambda\rangle_\bfG
=\begin{cases}
\frac{1}{2^{\deg(Z)-1}}(-1)^{\langle\Sigma,\Lambda\rangle}, & \text{if $\Lambda\in\cals^\epsilon_Z$};\\
0, & \text{otherwise}.
\end{cases}
\]
\end{prop}
Similar to the case of symplectic groups,
we have decompositions
\[
\calv(G)_1=\bigoplus_Z\calv(G)_Z,\qquad
\calv(G)^\sharp_1=\bigoplus_Z\calv(G)^\sharp_Z
\]
where $Z$ runs over all special symbols of  rank $n$ and defect $0$,
$\calv(G)_Z$ is the span of $\{\,\rho_\Lambda\mid\Lambda\in\cals^\epsilon_Z\,\}$
and $\calv(G)^\sharp_Z$ is the span of $\{\,R_\Sigma\mid\Sigma\in\cals_{Z,0}\,\}$.

Let $Z$ be a special symbol of defect $0$ and let $d=\deg(Z)$.
An \emph{arrangement} of $Z_\rmI$ is a partition $\Phi$
of the $2d$ singles in $Z_\rmI$ into $d$ pairs such that each pair contains one element
in the first row and one element in the second row of $Z_\rmI$.
Let the cell $C_{\Phi,\Psi}$ be defined as in (\ref{0302}).
The following two lemmas are from \cite{pan-finite-unipotent} \S 4.4:
\begin{lem}\label{0224}
Let $Z$ be a special symbol of rank $n$ and defect $0$,
and let $\Phi$ be an arrangement of $Z_\rmI$ with subsets of pairs $\Psi,\Psi'$.
\begin{enumerate}
\item[(i)] $\Lambda\in C_{\Phi,\Psi}$ if and only if $\Lambda^\rmt\in C_{\Phi,\Psi}$.

\item[(ii)] We have
\[
\cals^+_Z=\bigcup_{\Psi\leq\Phi,\ \#(\Phi\smallsetminus\Psi)\text{\ even}}C_{\Phi,\Psi}\quad\text{and}\quad
\cals^-_Z=\bigcup_{\Psi\leq\Phi,\ \#(\Phi\smallsetminus\Psi)\text{\ odd}}C_{\Phi,\Psi}
\]
where $\#(\Phi\smallsetminus\Psi)$ means the number of pairs in $\Phi\smallsetminus\Psi$.

\item[(iii)] The class function $\sum_{\Lambda\in C_{\Phi,\Psi}}\rho_\Lambda$ of\/ $\rmO^\epsilon_{2n}(q)$ is uniform.

\item[(iv)] If\/ $\Psi,\Psi'$ are two distinct subsets of pairs of\/ $\Phi$,
then $C_{\Phi,\Psi}\cap C_{\Phi,\Psi'}=\emptyset$.
\end{enumerate}
\end{lem}

A subset of pairs $\Psi$ of an arrangement $\Phi$ is called \emph{admissible} if
$\#(\Phi\smallsetminus\Psi)$ is even when $\epsilon=+$;
and $\#(\Phi\smallsetminus\Psi)$ is odd when $\epsilon=-$.

\begin{lem}\label{0231}
Let $\Lambda_1,\Lambda_2$ be two symbols in $\cals^\epsilon_Z$ such that
$\Lambda_1\neq\Lambda_2,\Lambda_2^\rmt$.
There exists an arrangement $\Phi$ of $Z_\rmI$ with admissible subsets of pairs $\Psi_1,\Psi_2$
such that $\Lambda_i,\Lambda_i^\rmt\in C_{\Phi,\Psi_i}$ for $i=1,2$ and
$C_{\Phi,\Psi_1}\cap C_{\Phi,\Psi_2}=\emptyset$.
\end{lem}

\begin{lem}\label{0228}
Suppose $Z$ is a special symbol with subsymbol of singles
$Z_\rmI=\binom{s_1,s_2,\ldots,s_d}{t_1,t_2,\ldots,t_d}$.
Let $\Phi_1,\Phi_2$ be two arrangements of $Z_\rmI$ defined by
\[
\textstyle\Phi_1=\bigl\{\binom{s_1}{t_1},\binom{s_2}{t_2},\ldots,\binom{s_d}{t_d}\bigr\},
\quad\text{and}\quad
\Phi_2=\bigl\{\binom{s_2}{t_1},\binom{s_3}{t_2},\ldots,\binom{s_d}{t_{d-1}},\binom{s_1}{t_d}\bigr\}
\]
Then for any subsets of pairs $\Psi_1,\Psi_2$ of\/ $\Phi_1,\Phi_2$ respectively such that $\#(\Psi_1)\equiv\#(\Psi_2)\pmod 2$,
we have $|C_{\Phi_1,\Psi_1}\cap C_{\Phi_2,\Psi_2}|=2$.
\end{lem}
\begin{proof}
Let $\Psi_1$ be a subset of pairs in $\Phi_1$,
$\Psi_2'$ a subset of pairs in $\Phi_2\smallsetminus\bigl\{\binom{s_1}{t_d}\bigr\}$,
and $\Psi_2''=\Psi_2'\cup\bigl\{\binom{s_1}{t_d}\bigr\}$.
Then $\#(\Psi_2'')=\#(\Psi_2')+1$.

Suppose that $\Lambda_M$ is in $C_{\Phi_1,\Psi_1}\cap C_{\Phi_2,\Psi_2}$ where $\Psi_2=\Psi_2'$ or $\Psi_2''$
for some $M\subset Z_\rmI$.
From the two conditions before Lemma~\ref{0226}, we have the following:
\begin{enumerate}
\item if $s_i\in M$ and $\binom{s_i}{t_i}\leq\Psi_1$, then $t_i\in M$;

\item if $s_i\not\in M$ and $\binom{s_i}{t_i}\leq\Psi_1$, then $t_i\not\in M$;

\item if $s_i\in M$ and $\binom{s_i}{t_i}\not\leq\Psi_1$, then $t_i\not\in M$;

\item if $s_i\not\in M$ and $\binom{s_i}{t_i}\not\leq\Psi_1$, then $t_i\in M$;

\item if $t_i\in M$ and $\binom{s_{i+1}}{t_i}\leq\Psi_2$, then $s_{i+1}\in M$;

\item if $t_i\not\in M$ and $\binom{s_{i+1}}{t_i}\leq\Psi_2$, then $s_{i+1}\not\in M$;

\item if $t_i\in M$ and $\binom{s_{i+1}}{t_i}\not\leq\Psi_2$, then $s_{i+1}\not\in M$;

\item if $t_i\not\in M$ and $\binom{s_{i+1}}{t_i}\not\leq\Psi_2$, then $s_{i+1}\in M$
\end{enumerate}
for $i=1,\ldots,d$.
This means that for any $\Psi_1,\Psi_2$,
the set $M$ is uniquely determined by the ``initial condition'' whether $s_1$ belongs to $M$ or not
before we apply the final condition ``whether $\binom{s_1}{t_d}\leq\Psi_2$ or not''.
Moreover, before we apply the final condition,
the two possible choices of $M$ are complement subsets to each other in $Z_\rmI$.
Now the final condition is either consistent with the other conditions or contradicts to
the other conditions.
This means that exactly one of the following two situations holds:
\begin{itemize}
\item $|C_{\Phi_1,\Psi_1}\cap C_{\Phi_2,\Psi'_2}|=2$ and $|C_{\Phi_1,\Psi_1}\cap C_{\Phi_2,\Psi''_2}|=0$; or

\item $|C_{\Phi_1,\Psi_1}\cap C_{\Phi_2,\Psi'_2}|=0$ and $|C_{\Phi_1,\Psi_1}\cap C_{\Phi_2,\Psi''_2}|=2$.
\end{itemize}
By (ii) of Lemma~\ref{0224} we know that $C_{\Phi_1,\Psi_1}\cap C_{\Phi_2,\Psi_2}=\emptyset$ if $\#(\Psi_1)\not\equiv\#(\Psi_2)\pmod 2$.
Therefore we must have $|C_{\Phi_1,\Psi_1}\cap C_{\Phi_2,\Psi_2}|=2$ if $\#(\Psi_1)\equiv\#(\Psi_2)\pmod 2$.
\end{proof}

\begin{lem}\label{0225}
For any $\Lambda\in\cals^\epsilon_Z$,
there exist two arrangements $\Phi_1,\Phi_2$ of $Z_\rmI$ with admissible subsets of pairs $\Psi_1,\Psi_2$ respectively
such that $C_{\Phi_1,\Psi_1}\cap C_{\Phi_2,\Psi_2}=\{\Lambda,\Lambda^\rmt\}$.
\end{lem}
\begin{proof}
Let $\Lambda\in\cals^\epsilon_Z$, and let $\Phi_1,\Phi_2$ be the arrangements given in Lemma~\ref{0228}.
Then by Lemma~\ref{0224}, there are admissible subsets of pairs $\Psi_1,\Psi_2$ of $\Phi_1,\Phi_2$
respectively such that $\Lambda,\Lambda^\rmt$ belong to both $C_{\Phi_1,\Psi_1}$ and
$C_{\Phi_2,\Psi_2}$.
Then the lemma follows from Lemma~\ref{0228} immediately.
\end{proof}

An example of the explicit computation for $C_{\Phi_1,\Psi_1}\cap C_{\Phi_2,\Psi_2}$ 
can be found in \cite{pan-finite-unipotent} \S 4.4.

\subsection{Howe correspondence and symbol correspondence}\label{0230}
Let $\lambda=[\lambda_1,\ldots,\lambda_k]$ and $\mu=[\mu_1,\ldots,\mu_l]$ be two partitions with
$\lambda_1\geq\cdots\geq\lambda_k\geq 0$ and $\mu_1\geq\cdots\geq\mu_l\geq 0$.
We denote
\[
\lambda\preccurlyeq\mu\quad\text{ if \ }\mu_1\geq\lambda_1\geq\mu_2\geq\lambda_2\geq\cdots.
\]
Then we define two relations on the set of symbols:
\begin{align*}
\calb^+ &=\{\,(\Lambda,\Lambda')\in\cals\times\cals
\mid\Upsilon(\Lambda_*)\preccurlyeq\Upsilon(\Lambda'^*),\ \Upsilon(\Lambda'_*)\preccurlyeq\Upsilon(\Lambda^*)\,\},\\
\calb^- &=\{\,(\Lambda,\Lambda')\in\cals\times\cals
\mid\Upsilon(\Lambda^*)\preccurlyeq\Upsilon(\Lambda'_*),\ \Upsilon(\Lambda'^*)\preccurlyeq\Upsilon(\Lambda_*)\,\}.
\end{align*}
Moreover, we define
\begin{align}\label{0227}
\begin{split}
\calb_{\Sp,\rmO^+} &=\{\,(\Lambda,\Lambda')\in\calb^+\cap(\cals_\Sp\times\cals_{\rmO^+})\mid {\rm def}(\Lambda')=-{\rm def}(\Lambda)+1\,\}, \\
\calb_{\Sp,\rmO^-} &=\{\,(\Lambda,\Lambda')\in\calb^-\cap(\cals_\Sp\times\cals_{\rmO^-})\mid {\rm def}(\Lambda')=-{\rm def}(\Lambda)-1\,\}, \\
\calb_{\Sp_{2n},\rmO_{2n'}^\epsilon} &=\calb_{\Sp,\rmO^\epsilon}\cap(\cals_{\Sp_{2n}}\times\cals_{\rmO_{2n'}^\epsilon})
\end{split}
\end{align}
where $\epsilon=+$ or $-$.

Recall that $\omega_{\bfG,\bfG',1}$ denotes the unipotent part of Weil character
$\omega^\psi_{\bfG,\bfG'}$ for the dual pair $(\bfG,\bfG')=(\Sp_{2n},\rmO^\epsilon_{2n'})$.
The following proposition is from \cite{pan-finite-unipotent} theorem 3.4:
\begin{prop}\label{0220}
Let $(\bfG,\bfG')=(\Sp_{2n},\rmO^\epsilon_{2n'})$.
We have the decomposition
\[
\omega_{\bfG,\bfG',1}
=\sum_{(\Lambda,\Lambda')\in\calb_{\bfG,\bfG'}}\rho_\Lambda\otimes\rho_{\Lambda'}.
\]
\end{prop}

A similar formulation of the Howe correspondence of irreducible unipotent characters for a unitary dual pair
will be given in Proposition~\ref{0519}.


\section{Uniform Projection of Unipotent Characters}

\subsection{Uniform projection}
Let $G=\bfG^F$ be a finite classical group.
Recall that the space of class functions $\calv(G)$ is an inner product space
with $\cale(G)$ as an orthonormal basis  with respect to the inner product
$\langle,\rangle_\bfG$.
Let $\calv(G)^\sharp$ denote the subspace of $\calv(G)$ spanned by all Deligne-Lusztig virtual characters
$R_{\bfT,\theta}$.
For $f\in\calv(G)$, let $f^\sharp$ denote the orthogonal projection of $f$ over $\calv(G)^\sharp$.
A class function $f\in\calv(G)$ is called \emph{uniform} if $f=f^\sharp$.

If $G$ is connected,
it is well known that the trivial character $\bf1_\bfG$ and the character $\chi_\reg^\bfG$ of regular representation
are both uniform.
For $\rmO_n^\epsilon(q)$,
because $\chi_\reg^{\rmO_n^\epsilon}
=\Ind_{\SO_n^\epsilon(q)}^{\rmO_n^\epsilon(q)}\chi_\reg^{\SO_n^\epsilon}$,
we see that $\chi_\reg^{\rmO_n^\epsilon}$ is uniform by (\ref{0213})
and we have the following identity
(\cf.~\cite{carter-finite} corollary 7.5.6):
\[
\chi_\reg^\bfG
=\frac{|G|_{p'}}{|\bfW_\bfG|}\sum_{w\in\bfW_\bfG}\frac{\epsilon_\bfG\epsilon_{\bfT_w}}{|T_w|}\sum_{\theta\in\cale(T_w)}
R_{\bfT_w,\theta}^\bfG.
\]

\begin{lem}\label{0207}
If $f\in\calv(G)$,
then $f^\sharp(1)=f(1)$.
In particular, if $\rho\in\cale(G)$,
then $\rho^\sharp\neq 0$.
\end{lem}
\begin{proof}
Let $\chi_\reg$ denote the character of the regular representation of $G$.
Then $\chi_\reg(1)=|G|$ and $\chi_\reg(g)=0$ if $g\neq 1$.
Moreover, we know that $\chi_\reg$ is uniform,
so
\[
f(1)
=\langle f,\chi_\reg\rangle_\bfG
=\langle f^\sharp,\chi_\reg\rangle_\bfG
=f^\sharp(1).
\]
If $\rho\in\cale(G)$,
then $\rho^\sharp(1)=\rho(1)\neq0$,
and hence $\rho^\sharp\neq 0$.
\end{proof}

The following lemma is well-known (\cf.~\cite{carter-finite} theorem~7.3.8).

\begin{lem}\label{0209}
Let $\rho$ be an irreducible character of $G$.
\begin{enumerate}
\item[(i)] If $\rho$ is unipotent,
then $\rho^\sharp\in\calv(G)_1$ where $\calv(G)_1$ denotes the span of unipotent characters of $G$.

\item[(ii)] If $\rho$ is not unipotent,
then $\langle\rho^\sharp,\eta\rangle_\bfG=0$ for any $\eta\in\cale(G)_1$.
\end{enumerate}
\end{lem}

\subsection{Uniform projection for a symplectic group}

\begin{prop}\label{0312}
Let $\rho_1,\rho_2$ be two irreducible unipotent characters of\/ $\Sp_{2n}(q)$.
If $\rho_1^\sharp=\rho_2^\sharp$,
then $\rho_1=\rho_2$.
\end{prop}
\begin{proof}
Let $\bfG=\Sp_{2n}$.
Suppose that $\rho_1=\rho_{\Lambda_1}$ and $\rho_2=\rho_{\Lambda_2}$ for some
$\Lambda_1,\Lambda_2\in\cals_\bfG$.
Since $\rho^\sharp_{\Lambda_1}=\rho^\sharp_{\Lambda_2}$,
by Proposition~\ref{0910} we have $\Lambda_1,\Lambda_2\in\cals_Z$
for some special symbol $Z$ of rank $n$ and defect $1$.
For any arrangement $\Phi$ of $Z_\rmI$ with a subset of pairs $\Psi$,
we have
\begin{align*}
\Biggl\langle\sum_{\Lambda\in C_{\Phi,\Psi}}\rho_\Lambda,\rho_{\Lambda_1}\Biggr\rangle_{\!\!\bfG}
=\Biggl\langle\sum_{\Lambda\in C_{\Phi,\Psi}}\rho_\Lambda,\rho_{\Lambda_1}^\sharp\Biggr\rangle_{\!\!\bfG}
=\Biggl\langle\sum_{\Lambda\in C_{\Phi,\Psi}}\rho_\Lambda,\rho_{\Lambda_2}^\sharp\Biggr\rangle_{\!\!\bfG}
=\Biggl\langle\sum_{\Lambda\in C_{\Phi,\Psi}}\rho_\Lambda,\rho_{\Lambda_2}\Biggr\rangle_{\!\!\bfG}
\end{align*}
because $\sum_{\Lambda\in C_{\Phi,\Psi}}\rho_\Lambda$ is uniform by Lemma~\ref{0226}.

Now if $\Lambda_1\neq\Lambda_2$,
by Lemma~\ref{0226} we can find $\Phi,\Psi$ such that $\Lambda_1\in C_{\Phi,\Psi}$ and $\Lambda_2\not\in C_{\Phi,\Psi}$,
and hence
\[
\Biggl\langle\sum_{\Lambda\in C_{\Phi,\Psi}}\rho_\Lambda,\rho_{\Lambda_1}\Biggr\rangle_{\!\!\bfG}=1
\quad\text{and}\quad
\Biggl\langle\sum_{\Lambda\in C_{\Phi,\Psi}}\rho_\Lambda,\rho_{\Lambda_2}\Biggr\rangle_{\!\!\bfG}=0.
\]
We get a contradiction.
\end{proof}

\begin{rem}
This proposition is known in \cite{DM-Lusztig} proposition 6.3, and \cite{GM-guide} theorem 4.4.23.
Here we provide a different proof.
\end{rem}

\subsection{Uniform projection for an orthogonal group}
For a nontrivial orthogonal group,
both ${\bf 1}_{\rmO_n^\epsilon}$ and $\sgn_{\rmO_n^\epsilon}$ are not uniform and
\[
{\bf 1}_{\rmO_n^\epsilon}^\sharp
=\tfrac{1}{2}({\bf 1}_{\rmO_n^\epsilon}+\sgn_{\rmO_n^\epsilon})
=\sgn_{\rmO_n^\epsilon}^\sharp
\]
by the following proposition:

\begin{prop}\label{0217}
Suppose $\rho$ is an irreducible character of\/ $\SO_n^\epsilon(q)$ such that
$\Ind_{\SO_n^\epsilon(q)}^{\rmO_n^\epsilon(q)}\rho$ is the direct sum $\rho_1+\rho_2$
of two irreducible characters $\rho_1,\rho_2$ of\/ $\rmO_n^\epsilon(q)$.
Then $\rho_1^{\sharp}=\rho_2^{\sharp}$.
\end{prop}
\begin{proof}
Fix an element $\sigma\in\rmO_n^\epsilon(q)\smallsetminus\SO_n^\epsilon(q)$.
Since $\SO_n^\epsilon(q)$ is a normal subgroup of index $2$ of $\rmO_n^\epsilon(q)$,
by the basic theory of induced characters (\cf.~\cite{sr-finite} proposition 20) we know that
\begin{align*}
\bigl[\Ind_{\SO_n^\epsilon(q)}^{\rmO_n^\epsilon(q)}\phi\bigr](y)
&=\frac{1}{|\SO_n^\epsilon(q)|}\sum_{t\in\rmO_n^\epsilon(q),\ t^{-1}yt\in\SO_n^\epsilon(q)}\phi(t^{-1}yt)\\
&=\begin{cases}
\phi(y)+\phi(\sigma^{-1} y\sigma), & \text{if $y\in\SO_n^\epsilon(q)$};\\
0, & \text{if $y\in\rmO_n^\epsilon(q)\smallsetminus\SO_n^\epsilon(q)$}
\end{cases}
\end{align*}
for any class function $\phi$ on $\SO_n^\epsilon(q)$.
In particular,
we have
$R_{\bfT,\theta}^{\rmO_n^\epsilon}(y)=0$ if $y\in\rmO_n^\epsilon(q)\smallsetminus\SO_n^\epsilon(q)$.
Now we know that $\rho_2=\rho_1\cdot\sgn$ and
$\rho_1|_{\SO_n^\epsilon(q)}=\rho_2|_{\SO_n^\epsilon(q)}=\eta.$
Thus
\begin{align*}
\langle\rho_1,R_{\bfT,\theta}^{\rmO_n^\epsilon}\rangle_{\rmO_n^\epsilon}
= \frac{1}{|\rmO_n^\epsilon(q)|}\sum_{y\in\rmO_n^\epsilon(q)}\rho_1(y)R_{\bfT,\theta}^{\rmO_n^\epsilon}(y^{-1})
&= \frac{1}{|\rmO_n^\epsilon(q)|}\sum_{y\in\SO_n^\epsilon(q)}\rho_1(y)R_{\bfT,\theta}^{\rmO_n^\epsilon}(y^{-1}) \\
&= \frac{1}{|\rmO_n^\epsilon(q)|}\sum_{y\in\SO_n^\epsilon(q)}\rho_2(y)R_{\bfT,\theta}^{\rmO_n^\epsilon}(y^{-1}) \\
&= \langle\rho_2,R_{\bfT,\theta}^{\rmO_n^\epsilon}\rangle_{\rmO_n^\epsilon}
\end{align*}
for any $(\bfT,\theta)$.
Because the subspace of uniform class function on $\rmO_n^\epsilon(q)$ is spanned the $R_{\bfT,\theta}^{\rmO_n^\epsilon}$'s,
the proposition is proved.
\end{proof}

\begin{prop}
Let $\rho_1,\rho_2$ be two irreducible unipotent characters of\/ $\rmO^\epsilon_{2n}(q)$ for $n\geq 1$.
If $\rho_1^\sharp=\rho_2^\sharp$,
then $\rho_1=\rho_2$ or $\rho_1=\rho_2\cdot\sgn$.
\end{prop}
\begin{proof}
Let $\bfG=\rmO^\epsilon_{2n}$.
Suppose that $\rho_1=\rho_{\Lambda_1}$ and $\rho_2=\rho_{\Lambda_2}$ for some $\Lambda_1,\Lambda_2\in\cals_\bfG$.
By the similar argument in the proof of Proposition~\ref{0312},
we know that $\Lambda_1,\Lambda_2\in\cals^\epsilon_Z$ for some special symbol $Z$ of rank $n$ and defect $0$.
Moreover,
we have
\[
\Biggl\langle\sum_{\Lambda\in C_{\Phi,\Psi}}\rho_\Lambda,\rho_{\Lambda_1}\Biggr\rangle_{\!\!\bfG}
=\Biggl\langle\sum_{\Lambda\in C_{\Phi,\Psi}}\rho_\Lambda,\rho_{\Lambda_2}\Biggr\rangle_{\!\!\bfG}
\]
for any arrangement $\Phi$ of $Z_\rmI$ with a subset of pairs $\Psi$.

If $\Lambda_1\neq\Lambda_2,\Lambda_2^\rmt$,
then by Lemma~\ref{0231} we can find $\Phi,\Psi$ such that $\Lambda_1\in C_{\Phi,\Psi}$ and $\Lambda_2\not\in C_{\Phi,\Psi}$,
and we get a contradiction.
Therefore we must have $\Lambda_1=\Lambda_2$ or $\Lambda_1=\Lambda_2^\rmt$.
If $\Lambda_1=\Lambda_2$, then $\rho_1=\rho_2$;
if $\Lambda_1=\Lambda_2^\rmt$, then $\rho_1=\rho_2\cdot\sgn$.
\end{proof}

\subsection{Uniform projection of $\chi_\bfG$}\label{0310}
Recall that in \cite{pan-odd} subsection 3.3, for $w\in\bfW_\bfG$ we define a character $\theta_w$ of $T_w$ of order two.
Then we show that the class function
\begin{equation}\label{0320}
\chi_\bfG=\frac{1}{|\bfW_\bfG|}\sum_{w\in\bfW_\bfG}R^\bfG_{\bfT_w,\theta_w}
\end{equation}
is in fact a linear character if $\bfG$ is a general linear, unitary or special orthogonal group,
in particular, $\chi_\bfG$ is uniform for these groups.
Note that $\chi_\bfG$ is the linear character of order $2$ corresponding to the element $-1$ in the center of $G^*$
(\cf.~\cite{DM} proposition 13.30).
We know that $\chi_\bfG|_{T_w}=\theta_w$,
and hence
\begin{equation}\label{0208}
R_{\bfT_w,\theta\theta_w}^\bfG=\chi_\bfG R_{\bfT_w,\theta}^\bfG
\end{equation}
for any $(\bfT_w,\theta)$.
Note that $\chi_\bfG$ is not a linear character if $\bfG$ is a symplectic group.

We know that $\Ind_{\SO_n^\epsilon(q)}^{\rmO_n^\epsilon(q)}\chi_{\SO_n^\epsilon}$
is the sum of two characters which differ by the $\sgn$ character.
We will denote these two characters by $\chi_{\rmO_n^\epsilon}$ and $\chi_{\rmO_n^\epsilon}\cdot\sgn$.
Hence, we have
\begin{align*}
\chi_{\rmO_n^\epsilon}+\chi_{\rmO_n^\epsilon}\cdot\sgn
=\Ind_{\SO_n^\epsilon(q)}^{\rmO_n^\epsilon(q)}
\Biggl[\frac{1}{|\bfW_{\SO_n^\epsilon}|}\sum_{w\in\bfW_{\SO_n^\epsilon}}R^{\SO_n^\epsilon}_{\bfT_w,\theta_w}\Biggr]
&=\frac{1}{|\bfW_{\SO_n^\epsilon}|}\sum_{w\in\bfW_{\SO_n^\epsilon}}
\Ind_{\SO_n^\epsilon(q)}^{\rmO_n^\epsilon(q)}R^{\SO_n^\epsilon}_{\bfT_w,\theta_w}\\
&=\frac{1}{|\bfW_{\SO_n^\epsilon}|}\sum_{w\in\bfW_{\SO_n^\epsilon}}R^{\rmO_n^\epsilon}_{\bfT_w,\theta_w}.
\end{align*}
The linear character $\chi_{\rmO_n^\epsilon}$ is not uniform,
in fact by Proposition~\ref{0217} we have
\[
\chi_{\rmO_n^\epsilon}^\sharp
=\frac{1}{2}(\chi_{\rmO_n^\epsilon}+\chi_{\rmO_n^\epsilon}\cdot\sgn)
=(\chi_{\rmO_n^\epsilon}\cdot\sgn)^\sharp.
\]
Note that (\ref{0208}) also holds for $\bfG=\rmO^\epsilon_n$, more precisely,
we have
\begin{align*}
\chi_{\rmO^\epsilon_n} R_{\bfT_w,\theta}^{\rmO^\epsilon_n}
=\chi_{\rmO^\epsilon_n} \Ind^{\rmO^\epsilon_n(q)}_{\SO^\epsilon_n(q)} R_{\bfT_w,\theta}^{\SO^\epsilon_n}
=\Ind^{\rmO^\epsilon_n(q)}_{\SO^\epsilon_n(q)}\bigl[\chi_{\SO^\epsilon_n} R_{\bfT_w,\theta}^{\SO^\epsilon_n}\bigr]
&=\Ind^{\rmO^\epsilon_n(q)}_{\SO^\epsilon_n(q)}R_{\bfT_w,\theta\theta_w}^{\SO^\epsilon_n} \\
&=R_{\bfT_w,\theta\theta_w}^{\rmO^\epsilon_n}.
\end{align*}

\begin{lem}\label{0303}
If $f$ is a uniform class function on $\rmO_n^\epsilon(q)$,
then $f\chi_{\rmO_n^\epsilon}$ is also uniform.
\end{lem}
\begin{proof}
Suppose that $f=\sum_{w\in\bfW_{\SO_n^\epsilon}}\sum_{\theta\in\cale(T_w)}a_{w,\theta}R_{\bfT_w,\theta}^{\rmO_n^\epsilon}$
for some coefficients $a_{w,\theta}$.
Then
\begin{align*}
f\chi_{\rmO_n^\epsilon}
=\sum_{w\in\bfW_{\SO_n^\epsilon}}\sum_{\theta\in\cale(T_w)}a_{w,\theta}\chi_{\rmO_n^\epsilon}
R_{\bfT_v,\theta}^{\rmO_n^\epsilon}
=\sum_{w\in\bfW_{\SO_n^\epsilon}}\sum_{\theta\in\cale(T_w)}a_{w,\theta}
R_{\bfT_w,\theta\theta_w}^{\rmO_n^\epsilon}.
\end{align*}
\end{proof}

\begin{lem}\label{0301}
Let $\rho$ be an irreducible character of\/ $\rmO_n^\epsilon(q)$.
Then
\[
(\rho\chi_{\rmO_n^\epsilon})^\sharp
=\rho^\sharp\chi_{\rmO_n^\epsilon}.
\]
\end{lem}
\begin{proof}
For any $w\in\bfW_{\SO_n^\epsilon}$ and $\theta\in\cale(T_w)$,
we have
\[
\langle\rho\chi_{\rmO_n^\epsilon},R^{\rmO_n^\epsilon}_{\bfT_w,\theta}\rangle_{\rmO_n^\epsilon}
=\frac{1}{|\rmO_n^\epsilon(q)|}\sum_{y\in \rmO_n^\epsilon(q)}\rho(y)\chi_{\rmO_n^\epsilon}(y)
R^{\rmO_n^\epsilon}_{\bfT_w,\theta}(y^{-1}).
\]
Because $\chi_{\rmO_n^\epsilon}(y)=\chi_{\rmO_n^\epsilon}(y^{-1})$ and
$\chi_{\rmO_n^\epsilon}R_{\bfT_w,\theta}^{\rmO_n^\epsilon}=R_{\bfT_w,\theta\theta_w}^{\rmO_n^\epsilon}$,
we have
\begin{align*}
\langle\rho\chi_{\rmO_n^\epsilon},R^{\rmO_n^\epsilon}_{\bfT_w,\theta}\rangle_{\rmO_n^\epsilon}
=\langle\rho,R^{\rmO_n^\epsilon}_{\bfT_w,\theta\theta_w}\rangle_{\rmO_n^\epsilon}
=\langle\rho^\sharp,R^{\rmO_n^\epsilon}_{\bfT_w,\theta\theta_w}\rangle_{\rmO_n^\epsilon}
=\langle\rho^\sharp\chi_{\rmO_n^\epsilon},R^{\rmO_n^\epsilon}_{\bfT_w,\theta}\rangle_{\rmO_n^\epsilon}.
\end{align*}
This means that $\rho\chi_{\rmO_n^\epsilon}$ and $\rho^\sharp\chi_{\rmO_n^\epsilon}$ have the same
uniform projection.
Now by Lemma~\ref{0303} we know that $\rho^\sharp\chi_{\rmO_n^\epsilon}$ is already uniform,
so the lemma is proved.
\end{proof}

\subsection{Uniform projection of $R_\bfG$}\label{0309}
Let $R_\bfG$ be the character of $G\times G$ given by
\[
R_{\bfG}=\sum_{\rho\in\cale(G)}\rho\otimes\rho
\]
with unipotent part
\[
R_{\bfG,1}=\sum_{\rho\in\cale(G)_1}\rho\otimes\rho.
\]
Hence we have the uniform projection $R_{\bfG}^\sharp=\sum_{\rho\in\cale(G)}\rho^\sharp\otimes\rho^\sharp$.
From \cite{amr} proposition~1.9,
when $G$ is connected,
we have
\[
R_{\bfG}^\sharp
=\frac{1}{|\bfW_\bfG|}\sum_{w\in\bfW_\bfG}\sum_{\theta\in\cale(T_w)}
R_{\bfT_w,\theta}^\bfG\otimes R_{\bfT_w,\theta}^\bfG.
\]
Then the unipotent part of $R_{\bfG}^\sharp$ is
\begin{equation}\label{0305}
R_{\bfG,1}^\sharp=\frac{1}{|\bfW_\bfG|}\sum_{w\in\bfW_\bfG}R_{\bfT_w,\bf1}^\bfG\otimes R_{\bfT_w,\bf1}^\bfG.
\end{equation}

\subsubsection{For a symplectic group}
Let $\bfG=\Sp_{2n}$ and let $Z$ be a special symbol of rank $n$ and defect $1$.
\[
R_{\bfG,1}=\sum_{\rho\in\cale(G)_1}\rho\otimes\rho
=\sum_Z R_{\bfG,Z}
\]
where $R_{\bfG,Z}:=\sum_{\Lambda\in\cals_Z}\rho_\Lambda\otimes\rho_\Lambda$
and $Z$ runs over special symbols of rank $n$ and defect $1$.

\begin{lem}\label{0313}
Let $\bfG=\Sp_{2n}$, and let $Z$ be a special symbol of rank $n$ and defect $1$.
Then
\[
\Biggl(\sum_{\Lambda\in\cals_Z}\rho_\Lambda\otimes\rho_\Lambda\Biggr)^\sharp
=\sum_{\Sigma\in\cals_{Z,1}}R^\bfG_\Sigma\otimes R^\bfG_\Sigma.
\]
\end{lem}
\begin{proof}
For $\Lambda\in\cals_Z$, we know that
\[
\rho_\Lambda^\sharp
=\frac{1}{2^{\deg(Z)}}\sum_{\Sigma\in\cals_{Z,1}}(-1)^{\langle\Sigma,\Lambda\rangle}R^\bfG_\Sigma
\]
by Proposition~\ref{0901}.
Therefore, we can write
\[
\sum_{\Lambda\in\cals_Z}\rho^\sharp_\Lambda\otimes\rho^\sharp_{\Lambda}
=\sum_{\Sigma\in\cals_{Z,1}}\sum_{\Sigma'\in\cals_{Z,1}}c_{\Sigma,\Sigma'}R^\bfG_{\Sigma}\otimes R^\bfG_{\Sigma'}
\]
where
\[
c_{\Sigma,\Sigma'}
=\frac{1}{2^{2\deg(Z)}}\sum_{\Lambda\in\cals_Z}(-1)^{\langle\Sigma,\Lambda\rangle+\langle\Sigma',\Lambda\rangle}.
\]
Hence by \cite{pan-uniform} lemma~4.4,
we have
\[
c_{\Sigma,\Sigma'}=\begin{cases}
1, & \text{if $\Sigma'=\Sigma$};\\
0, & \text{otherwise}.
\end{cases}
\]
Then the lemma is proved.
\end{proof}

Therefore, by Lemma~\ref{0313},
we have
\[
R_{\bfG,1}^\sharp
=\sum_Z R^\sharp_{\bfG,Z}
=\sum_Z\sum_{\Sigma\in\cals_{Z,1}} R^\bfG_{\Sigma}\otimes R^\bfG_{\Sigma}
=\sum_{\Sigma\in\cals_{n,1}}R^\bfG_{\Sigma}\otimes R^\bfG_{\Sigma}.
\]

\begin{lem}\label{0314}
Let $\bfG=\Sp_{2n}$, and let $Z$ be a special symbol of rank $n$ and defect $1$.
Suppose that $\Omega_Z$ is a unipotent character of $G\times G$ such that $\Omega_Z^\sharp=R^\sharp_{\bfG,Z}$.
Then $\Omega_Z=R_{\bfG,Z}$.
\end{lem}
\begin{proof}
From the assumption and Lemma~\ref{0313},
we have
\[
\Omega_Z^\sharp
=R_{\bfG,Z}^\sharp
=\Biggl(\sum_{\Lambda\in\cals_Z}\rho_\Lambda\otimes\rho_\Lambda\Biggr)^\sharp
=\sum_{\Sigma\in\cals_{Z,1}}R^\bfG_\Sigma\otimes R^\bfG_\Sigma.
\]
By Proposition~\ref{0901},
we can write
\[
\Omega_Z=\sum_{\Lambda,\Lambda'\in\cals_Z}a_{\Lambda,\Lambda'}\rho_\Lambda\otimes\rho_{\Lambda'}
\]
where $a_{\Lambda,\Lambda'}=\langle\rho_\Lambda\otimes\rho_{\Lambda'},\Omega_Z\rangle$.
By (i) of Lemma~\ref{0226},
the class function
\[
\sum_{\Lambda\in C_{\Phi,\Psi}}\sum_{\Lambda'\in C_{\Phi',\Psi'}}\rho_\Lambda\otimes\rho_{\Lambda'}
\]
on $G\times G$ is uniform for any arrangements $\Phi,\Phi'$ of $Z_\rmI$ with subsets of pairs
$\Psi,\Psi'$ respectively.
Then we have
\begin{align*}
\sum_{\Lambda\in C_{\Phi,\Psi}}\sum_{\Lambda'\in C_{\Phi',\Psi'}}a_{\Lambda,\Lambda'}
&=\Biggl\langle \sum_{\Lambda\in C_{\Phi,\Psi}}\sum_{\Lambda'\in C_{\Phi',\Psi'}}\rho_\Lambda\otimes\rho_{\Lambda'},\Omega_Z\Biggr\rangle_{\!\!\bfG} \\
&=\Biggl\langle \sum_{\Lambda\in C_{\Phi,\Psi}}\sum_{\Lambda'\in C_{\Phi',\Psi'}}\rho_\Lambda\otimes\rho_{\Lambda'},\Omega_Z^\sharp\Biggr\rangle_{\!\!\bfG} \\
&=\Biggl\langle \sum_{\Lambda\in C_{\Phi,\Psi}}\sum_{\Lambda'\in C_{\Phi',\Psi'}}\rho_\Lambda\otimes\rho_{\Lambda'}, \sum_{\Lambda\in\cals_Z}\rho_{\Lambda''}^\sharp\otimes\rho_{\Lambda''}^\sharp\Biggr\rangle_{\!\!\bfG} \\
&=\Biggl\langle \sum_{\Lambda\in C_{\Phi,\Psi}}\sum_{\Lambda'\in C_{\Phi',\Psi'}}\rho_\Lambda\otimes\rho_{\Lambda'}, \sum_{\Lambda\in\cals_Z}\rho_{\Lambda''}\otimes\rho_{\Lambda''}\Biggr\rangle_{\!\!\bfG}
\end{align*}
by Lemma~\ref{0313}.
For a symbol $\Lambda''\in\cals_Z$ to contribute a multiplicity,
it is required that $\Lambda''=\Lambda$ for some $\Lambda\in C_{\Phi,\Psi}$ and
$\Lambda''=\Lambda'$ for some $\Lambda'\in C_{\Phi',\Psi'}$, i.e.,
$\Lambda''\in C_{\Phi,\Psi}\cap C_{\Phi',\Psi'}$.
Therefore, we have
\begin{equation}\label{0319}
\sum_{\Lambda\in C_{\Phi,\Psi}}\sum_{\Lambda'\in C_{\Phi',\Psi'}}a_{\Lambda,\Lambda'}
=|C_{\Phi,\Psi}\cap C_{\Phi',\Psi'}|
\end{equation}

Now suppose that $\Lambda'\neq\Lambda$,
by (ii) of Lemma~\ref{0226}, we can find two arrangements $\Phi,\Phi'$ with subsets of pairs $\Psi,\Psi'$
respectively such that $\Lambda\in C_{\Phi,\Psi}$, $\Lambda'\in C_{\Phi',\Psi'}$ and
$|C_{\Phi,\Psi}\cap C_{\Phi',\Psi'}|=\emptyset$.
Because each $a_{\Lambda,\Lambda'}$ is a non-negative integer,
we conclude that $a_{\Lambda,\Lambda'}=0$ by (\ref{0319}).
Moreover, for any $\Lambda\in\cals_Z$, by (iii) of Lemma~\ref{0226},
(\ref{0319}) can be reduced to
$a_{\Lambda,\Lambda}=1$.
\end{proof}

\begin{prop}\label{0321}
Let $\bfG=\Sp_{2n}$.
If\/ $\Omega$ is a unipotent character of $G\times G$ such that
$\Omega^\sharp=R_{\bfG,1}^\sharp$,
then $\Omega=R_{\bfG,1}$.
\end{prop}
\begin{proof}
Note that $\calv(G)_1=\bigoplus_Z\calv(G)_Z$ and hence
\[
\calv(G)_1\otimes\calv(G)_1=\bigoplus_{Z,Z'}\calv(G)_Z\otimes\calv(G)_{Z'}
\]
where $Z,Z'$ run over special symbols of rank $n$ and defect $1$.
Then, as an element in $\calv(G)_1\otimes\calv(G)_1$,
write $\Omega=\sum_{Z,Z'}\Omega_{Z,Z'}$ where $\Omega_{Z,Z'}$ is the projection of $\Omega$
over the subspace $\calv(G)_Z\otimes\calv(G)_{Z'}$.
Then
\[
\sum_Z R_{\bfG,Z}^\sharp
=R_{\bfG,1}^\sharp
=\Omega^\sharp=\sum_{Z,Z'}\Omega_{Z,Z'}^\sharp.
\]
Therefore, for $Z\neq Z'$ we have $\Omega_{Z,Z'}^\sharp=0$ and hence $\Omega_{Z,Z'}=0$ by Lemma~\ref{0207},
moreover, we have $\Omega_{Z,Z}^\sharp=R^\sharp_{\bfG,Z}$ and hence $\Omega_{Z,Z}=R_{\bfG,Z}$ by Lemma~\ref{0314}.
Therefore
\[
\Omega=\sum_Z\Omega_{Z,Z}=\sum_Z R_{\bfG,Z}=R_{\bfG,1}.
\]
\end{proof}

\subsubsection{For an even orthogonal group}
Let $\bfG=\rmO^\epsilon_{2n}$.
From \cite{amr} proposition~1.15,
we have
\begin{equation}\label{0306}
R_{\rmO_{2n}^\epsilon,1}^\sharp
=\frac{1}{|W_n|}\sum_{w\in W_n^\epsilon}
R_{\bfT_w,\bf1}^{\rmO_{2n}^\epsilon}\otimes R_{\bfT_w,\bf1}^{\rmO_{2n}^\epsilon}.
\end{equation}
Let $\bar\cals_{n,0}$ be a complete set of representatives of subsets $\{\Sigma,\Sigma^\rmt\}$ in $\cals_{n,0}$,
Similarly, let $\bar\cals_{\rmO^\epsilon_{2n}}$ be a complete set of representatives of subsets
$\{\Lambda,\Lambda^\rmt\}$ in $\cals_{\rmO^\epsilon_{2n}}$.

\begin{lem}\label{0316}
Let $\bfG=\rmO^\epsilon_{2n}$,
and let $Z$ be a special symbol of rank $n$ and defect $0$.
Then
\[
\Biggl(\sum_{\Lambda\in\cals^\epsilon_Z}\rho_\Lambda\otimes\rho_\Lambda\Biggr)^\sharp
=\frac{1}{2}\sum_{\Sigma\in\bar\cals_{Z,0}}R^\bfG_\Sigma\otimes R^\bfG_\Sigma.
\]
\end{lem}
\begin{proof}
For $\Lambda\in\cals^\epsilon_Z$, we know that
\[
\rho_\Lambda^\sharp
=\frac{1}{2^{\deg(Z)}}\sum_{\Sigma\in\bar\cals_{Z,0}}(-1)^{\langle\Sigma,\Lambda\rangle}R^\bfG_\Sigma
\]
by Proposition~\ref{0315}.
Therefore, we can write
\[
\sum_{\Lambda\in\cals^\epsilon_Z}\rho^\sharp_\Lambda\otimes\rho^\sharp_{\Lambda}
=\sum_{\Sigma\in\bar\cals_{Z,0}}\sum_{\Sigma'\in\bar\cals_{Z,0}}
c_{\Sigma,\Sigma'}R^\bfG_{\Sigma}\otimes R^\bfG_{\Sigma'}
\]
where
\[
c_{\Sigma,\Sigma'}
=\frac{1}{2^{2\deg(Z)}}\sum_{\Lambda\in\cals^\epsilon_Z}(-1)^{\langle\Sigma,\Lambda\rangle+\langle\Sigma',\Lambda\rangle}.
\]
Hence by \cite{pan-uniform} lemma~4.10,
we have
\[
c_{\Sigma,\Sigma'}=\begin{cases}
\frac{1}{2}, & \text{if $\Sigma'=\Sigma$};\\
0, & \text{otherwise}.
\end{cases}
\]
Note that $\Sigma,\Sigma'$ are chosen from $\bar\cals_{Z,0}$,
so $\Sigma'\neq\Sigma^\rmt$.
Then the proposition is proved.
\end{proof}

Similar to the case of symplectic groups,
for $\bfG=\rmO^\epsilon_{2n}$,
by Lemma~\ref{0316},
we have
\[
R_{\bfG,1}^\sharp
=\sum_Z R_{\bfG,Z}^\sharp
=\sum_Z\Biggl(\sum_{\Lambda\in\cals_Z^\epsilon}\rho_\Lambda\otimes\rho_\Lambda\Biggr)^\sharp
=\frac{1}{2}\sum_Z\sum_{\Sigma\in\bar\cals_{Z,0}}R_\Sigma^\bfG\otimes R_\Sigma^\bfG
=\frac{1}{2}\sum_{\Sigma\in\bar\cals_{n,0}}R^\bfG_{\Sigma}\otimes R^\bfG_{\Sigma}.
\]

\begin{lem}\label{0318}
Let $\bfG=\rmO^\epsilon_{2n}$, and let $Z$ be a special symbol of rank $n$ and defect $0$.
Suppose that $\Omega_Z$ is a unipotent character of $G\times G$ such that
$\Omega_Z^\sharp=\frac{1}{2}\sum_{\Sigma\in\bar\cals_{Z,0}}R^\bfG_\Sigma\otimes R^\bfG_\Sigma$.
Then
\[
\Omega_Z=\sum_{\Lambda\in\bar\cals^\epsilon_Z}
\bigl[a_{\Lambda,\Lambda}\rho_\Lambda\otimes\rho_\Lambda
+a_{\Lambda^\rmt,\Lambda}\rho_{\Lambda^\rmt}\otimes\rho_\Lambda
+a_{\Lambda,\Lambda^\rmt}\rho_\Lambda\otimes\rho_{\Lambda^\rmt}
+a_{\Lambda^\rmt,\Lambda^\rmt}\rho_{\Lambda^\rmt}\otimes\rho_{\Lambda^\rmt}\bigr]
\]
where each $a_{\Lambda,\Lambda'}$ is a non-negative integer and
$a_{\Lambda,\Lambda}+a_{\Lambda^\rmt,\Lambda}+a_{\Lambda,\Lambda^\rmt}+a_{\Lambda^\rmt,\Lambda^\rmt}=2$.
\end{lem}
\begin{proof}
Write $\Omega_Z=\sum_{\Lambda,\Lambda'\in\cals_Z^\epsilon}a_{\Lambda,\Lambda'}\rho_\Lambda\otimes\rho_{\Lambda'}$
where $a_{\Lambda,\Lambda'}=\langle\rho_\Lambda\otimes\rho_\Lambda',\Omega_Z\rangle_\bfG$.
By the same argument in the proof of Lemma~\ref{0314},
for any arrangements $\Phi,\Phi'$ of $Z_\rmI$ with admissible subsets of pairs $\Psi,\Psi'$ respectively,
we have
\begin{align*}
\sum_{\Lambda\in C_{\Phi,\Psi}}\sum_{\Lambda'\in C_{\Phi',\Psi'}}a_{\Lambda,\Lambda'}
&=\Biggl\langle \sum_{\Lambda\in C_{\Phi,\Psi}}\sum_{\Lambda'\in C_{\Phi',\Psi'}}\rho_\Lambda\otimes\rho_{\Lambda'}, \sum_{\Lambda''\in\cals_Z}\rho_{\Lambda''}\otimes\rho_{\Lambda''}\Biggr\rangle_{\!\!\bfG}
\end{align*}
by Lemma~\ref{0316}.
Therefore, we have
\begin{equation}\label{0317}
\sum_{\Lambda\in C_{\Phi,\Psi}}\sum_{\Lambda'\in C_{\Phi',\Psi'}}a_{\Lambda,\Lambda'}
=|C_{\Phi,\Psi}\cap C_{\Phi',\Psi'}|
\end{equation}
for any arrangements $\Phi,\Phi'$ of $Z_\rmI$ with admissible subsets of pairs $\Psi,\Psi'$ respectively.

Now suppose that $\Lambda'\neq\Lambda,\Lambda^\rmt$,
by Lemma~\ref{0224}, we can find two arrangements $\Phi,\Phi'$ with admissible subset of pairs
$\Psi,\Psi'$ respectively such that $\Lambda,\Lambda^\rmt\in C_{\Phi,\Psi}$ and
$\Lambda',\Lambda'^\rmt\in C_{\Phi',\Psi'}$ and $|C_{\Phi,\Psi}\cap C_{\Phi',\Psi'}|=\emptyset$.
Because each $a_{\Lambda,\Lambda'}$ is non-negative,
we conclude that
$a_{\Lambda,\Lambda'}=a_{\Lambda^\rmt,\Lambda'}=a_{\Lambda,\Lambda'^\rmt}=a_{\Lambda^\rmt,\Lambda'^\rmt}=0$
by (\ref{0317}).
Moreover, for any $\Lambda\in\cals^\epsilon_Z$, by Lemma~\ref{0225},
(\ref{0317}) can be reduced to
$a_{\Lambda,\Lambda}+a_{\Lambda^\rmt,\Lambda}+a_{\Lambda,\Lambda^\rmt}+a_{\Lambda^\rmt,\Lambda^\rmt}=2$.
\end{proof}

\begin{prop}\label{0311}
Let $\bfG=\rmO^\epsilon_{2n}$.
If\/ $\Omega$ is a unipotent character of $G\times G$ such that
$\Omega^\sharp=R_{\bfG,1}^\sharp$,
then
\[
\Omega=\sum_{\Lambda\in\bar\cals_{\rmO^\epsilon_{2n}}}
\bigl[a_{\Lambda,\Lambda}\rho_\Lambda\otimes\rho_\Lambda
+a_{\Lambda^\rmt,\Lambda}\rho_{\Lambda^\rmt}\otimes\rho_\Lambda
+a_{\Lambda,\Lambda^\rmt}\rho_\Lambda\otimes\rho_{\Lambda^\rmt}
+a_{\Lambda^\rmt,\Lambda^\rmt}\rho_{\Lambda^\rmt}\otimes\rho_{\Lambda^\rmt}\bigr]
\]
where each $a_{\Lambda,\Lambda'}$ is a non-negative integer and
$a_{\Lambda,\Lambda}+a_{\Lambda^\rmt,\Lambda}+a_{\Lambda,\Lambda^\rmt}+a_{\Lambda^\rmt,\Lambda^\rmt}=2$.
\end{prop}
\begin{proof}
As in the proof of Proposition~\ref{0321}, we know that
$\Omega=\sum_Z\Omega_{Z,Z}$ where $Z$ runs over special symbols of rank $n$ and defect $0$,
and $\Omega_{Z,Z}$ is the projection of $\Omega$ over the subspace $\calv(G)_Z\otimes\calv(G)_Z$.
The the proposition follows from Lemma~\ref{0318} immediately.
\end{proof}


\section{Uniform Projection of the Weil Character}

\subsection{Lusztig correspondence}\label{0204}
Let $\bfG^*$ denote the dual group of $\bfG$.
The Frobenius endomorphism of $\bfG^*$ is still denoted by $F$,
and $G^*=\bfG^{*F}$ the group of rational points.
It is known that there is a bijection between the set of $G$-conjugacy classes of $(\bfT,\theta)$ where
$\theta\in\cale(T)$ and the set of $G^*$-conjugacy classes of $(\bfT^*,s)$ where
$\bfT^*$ is a rational maximal torus in $\bfG^*$ and $s\in T^*=\bfT^{*F}$.
If $(\bfT,\theta)$ is corresponding to $(\bfT^*,s)$,
then $R_{\bfT,\theta}$ is also denoted by $R_{\bfT^*,s}$

For a semisimple element $s\in (G^*)^0$,
define
\[
\cale(G)_s=\{\,\rho\in\cale(G)\mid\langle\rho,R_{\bfT^*,s}\rangle_\bfG\neq 0
\text{ for some }\bfT^*\text{ containing }s\,\}.
\]
The set $\cale(G)_s$ is called a \emph{Lusztig series},
and it is known that $\cale(G)$ is partitioned into Lusztig series indexed by the conjugacy classes $(s)$
of semisimple elements $s$ in $(G^*)^0$,
i.e.,
\begin{equation}\label{0219}
\cale(G)=\bigcup_{(s)\subset (G^*)^0}\cale(G)_s.
\end{equation}
An irreducible character of $G$ is unipotent if it is in $\cale(G)_1$.
For any semisimple element $s\in G^*$,
let $\calv(G)_s$ denote the linear span of the set $\cale(G)_s$.
From (\ref{0219}), we have an orthogonal decomposition
\[
\calv(G)=\bigoplus_{(s)\subset (G^*)^0}\calv(G)_s.
\]

The following result (\cf.~\cite{DM} theorem 13.23, remark 13.24) is fundamental for
the classification of $\cale(G)$:

\begin{prop}[Lusztig]\label{0201}
There is a bijection
\[
\grL_s\colon\cale(G)_s\longrightarrow \cale(C_{G^*}(s))_1
\]
satisfying the condition
\begin{equation}\label{0218}
\langle\rho,\epsilon_\bfG R^\bfG_{\bfT^*,s}\rangle_G
=\langle\grL_s(\rho),\epsilon_{C_{\bfG^*}(s)}R^{C_{\bfG^*}(s)}_{\bfT^*,{\bf1}}\rangle_{C_{G^*}(s)}
\end{equation}
for any rational maximal torus $\bfT^*$ containing $s$
where $C_{G^*}(s)$ denotes the centralizer of $s$ in\/ $G^*$.
Moreover, we have
\begin{equation}\label{0203}
\dim\rho=\frac{|G|_{p'}}{|C_{G^*}(s)|_{p'}}\dim\grL_s(\rho)
\end{equation}
where $|G|_{p'}$ denotes greatest factor of\/ $|G|$ not divided by $p$,
$\epsilon_G=(-1)^{r}$ where $r$ is the $\bfF_q$-rank of\/ $G$.
\end{prop}
Such a bijection $\grL_s$ is called a \emph{Lusztig correspondence}.
Note that $\grL_s$ is usually not uniquely determined.
A discuss on the ambiguity of $\grL_s$ can be found in \cite{pan-ambiguity}.
The following proposition follows from \cite{lusztig-book} (9.9.1) (see also \cite{pan-chain01} proposition 3.3):

\begin{prop}[Lusztig]\label{0202}
The Lusztig series $\cale(G)_s$ has a cuspidal representation
if and only if $\cale(C_{G^*}(s))_1$ has a cuspidal representation
and the largest $\bfF_q$-split torus in the center of $C_{G^*}(s)$ coincides with
the largest $\bfF_q$-split torus in the center of\/ $G^*$.
In this case, $\rho$ is cuspidal if and only if $\grL_s(\rho)$ is cuspidal.
\end{prop}

From (\ref{0218}) we see that the bijection $\grL_s$ can be extended by linearity to be an isometry
\[
\grL_s\colon\calv(G)_s\longrightarrow\calv(C_{G^*}(s))_1,
\]
still denoted by $\grL_s$.
Moreover, by (\ref{0218}) the isometry maps $\calv(G)_s^\sharp$ onto $\calv(C_{G^*}(s))_1^\sharp$,
and we have
\begin{equation}\label{0410}
\grL_s(\rho^\sharp)=\grL_s(\rho)^\sharp
\end{equation}
for any $\rho\in\cale(G)_s$.

Suppose $G$ is special orthogonal or orthogonal.
Note that $\chi_\bfG R^\bfG_{\bfT^*,s}=R^\bfG_{\bfT^*,-s}$ from (\ref{0208}).
Hence, if $\rho$ is in $\cale(G)_s$ for some $s$,
then $\rho\chi_\bfG$ is also an irreducible character and is in $\cale(G)_{-s}$.
It is clear that $\epsilon_{C_{G^*}(s)}=\epsilon_{C_{G^*}(-s)}$,
and (\ref{0218}) implies that
\[
\grL_s(R_{\bfT^*,s}^\bfG)
=\epsilon_\bfG\epsilon_{C_{\bfG^*}(s)}R_{\bfT^*,1}^{C_{\bfG^*}(s)}
=\grL_{-s}(R_{\bfT^*,-s}^\bfG)
=\grL_{-s}(R_{\bfT^*,s}^\bfG\chi_\bfG).
\]
Therefore the composition $\grL_{-s}^{-1}\circ\grL_s$ gives a bijection $\cale(G)_s\rightarrow\cale(G)_{-s}$
by $\rho\mapsto\rho\chi_\bfG$.
Clearly, $\rho$ is cuspidal if and only if $\rho\chi_\bfG$ is cuspidal.

\subsection{Uniform projection of the Weil character}
Recall that we fix a nontrivial additive character $\psi$ of $\bfF_q$.
Other additive characters of $\bfF_q$ are the form $\psi_a$ for $a\in\bfF_q^\times$,
where $\psi_a(x):=\psi(ax)$.
The Weil character $\omega^\psi_{\Sp_{2n}}$ really depends on $\psi$.
It is known that $\omega_{\Sp_{2n}}^\psi=\omega_{\Sp_{2n}}^{\psi_a}$
if $a$ is a square in $\bfF_q^\times$,
and $\omega_{\Sp_{2n}}^\psi\neq\omega_{\Sp_{2n}}^{\psi_a}$ otherwise.
We shall see that the uniform projection $(\omega_{\Sp_{2n}}^\psi)^\sharp$ does not depend on
the character $\psi$.

It is known that $\omega_{\Sp_{2n}}^\psi$ is the sum
$\chi^{(1)}_{\frac{q^n+1}{2}}+\chi^{(1)}_{\frac{q^n-1}{2}}$ of two irreducible characters of degree
$\frac{q^n+1}{2}$ and $\frac{q^n-1}{2}$ respectively.
Similarly,
$\omega_{\Sp_{2n}}^{\psi_a}=\chi^{(2)}_{\frac{q^n+1}{2}}+\chi^{(2)}_{\frac{q^n-1}{2}}$.

\begin{lem}\label{0411}
The irreducible character $\chi^{(i)}_{\frac{q^n\pm1}{2}}$ of\/ $\Sp_{2n}(q)$ is not unipotent
for $i=1,2$.
\end{lem}
\begin{proof}
Suppose that $\chi^{(i)}_{\frac{q^n\pm1}{2}}$ is unipotent.
Then $\chi^{(i)}_{\frac{q^n\pm1}{2}}=\rho_\Lambda$ for some symbol $\Lambda$ of rank $n$
and defect $\equiv 1\pmod 4$.
Let $Z$ be the special symbol associated to $\Lambda$.
Now we know that
\begin{equation}\label{0911}
2^{\deg(Z)}\deg\rho_\Lambda\equiv q^{a_\Lambda}\pmod{q^{a_\Lambda+1}}
\end{equation}
by \cite{lg-symplectic} lemma 5.2.
So we must have $\deg(Z)=1$ and $a_\Lambda=0$.
But now $m\geq\deg(Z)=1$,
we have $a_\Lambda=a_Z>0$ by Lemma~\ref{0912},
and we get a contradiction.
Therefore $\chi^{(i)}_{\frac{q^n\pm1}{2}}$ is not unipotent.
\end{proof}

\begin{prop}\label{0304}
Let $a$ be a non-square element in $\bfF_q^\times$.
Then
\[
(\omega_{\Sp_{2n}}^\psi)^\sharp
=(\omega_{\Sp_{2n}}^{\psi_a})^\sharp
=\frac{1}{2}\omega_{\Sp_{2n}}^\psi+\frac{1}{2}\omega_{\Sp_{2n}}^{\psi_a}.
\]
\end{prop}
\begin{proof}
Let $\rho=\chi^{(1)}_{\frac{q^n+1}{2}}$ or $\chi^{(2)}_{\frac{q^n+1}{2}}$.
Then $\rho$ is in the Lusztig series $\cale(\Sp_{2n}(q))_{s^+}$ for some semisimple element $s^+$ in $\SO_{2n+1}(q)$.
Hence by (\ref{0203})
\[
\frac{q^n+1}{2}=\frac{|\Sp_{2n}(q)|_{p'}}{|C_{\SO_{2n+1}(q)}(s^+)|_{p'}}\deg\grL_{s^+}(\rho).
\]
We know that $\rho$ is not unipotent by Lemma~\ref{0411},
so $s^+$ is not in the center of $\SO_{2n+1}(q)$ and hence
$C_{\SO_{2n+1}(q)}(s^+)\neq\SO_{2n+1}(q)$.
Because $\deg\rho$ is so small,
the only possibility is that $C_{\SO_{2n+1}(q)}(s^+)\simeq \rmO_{2n}^+(q)$
and $\deg\grL_{s^+}(\rho)=1$.
Moreover,
we know that the only degree one irreducible unipotent characters of $\rmO_{2n}^+(q)$ are
${\bf 1}_{\rmO_{2n}^+}$ and $\sgn_{\rmO_{2n}^+}$.
Hence we obtain a bijection
\begin{equation}\label{0914}
\grL_{s^+}\colon\Bigl\{\chi^{(1)}_{\frac{q^n+1}{2}},\chi^{(2)}_{\frac{q^n+1}{2}}\Bigr\}
\longrightarrow\{{\bf1}_{\rmO_{2n}^+},\sgn_{\rmO_{2n}^+}\}.
\end{equation}
Because ${\bf1}_{\rmO_{2n}^+}+\sgn_{\rmO_{2n}^+}=\Ind^{\rmO^+_{2n}(q)}_{\SO^+_{2n}(q)}{\bf1}_{\SO_{2n}^+}$ and the trivial character ${\bf1}_{\SO_{2n}^+}$ of $\SO_{2n}^+(q)$
is uniform,
the class function ${\bf1}_{\rmO_{2n}^+}+\sgn_{\rmO_{2n}^+}$ of $\rmO^+_{2n}(q)$ is then also uniform,
and from Proposition~\ref{0217},
we know that
\[
{\bf 1}_{\rmO_{2n}^+}^\sharp=\sgn^\sharp_{\rmO_{2n}^+}
=\frac{1}{2}({\bf 1}_{\rmO_{2n}^+}+\sgn_{\rmO_{2n}^+}).
\]
Hence by (\ref{0410}) the class function $\chi^{(1)}_{\frac{q^n+1}{2}}+\chi^{(2)}_{\frac{q^n+1}{2}}$ of $\Sp_{2n}(q)$ is uniform,
and
\[
\Bigl(\chi^{(1)}_{\frac{q^n+1}{2}}\Bigr)^\sharp
=\Bigl(\chi^{(2)}_{\frac{q^n+1}{2}}\Bigr)^\sharp
=\frac{1}{2}\Bigl(\chi^{(1)}_{\frac{q^n+1}{2}}+\chi^{(2)}_{\frac{q^n+1}{2}}\Bigr).
\]

Similarly,
$\chi^{(1)}_{\frac{q^n-1}{2}}$ and $\chi^{(2)}_{\frac{q^n-1}{2}}$ are in the Lusztig series $\cale(\Sp_{2n}(q))_{s^-}$
where the centralizer $C_{\SO_{2n+1}(q)}(s^-)\simeq\rmO_{2n}^-(q)$,
and hence the class function $\chi^{(1)}_{\frac{q^n-1}{2}}+\chi^{(2)}_{\frac{q^n-1}{2}}$ of $\Sp_{2n}(q)$ is uniform.
Moreover,
\[
\Bigl(\chi^{(1)}_{\frac{q^n-1}{2}}\Bigr)^\sharp
=\Bigl(\chi^{(2)}_{\frac{q^n-1}{2}}\Bigr)^\sharp
=\frac{1}{2}\Bigl(\chi^{(1)}_{\frac{q^n-1}{2}}+\chi^{(2)}_{\frac{q^n-1}{2}}\Bigr).
\]
Therefore, we have
\begin{align*}
(\omega_{\Sp_{2n}}^\psi)^\sharp
=\Bigl(\chi^{(1)}_{\frac{q^n+1}{2}}+\chi^{(1)}_{\frac{q^n-1}{2}}\Bigr)^\sharp
& =\frac{1}{2}\Bigl(\chi^{(1)}_{\frac{q^n+1}{2}}+\chi^{(2)}_{\frac{q^n+1}{2}}\Bigr)
+\frac{1}{2}\Bigl(\chi^{(1)}_{\frac{q^n-1}{2}}+\chi^{(2)}_{\frac{q^n-1}{2}}\Bigr) \\
& =\frac{1}{2}\omega_{\Sp_{2n}}^\psi+\frac{1}{2}\omega_{\Sp_{2n}}^{\psi_a}
=(\omega_{\Sp_{2n}}^{\psi_a})^\sharp.
\end{align*}
\end{proof}

The following result is originally proved in \cite{pan-odd} via some tedious computation.
Now we have an easier proof via the Lusztig correspondence.

\begin{cor}
For $n\geq 1$,
we have the decomposition
\[
\omega_{\Sp_{2n}}^\sharp
=\frac{1}{|\bfW_{\Sp_{2n}}|}\sum_{w\in\bfW_{\Sp_{2n}}}\epsilon_w R_{\bfT_w,\theta_w}^{\Sp_{2n}}.
\]
\end{cor}
\begin{proof}
Identify $\bfW_{\Sp_{2n}}=W_n=W_n^+\cup W_n^-$.
Then by \cite{carter-finite} corollary 7.6.5, we have
\[
{\bf 1}_{\SO_{2n}^\epsilon}
=\frac{1}{|W_n^\epsilon|}\sum_{w\in\bfW_n^\epsilon}R_{\bfT_w,\bf1}^{\SO_{2n}^\epsilon}.
\]
Therefore
\[
{\bf 1}_{\rmO_{2n}^\epsilon}+\sgn_{\rmO_{2n}^\epsilon}
=\Ind^{\rmO_{2n}^\epsilon(q)}_{\SO_{2n}^\epsilon(q)}{\bf 1}_{\SO_{2n}^\epsilon}
=\frac{1}{|W_n^\epsilon|}\sum_{w\in W_n^\epsilon}R_{\bfT_w,\bf1}^{\rmO_{2n}^\epsilon}.
\]
Now $\epsilon_{\Sp_{2n}}\epsilon_{\rmO_{2n}^+}=1$, then by Proposition~\ref{0201} and (\ref{0914}),
we have
\[
\chi^{(1)}_{\frac{q^n+1}{2}}+\chi^{(2)}_{\frac{q^n+1}{2}}
=\frac{1}{|W_n^+|}\sum_{w\in W_n^+}R_{\bfT_w,\theta_w}^{\Sp_{2n}}.
\]
Because $\epsilon_{\Sp_{2n}}\epsilon_{\rmO_{2n}^-}=-1$,
we have
\[
\chi^{(1)}_{\frac{q^n-1}{2}}+\chi^{(2)}_{\frac{q^n-1}{2}}
=\frac{1}{|W_n^-|}\sum_{w\in W_n^-}-R_{\bfT_w,\theta_w}^{\Sp_{2n}}.
\]
Since $\epsilon_w=\epsilon$ if $w\in W_n^\epsilon$ and $|W_n|=2|W_n^\epsilon|$,
the corollary follows form Proposition~\ref{0304} immediately.
\end{proof}


\section{Compatibility for Unitary Groups}
In this section,
we consider a dual pair of two unitary groups over a finite field of odd characteristic.

\subsection{Unipotent characters of a unitary group}\label{0520}
We first recall some notations from \cite{FS}.
The \emph{$\beta$-set} $X$ associated to a partition $\lambda=[\lambda_1,\ldots,\lambda_m]$ is defined to be
\[
X=\{\lambda_1+(m-1),\lambda_2+(m-2),\ldots,\lambda_{m-1}+1,\lambda_m\},
\]
i.e., $\lambda=\Upsilon(X)$ where $\Upsilon$ is given in Subsection~\ref{0229}.
A \emph{$2$-hook} of a $\beta$-set $X$ is a pair $(y,x)$ of non-negative integers such that $y\not\in X$, $x\in X$ and $y+2=x$.
If $(y,x)$ is a $2$-hook of $X$,
then the $\beta$-set $X_1:=\{y\}\cup(X\smallsetminus\{x\})$ is said to be obtained from $X$ by
removing the $2$-hook $(y,x)$.
It is easy to see that $\|\Upsilon(X_1)\|=\|\Upsilon(X)\|-2$.
For a $\beta$-set $X$, we obtain a sequence of $\beta$-sets $X=X_0,X_1,\ldots,X_h$ by removing a $2$-hook successively,
until we can not remove any $2$-hook.
The resulting $\beta$-set is denoted by $X_\infty$.
The partition $\lambda_\infty:=\Upsilon(X_\infty)$ is called the \emph{$2$-core} of $\lambda$.
It is easy to see that $X_\infty$ is of the form
\begin{equation}\label{0513}
\{2k-1,2k-3,\ldots,3,1,2l-2,2l-4,\ldots,2,0\}
\end{equation}
for some non-negative integers $k,l$,
and $\lambda_\infty=[d,d-1,\ldots,1]$ where
\[
d=\begin{cases}
k-l, & \text{if $k\geq l$};\\
l-k-1, & \text{if $k<l$}.
\end{cases}
\]

Let $X$ be the $\beta$-set associate to a partition $\lambda$,
we define two $\beta$-sets
\begin{align*}
X(0) &=\biggl\{\,\frac{x}{2}\mid x\in X,\ x\equiv 0\pmod 2\,\biggr\},\\
X(1) &=\biggl\{\,\frac{x-1}{2}\mid x\in X,\ x\equiv 1\pmod 2\,\biggr\},
\end{align*}
and the symbol
\begin{equation}
\Lambda_\lambda=\begin{cases}
\binom{X(0)}{X(1)}, & \text{if $\ell(\lambda)$ is odd};\\
\binom{X(1)}{X(0)}, & \text{if $\ell(\lambda)$ is even}
\end{cases}
\end{equation}
where $\ell(\lambda)$ is given in Subsection~\ref{0229}. 
It is easy to check that (\cf.~\cite{FS} p.223)
\begin{equation}\label{0516}
\|\lambda\|=\|\lambda_\infty\|+2\|\Upsilon(\Lambda_\lambda)\|.
\end{equation}

\begin{exam}
Consider the partition $\lambda=[7,5,4,3,2,2]$ of $23$.
Then $X=\{12,9,7,5,3,2\}$, $X(0)=\{6,1\}$, $X(1)=\{4,3,2,1\}$,
$X_\infty=\{7,5,3,2,1,0\}$, and hence $\lambda_\infty=[2,1]$.
Now $\ell(\lambda)=6$ is even.
Therefore $\Lambda_\lambda=\binom{4,3,2,1}{6,1}$ and $\Upsilon(\Lambda_\lambda)=\sqbinom{1,1,1,1}{5,1}$.
\end{exam}

\begin{lem}\label{0512}
Let $\lambda$ be a partition.
Then
\[
{\rm def}(\Lambda_\lambda)=\begin{cases}
\ell(\lambda_\infty), & \text{if both $\ell(\lambda),\ell(\lambda_\infty)$ are even};\\
\ell(\lambda_\infty)+1, & \text{if $\ell(\lambda)$ is odd and $\ell(\lambda_\infty)$ is even};\\
-\ell(\lambda_\infty), & \text{if both $\ell(\lambda),\ell(\lambda_\infty)$ are odd};\\
-\ell(\lambda_\infty)-1, & \text{if $\ell(\lambda)$ is even and $\ell(\lambda_\infty)$ is odd}.
\end{cases}
\]
\end{lem}
\begin{proof}
After removing all possible $2$-hooks,
we obtain a $\beta$-set $X_\infty$ of the form (\ref{0513}) for some non-negative integers $k,l$.
Moreover, we know that $\ell(\lambda)=|X|=|X_\infty|=k+l$, $|X(1)|=k$ and $|X(0)|=l$.
Now we have the following situations:
\begin{enumerate}
\item If $k+l$ is even and $k\geq l$,
then ${\rm def}(\Lambda_\lambda)=k-l$ and $\ell(\lambda_\infty)=k-l$.
Now $\ell(\lambda_\infty)$ is even and ${\rm def}(\Lambda_\lambda)=\ell(\lambda_\infty)$.

\item If $k+l$ is even and $k<l$,
then ${\rm def}(\Lambda_\lambda)=k-l$ and $\ell(\lambda_\infty)=l-k-1$.
Now $\ell(\lambda_\infty)$ is odd and ${\rm def}(\Lambda_\lambda)=-\ell(\lambda_\infty)-1$.

\item If $k+l$ is odd and $k\geq l$,
then ${\rm def}(\Lambda_\lambda)=l-k$ and $\ell(\lambda_\infty)=k-l$.
Now $\ell(\lambda_\infty)$ is odd and ${\rm def}(\Lambda_\lambda)=-\ell(\lambda_\infty)$.

\item If $k+l$ is odd and $k<l$,
then ${\rm def}(\Lambda_\lambda)=l-k$ and $\ell(\lambda_\infty)=l-k-1$.
Now $\ell(\lambda_\infty)$ is even and ${\rm def}(\Lambda_\lambda)=\ell(\lambda_\infty)+1$.
\end{enumerate}
\end{proof}

Now given a non-negative integer $d$ and a bi-partition $\sqbinom{\mu}{\nu}$ such that
$n=\frac{d(d+1)}{2}+2\|\sqbinom{\mu}{\nu}\|$,
we can associate a unique partition $\lambda$ of $n$ as follows.
Write $\mu=[\mu_1,\ldots,\mu_r]$ and $\nu=[\nu_1,\ldots,\nu_s]$.
\begin{enumerate}
\item If $d$ is even and $r-s>d$,
we define
\begin{align*}
X(0) &=\{\mu_1+(r-1),\ldots,\mu_{r-1}+1,\mu_r\}, \\
X(1) &=\{\nu_1+(r-d-2),\ldots,\nu_s+(r-s-d-1),r-s-d-2,\ldots,1,0\}.
\end{align*}

\item If $d$ is even and $r-s\leq d$,
we define
\begin{align*}
X(0) &=\{\nu_1+(s-1),\ldots,\nu_{s-1}+1,\nu_s\}, \\
X(1) &=\{\mu_1+(s+d-1),\ldots,\mu_r+(s-r+d),s-r+d-1,\ldots,1,0\}.
\end{align*}

\item If $d$ is odd and $r-s\geq -d$,
we define
\begin{align*}
X(0) &=\{\mu_1+(r-1),\ldots,\mu_{r-1}+1,\mu_r\}, \\
X(1) &=\{\nu_1+(r+d-1),\ldots,\nu_s+(r-s+d),r-s-d-1,\ldots,1,0\}.
\end{align*}

\item If $d$ is odd and $r-s<-d$,
we define
\begin{align*}
X(0) &=\{\nu_1+(s-1),\ldots,\nu_{s-1}+1,\nu_s\}, \\
X(1) &=\{\mu_1+(s-d-2),\ldots,\mu_r+(s-r-d-1),s-r-d-2,\ldots,1,0\}.
\end{align*}
\end{enumerate}
Let
\[
X=\{\,2x\mid x\in X(0)\,\}\cup\{\,2x'+1\mid x'\in X(1)\,\}\quad\text{ and }\quad
\lambda=\Upsilon(X).
\]
It is not difficult to check that $\|\lambda\|=n$,
$\lambda_\infty=[d,d-1,\ldots,1]$ and
$\Upsilon(\Lambda_\lambda)=\sqbinom{\mu}{\nu}$.
So we establish a bijection
\[
\textstyle
\calp(n)\longleftrightarrow
\Bigl\{\,\bigl(d,\sqbinom{\mu}{\nu}\bigr)\mid d\in\bbZ_{\geq 0},
\ \frac{d(d+1)}{2}+2\|\sqbinom{\mu}{\nu}\|=n\,\Bigr\}.
\]

\begin{exam}
Suppose that $d=2$ and the bi-partition is $\sqbinom{1,1,1,1}{5,1}$.
Now $d$ is even, $r=4$, $s=2$ and $r-s\leq d$.
Then $X(1)=\{4,3,2,1\}$, $X(0)=\{6,1\}$, $X=\{12,9,7,5,3,2\}$ and $\lambda=[7,5,4,3,2,2]$.
\end{exam}

\subsection{Howe correspondence for a dual pair of unitary groups}
Define $\calb^+_{\rmU,\rmU}$ to be the set consisting of pairs of symbols $(\Lambda,\Lambda')\in\calb^+$ (\cf.~Subsection~\ref{0230})
such that
\begin{equation}\label{0517}
{\rm def}(\Lambda')=\begin{cases}
-{\rm def}(\Lambda),-{\rm def}(\Lambda)+1, & \text{if ${\rm def}(\Lambda)$ is even};\\
-{\rm def}(\Lambda)+1,-{\rm def}(\Lambda)+2, & \text{if ${\rm def}(\Lambda)$ is odd},
\end{cases}
\end{equation}
Similarly, define $\calb^-_{\rmU,\rmU}$ to be the set consisting of pairs of symbols $(\Lambda,\Lambda')\in\calb^-$ such that
\begin{equation}\label{0518}
{\rm def}(\Lambda')=\begin{cases}
-{\rm def}(\Lambda)-2,-{\rm def}(\Lambda)-1, & \text{if ${\rm def}(\Lambda)$ is even};\\
-{\rm def}(\Lambda)-1,-{\rm def}(\Lambda), & \text{if ${\rm def}(\Lambda)$ is odd},
\end{cases}
\end{equation}
For $\epsilon=+,-$,
it is clear that the relation $\calb^\epsilon_{\rmU,\rmU}$ is symmetric, i.e.,
$(\Lambda,\Lambda')\in\calb^\epsilon_{\rmU,\rmU}$ if and only if $(\Lambda',\Lambda)\in\calb^\epsilon_{\rmU,\rmU}$.
Finally, we define
\begin{equation}\label{0521}
\calb_{\rmU_n,\rmU_{n'}}=\{(\Lambda_\lambda,\Lambda_{\lambda'})\in\calb^+_{\rmU,\rmU}\cup\calb^-_{\rmU,\rmU}\mid\|\lambda\|=n,\ \|\lambda'\|=n'\,\}.
\end{equation}

\begin{lem}\label{0522}
Let $\lambda,\lambda'$ be two partitions.
\begin{enumerate}
\item[(i)] Suppose that $(\Lambda_\lambda,\Lambda_{\lambda'})\in\calb^+_{\rmU,\rmU}$.
Then
\[
\ell(\lambda'_\infty)=\begin{cases}
\ell(\lambda_\infty), & \text{if $\ell(\lambda_\infty)=0$};\\
\ell(\lambda_\infty)-1, & \text{if $\ell(\lambda_\infty)$ is even and nonzero};\\
\ell(\lambda_\infty)+1, & \text{if $\ell(\lambda_\infty)$ is odd}.
\end{cases}
\]

\item[(ii)] Suppose that $(\Lambda_\lambda,\Lambda_{\lambda'})\in\calb^-_{\rmU,\rmU}$.
Then
\[
\ell(\lambda'_\infty)=\begin{cases}
\ell(\lambda_\infty)+1, & \text{if $\ell(\lambda_\infty)$ is even};\\
\ell(\lambda_\infty)-1, & \text{if $\ell(\lambda_\infty)$ is odd}.
\end{cases}
\]
\end{enumerate}
\end{lem}
\begin{proof}
Suppose that $(\Lambda_\lambda,\Lambda_{\lambda'})\in\calb^+_{\rmU,\rmU}$.
We have the following situations:
\begin{enumerate}
\item If $\ell(\lambda_\infty)=0$,
then ${\rm def}(\Lambda_\lambda)=0,1$ by Lemma~\ref{0512},
and hence ${\rm def}(\Lambda_{\lambda'})=1,0$ by (\ref{0517}).
Therefore $\ell(\lambda'_\infty)=0$ by Lemma~\ref{0512}, again.

\item If $\ell(\lambda_\infty)=2m$ for $m\geq 1$,
then ${\rm def}(\Lambda_\lambda)=2m,2m+1$ by Lemma~\ref{0512}
and hence ${\rm def}(\Lambda_{\lambda'})=-2m+1,-2m$ by (\ref{0517}).
Therefore $\ell(\lambda'_\infty)=2m-1=\ell(\lambda_\infty)-1$ by Lemma~\ref{0512}, again.

\item If $\ell(\lambda_\infty)=2m-1$ for $m\geq 1$,
then ${\rm def}(\Lambda_\lambda)=-2m+1,-2m$ by Lemma~\ref{0512}
and hence ${\rm def}(\Lambda_{\lambda'})=2m,2m+1$ by (\ref{0517}).
Therefore $\ell(\lambda'_\infty)=2m=\ell(\lambda_\infty)+1$ by Lemma~\ref{0512}, again.
\end{enumerate}

Suppose that $(\Lambda_\lambda,\Lambda_{\lambda'})\in\calb^-_{\rmU,\rmU}$.
We have the following situations:
\begin{enumerate}
\item[(4)] If $\ell(\lambda_\infty)=2m$ for some $m\geq 0$,
then ${\rm def}(\Lambda_\lambda)=2m,2m+1$ by Lemma~\ref{0512}
and hence ${\rm def}(\Lambda_{\lambda'})=-2m-1,-2m-2$ by (\ref{0518}).
Therefore $\ell(\lambda'_\infty)=2m+1=\ell(\lambda_\infty)+1$ by Lemma~\ref{0512}.

\item[(5)] Suppose that $\ell(\lambda_\infty)=2m-1$ for some $m\geq 1$,
then ${\rm def}(\Lambda_\lambda)=-2m+1,-2m$ by Lemma~\ref{0512}
and hence ${\rm def}(\Lambda_{\lambda'})=2m-1,2m-2$ by (\ref{0518}).
Therefore $\ell(\lambda'_\infty)=2m-2=\ell(\lambda_\infty)-1$ by Lemma~\ref{0512}.
\end{enumerate}
\end{proof}

\begin{lem}\label{0523}
We have
\[
\calb_{\rmU_n,\rmU_{n'}}\subset
\begin{cases}
\calb^+_{\rmU,\rmU}, & \text{if $n+n'$ is even};\\
\calb^-_{\rmU,\rmU}, & \text{if $n+n'$ is odd}.
\end{cases}
\]
\end{lem}
\begin{proof}
Suppose that $\lambda\in\calp(n)$, $\lambda'\in\calp(n')$ and $(\Lambda_\lambda,\Lambda_{\lambda'})\in\calb^+_{\rmU_n,\rmU_{n'}}$.
Note that we have $\|\lambda\|\equiv\|\lambda_\infty\|\pmod 2$ and $\|\lambda'\|\equiv\|\lambda'_\infty\|\pmod 2$.
\begin{enumerate}
\item If both $n$ and $n'$ are even,
then $\ell(\lambda_\infty)\equiv 0,3\pmod 4$ and $\ell(\lambda'_\infty)\equiv 0,3\pmod 4$.
We have $(\Lambda_\lambda,\Lambda_{\lambda'})\in\calb^+_{\rmU,\rmU}$ by Lemma~\ref{0522}.

\item If both $n$ and $n'$ are odd,
the $\ell(\lambda_\infty)\equiv 1,2\pmod 4$ and $\ell(\lambda'_\infty)\equiv 1,2\pmod 4$.
We have $(\Lambda_\lambda,\Lambda_{\lambda'})\in\calb^+_{\rmU,\rmU}$ by Lemma~\ref{0522}.

\item If $n$ is even and $n'$ is odd,
then $\ell(\lambda_\infty)\equiv 0,3\pmod 4$ and $\ell(\lambda'_\infty)\equiv 1,2\pmod 4$.
We have $(\Lambda_\lambda,\Lambda_{\lambda'})\in\calb^-_{\rmU,\rmU}$ by Lemma~\ref{0522}.

\item If $n$ is odd and $n'$ is even,
then $\ell(\lambda_\infty)\equiv 1,2\pmod 4$ and $\ell(\lambda'_\infty)\equiv 0,3\pmod 4$.
We see that $(\Lambda_\lambda,\Lambda_{\lambda'})\in\calb^-_{\rmU,\rmU}$ by Lemma~\ref{0522}.
\end{enumerate}
\end{proof}

It is known that the irreducible unipotent characters of $\rmU_n(q)$
are parametrized by the set $\calp(n)$ (\cf.~\cite{FS} p.223),
so we denote $\rho=\rho_\lambda$ for $\rho\in\cale(\rmU_n(q))_1$ and $\lambda\in\calp(n)$.
It is well known that $\rho_\lambda$ is cuspidal if and only if
$\lambda=\lambda_\infty$, i.e.,
$\rho_\lambda$ is cuspidal if and only if
\[
\Lambda_\lambda=\begin{cases}
\binom{d-1,d-2,\ldots,1,0}{-}, & \text{if $d$ is even};\\
\binom{-}{d-1,d-2,\ldots,1,0}, & \text{if $d$ is odd}
\end{cases}
\]
for some non-negative integer $d$.

\begin{rem}
In the terminology of \cite{amr} proposition 5.14, for $\lambda\in\calp(n)$,
one associates $\rho_\lambda\in\cale(\rmU_n(q))_1$ a bi-partition $\lambda(0)\boxtimes\lambda(1)$
where $\lambda(0),\lambda(1)$ are the \emph{$2$-quotients of parameter $1$} of $\lambda$.
In our notation, write $\Upsilon(\Lambda_\lambda)=\sqbinom{\mu}{\nu}$.
Then $\mu=\Upsilon(X(0)),\nu=\Upsilon(X(1))$ are the $2$-quotient of parameter $\ell(\lambda)$ if $\ell(\lambda)$ is odd;
and $\mu=\Upsilon(X(1)),\nu=\Upsilon(X(0))$ are the $2$-quotient of parameter $\ell(\lambda)$ if $\ell(\lambda)$ is even.
Since the $2$-quotient depends only on parameter $t\pmod 2$,
we see that two notations are consistent, i.e., $\lambda(0)=\mu$ and $\lambda(1)=\nu$.
\end{rem}

The following proposition on the Howe correspondence of unipotent characters for $(\rmU_n,\rmU_{n'})$
is rephrased from \cite{amr} th\'eor\`eme 5.15.
This form is consistent with Proposition~\ref{0220}.

\begin{prop}[Aubert-Michel-Rouquier]\label{0519}
Let $(\bfG,\bfG')=(\rmU_n,\rmU_{n'})$.
We have the decomposition
\[
\omega_{\bfG,\bfG',1}
=\sum_{(\Lambda_\lambda,\Lambda_{\lambda'})\in\calb_{\bfG,\bfG'}}\rho_\lambda\otimes\rho_{\lambda'}.
\]
\end{prop}
\begin{proof}
Let $\lambda,\lambda'$ be two partitions of $n,n'$ respectively.
Write $\ell(\lambda_\infty)=d$, $\ell(\lambda'_\infty)=d'$,
$\Upsilon(\Lambda_\lambda)=\sqbinom{\mu}{\nu}$, and $\Upsilon(\Lambda'_\lambda)=\sqbinom{\mu'}{\nu'}$.
Now by \cite{amr} proposition 5,6, we know that $\rho_\lambda\otimes\rho_{\lambda'}$ occurs in $\omega_{\rmU_n,\rmU_{n'},1}$
if and only if $\lambda,\lambda'$ are \emph{$2$-transverse} (\cf.~\cite{amr} p.389).

Suppose that $\lambda,\lambda'$ are $2$-transverse.
Then by \cite{amr} lemme 5.9, we know that $|d-d'|=1$ or $d=d'=0$.
Without loss of generality, we may assume that $d'=d+1$ or $d=d'=0$.
From the proof of \cite{amr} th\'eor\`eme 5.15, we have:
\begin{enumerate}
\item If either $d$ is odd or $(d,d')=(0,0)$,
then we have $\lambda(1)\preccurlyeq\lambda'(0)$ and $\lambda'(1)\preccurlyeq\lambda(0)$,
i.e., $\nu\preccurlyeq\mu'$ and $\nu'\preccurlyeq\mu$,
which implies that $(\Lambda_{\lambda},\Lambda_{\lambda'})\in\calb^+_{\rmU,\rmU}$.

\item If $d$ is even and $(d,d')\neq(0,0)$,
then we have $\lambda(0)\preccurlyeq\lambda'(1)$ and $\lambda'(0)\preccurlyeq\lambda(1)$,
i.e., $\mu\preccurlyeq\nu'$ and $\mu'\preccurlyeq\nu$,
which implies that $(\Lambda_{\lambda},\Lambda_{\lambda'})\in\calb^-_{\rmU,\rmU}$.
\end{enumerate}

Conversely, suppose that $(\Lambda_{\lambda},\Lambda_{\lambda'})\in\calb_{\rmU_n,\rmU_{n'}}$.
\begin{enumerate}
\item If $n+n'$ is even,
then by Lemmas~\ref{0523} and \ref{0522},
we have either $(d,d')=(0,0)$ or $|d-d'|=1$ and $\nu\preccurlyeq\mu'$ and
$\nu'\preccurlyeq\mu$.

\item If $n+n'$ is odd,
then by Lemmas~\ref{0523} and \ref{0522},
we have $|d-d'|=1$ and $\nu\preccurlyeq\mu'$ and $\nu'\preccurlyeq\mu$.
\end{enumerate}
Then we see that $\Upsilon(\Lambda_\lambda)$ and $\Upsilon(\Lambda_{\lambda'})$ are $2$-transverse by
\cite{amr} proposition 5.12 and the proof of \cite{amr} th\'eor\`eme 5.15.
\end{proof}

\subsection{Lusztig correspondence for unitary groups}\label{0509}
Let $\bfG$ be a unitary group.
Then, for $s\in G^*$, we know that
$C_{\bfG^*}(s)=\prod_{\langle\lambda\rangle}\bfG_{[\lambda]}(s)$ where each $\bfG_{[\lambda]}(s)$ 
is a unitary group or a general linear group.
Define $\bfG^{(0)}=\prod_{\langle\lambda\rangle,\ \lambda\neq 1}\bfG_{[\lambda]}(s)$ 
and $\bfG^{(1)}(s)=\bfG_{[1]}(s)$
where $\bfG_{[\lambda]}(s)$ is given in Subsection~\ref{0501}.
Then $C_{\bfG^*}(s)=\bfG^{(0)}\times\bfG^{(1)}$ and hence there exists a one-to-one
correspondence
\begin{equation}\label{0507}
\Xi_s\colon \cale(G)_s\longrightarrow\cale(G^{(0)}\times G^{(1)})_1
\end{equation}
by Proposition~\ref{0201}.
So if $\rho\in\cale(G)_s$,
then $\Xi_s(\rho)$ can be written as $\rho^{(0)}\otimes\rho^{(1)}$
where $\rho^{(j)}$ is an irreducible unipotent character of $G^{(j)}$ for $j=0,1$.
Then $s$ can be written as
\begin{equation}\label{0524}
s=s^{(0)}\times s^{(1)}
\end{equation}
where $s^{(1)}$ is the part whose eigenvalues are all equal to $1$,
and $s^{(0)}$ is the part whose eigenvalues do not contain $1$.
In particular, $s^{(j)}$ is in the center of $\bfG^{(j)}$ for $j=0,1$.

The following can be extracted from \cite{amr} th\'eor\`eme 2.6 (\cf.~\cite{pan-chain01} theorem 3.10):

\begin{prop}\label{0510}
Let $(\bfG,\bfG')=(\rmU_n,\rmU_{n'})$.
Let $\rho\in\cale(G)_s$ and $\rho'\in\cale(G')_{s'}$ for some $s,s'$.
Then $\rho\otimes\rho'$ occurs in $\omega^\psi_{\bfG,\bfG'}$ if and only if the following conditions hold:
\begin{itemize}
\item $s^{(0)}=s'^{(0)}$, $G^{(0)}\simeq G'^{(0)}$ and $\rho^{(0)}=\rho'^{(0)}$,

\item $\rho^{(1)}\otimes\rho'^{(1)}$ occurs in $\omega_{\bfG^{(1)},\bfG'^{(1)},1}$.
\end{itemize}
\end{prop}


\section{Compatibility for Even Orthogonal Groups}
In this section we consider the dual pair $(\bfG,\bfG')=(\Sp_{2n},\rmO^\epsilon_{2n'})$
over a finite field of odd characteristic.

\subsection{Decomposition of the Weil character}
The following result is modified from the main theorem of \cite{srinivasan}.
See also \cite{amr} proposition 2.1.
Note that from \cite{pan-uniform} theorem 3.10,
we do not assume the cardinality of the finite field to be large.

\begin{prop}\label{0307}
Consider the dual pair $(\Sp_{2n},\rmO^\epsilon_{2n'})$ over a finite field of odd characteristic.
We have the decomposition
\begin{itemize}
\item if $n'>n$,
\begin{multline*}
\omega_{\Sp_{2n},\rmO^\epsilon_{2n'}}^\sharp
=\frac{1}{2}\sum_{k=0}^n \frac{1}{|W_k|}\frac{1}{|\bfW_{\Sp_{2(n-k)}}|}\frac{1}{|\bfW_{\SO_{2(n'-k)}}|}\sum_{v\in\bfW_k}
\sum_{\theta\in\cale(T_v)} \\
\sum_{w\in\bfW_{\Sp_{2(n-k)}}}\sum_{w'\in\bfW_{\SO_{2(n'-k)}^{\epsilon_v\epsilon}}}
\epsilon_v R^{\Sp_{2n}}_{\bfT_v\times\bfT_w,\theta\otimes\bf 1}\otimes
R^{\rmO^{\epsilon}_{2n'}}_{\bfT_v\times\bfT_{w'},\theta\otimes\bf 1};
\end{multline*}
\item if $n'\leq n$,
\begin{multline*}
\omega_{\Sp_{2n},\rmO^\epsilon_{2n'}}^\sharp
=\frac{1}{2}\sum_{k=0}^{n'-1}\frac{1}{|W_k|}\frac{1}{|\bfW_{\Sp_{2(n-k)}}|}\frac{1}{|\bfW_{\SO_{2(n'-k)}}|}\sum_{v\in W_k}
\sum_{\theta\in\cale(T_v)} \\
\sum_{w\in\bfW_{\Sp_{2(n-k)}}}\sum_{w'\in\bfW_{\SO_{2(n'-k)}^{\epsilon_v\epsilon}}}
\epsilon_v R^{\Sp_{2n}}_{\bfT_v\times\bfT_w,\theta\otimes\bf 1}\otimes R^{\rmO^{\epsilon}_{2n'}}_{\bfT_v\times\bfT_{w'},\theta\times\bf 1}\\
+\frac{1}{2}\frac{1}{|W_{n'}^\epsilon|}\frac{1}{|\bfW_{\Sp_{2(n-n')}}|}
\sum_{v\in W_{n'}^\epsilon}\sum_{\theta\in\cale(T_v)}\sum_{w\in\bfW_{\Sp_{2(n-n')}}}
\epsilon R^{\Sp_{2n}}_{\bfT_v\times\bfT_w,\theta\otimes\bf 1}\otimes R^{\rmO^\epsilon_{2n'}}_{\bfT_v,\theta}.
\end{multline*}
\end{itemize}
\end{prop}
\begin{proof}
The proof is the same as that of \cite{amr} proposition 2.1.
\end{proof}

If $s$ is a semisimple element in a maximal torus $T^*$ of $\rmO_{2k}^{\epsilon'}(q)$ for some $\epsilon'$,
then the semisimple element $(s,1)$ in $\rmO_{2k}^{\epsilon'}(q)\times\SO_{2(n-k)+1}(q)\subset\SO_{2n+1}(q)$
is denoted by $s_n$ for $n\geq k$;
similarly, $s_{n'}=(s,1)$ is a semisimple element in
$\rmO_{2k}^{\epsilon'}(q)\times\rmO_{2(n'-k)}^{\epsilon''}(q)\subseteq\rmO_{2n'}^\epsilon(q)$ for $n'\geq k$ with
$\epsilon'\epsilon''=\epsilon$.
From Subsection~\ref{0204}, we know that $R_{\bfT,\theta}^\bfG=R^\bfG_{\bfT^*,s}$
where $(\bfT^*,s)$ is corresponding to $(\bfT,\theta)$.
Hence we can rewrite the decompositions in Proposition~\ref{0307} as follows:
\begin{itemize}
\item if $n'>n$,
\begin{multline*}
\omega_{\Sp_{2n},\rmO^\epsilon_{2n'}}^\sharp
=\frac{1}{2}\sum_{k=0}^n \frac{1}{|W_k|}\frac{1}{|\bfW_{\Sp_{2(n-k)}}|}\frac{1}{|\bfW_{\SO_{2(n'-k)}}|}\sum_{v\in W_k}
\sum_{s\in T_v^*} \\
\sum_{w\in\bfW_{\Sp_{2(n-k)}}}\sum_{w'\in\bfW_{\SO_{2(n'-k)}^{\epsilon_v\epsilon}}}
\epsilon_v R^{\Sp_{2n}}_{\bfT_v^*\times\bfT_w^*,s_n}\otimes
R^{\rmO^{\epsilon}_{2n'}}_{\bfT_v^*\times\bfT_{w'}^*,s_{n'}};
\end{multline*}
\item if $n'\leq n$,
\begin{multline*}
\omega_{\Sp_{2n},\rmO^\epsilon_{2n'}}^\sharp
=\frac{1}{2}\sum_{k=0}^{n'-1}\frac{1}{|W_k|}\frac{1}{|\bfW_{\Sp_{2(n-k)}}|}\frac{1}{|\bfW_{\SO_{2(n'-k)}}|}\sum_{v\in W_k}
\sum_{s\in T_v^*} \\
\sum_{w\in\bfW_{\Sp_{2(n-k)}}}\sum_{w'\in\bfW_{\SO_{2(n'-k)}^{\epsilon_v\epsilon}}}
\epsilon_v R^{\Sp_{2n}}_{\bfT_v^*\times\bfT_w^*,s_n}\otimes R^{\rmO^{\epsilon}_{2n'}}_{\bfT_v^*\times\bfT_{w'}^*,s_{n'}}\\
+\frac{1}{2}\frac{1}{|W_{n'}^\epsilon|}\frac{1}{|\bfW_{\Sp_{2(n-n')}}|}
\sum_{v\in W_{n'}^\epsilon}\sum_{s\in T_v^*}
\sum_{w\in\bfW_{\Sp_{2(n-n')}}}
\epsilon R^{\Sp_{2n}}_{\bfT_v^*\times\bfT_w^*,s_n}\otimes R^{\rmO^\epsilon_{2n'}}_{\bfT_v^*,s}.
\end{multline*}
\end{itemize}

For $k\leq\min(n,n')$, $v\in W_k$, and $l\leq k$,
we define
\begin{align*}
T_{v,l}^* &=\{\,s\in T_v^*\mid \nu_{1}(s)=k-l\,\},\\
X_{k,l} &=\begin{cases}
\bigcup_{v\in W_k}T^*_{v,l} & \text{if $k\leq n$ and $k<n'$;} \\
\bigcup_{v\in W_k^\epsilon}T^*_{v,l} & \text{if $k\leq n$ and $k=n'$}
\end{cases}
\end{align*}
where $\nu_1(s)$ is given in Subsection~\ref{0501}.
The set $X_{k,l}$ is a subset of $\rmO_{2k}^\pm(q)$ when $k\leq n$ and $k<n'$ 
(or $\rmO^\epsilon_{2k}(q)$ when $k=n'\leq n$)
with an action by $W_k$ (or $W_k^\epsilon$).
For $s\in X_{k,l}$, let $(s)$ denote the orbit of $s$ under the action.
If $s\in X_{l,l}$ where $l\leq\min(n,n')$,
let $\omega^\psi_{\Sp_{2n},\rmO^\epsilon_{2n'},s}$ denote the orthogonal projection of $\omega^\psi_{\Sp_{2n},\rmO_{2n'}^\epsilon}$
over $\calv(\Sp_{2n}(q))_{s_{n}}\otimes\calv(\rmO^\epsilon_{2n'}(q))_{s_{n'}}$.
Then we have (\cf.~\cite{amr} proposition~2.4)
\begin{equation}\label{0603}
\omega^\psi_{\Sp_{2n},\rmO_{2n'}^\epsilon}
=\sum_{l=0}^{\min(n,n')}\sum_{(s)\subset X_{l,l}}\omega^\psi_{\Sp_{2n},\rmO_{2n'}^\epsilon,s}.
\end{equation}

The following lemma is a reformulation of \cite{amr} proposition 2.5.
Here we consider the two Weil characters $\omega^\psi_{\Sp_{2n},\rmO_{2n'}^+}$ and
$\omega^\psi_{\Sp_{2n},\rmO_{2n'}^-}$ separately
but in \cite{amr} $\omega_{n,n'}$ denotes the formal sum $\omega_{\Sp_{2n},\rmO_{2n'}^+}+\omega_{\Sp_{2n},\rmO_{2n'}^-}$.

\begin{lem}\label{0504}
Suppose $s\in X_{l,l}$ for some $l\leq\min(n,n')$.
Then
\begin{itemize}
\item if $n'>n$,
\begin{multline*}
\omega_{\Sp_{2n},\rmO^\epsilon_{2n'},s}^\sharp
=\frac{1}{2}\sum_{k=l}^{n}\frac{1}{|W_k(s_k)|}\frac{1}{|\bfW_{\Sp_{2(n-k)}}|}\frac{1}{|\bfW_{\SO_{2(n'-k)}}|}
\sum_{v\in W_k(s_k)}\\
\sum_{w\in\bfW_{\Sp_{2(n-k)}}}\sum_{w'\in\bfW_{\SO_{2(n'-k)}^{\epsilon_v\epsilon}}}
\epsilon_v R^{\Sp_{2n}}_{\bfT_v\times\bfT_w,s_n}\otimes R^{\rmO_{2n'}^\epsilon}_{\bfT_v\times\bfT_{w'},s_{n'}};
\end{multline*}
\item if $n'\leq n$,
\begin{multline*}
\omega_{\Sp_{2n},\rmO^\epsilon_{2n'},s}^\sharp
=\frac{1}{2}\sum_{k=l}^{n'-1}\frac{1}{|W_k(s_k)|}\frac{1}{|\bfW_{\Sp_{2(n-k)}}|}\frac{1}{|\bfW_{\SO_{2(n'-k)}}|}
\sum_{v\in W_k(s_k)} \\
\sum_{w\in\bfW_{\Sp_{2(n-k)}}}\sum_{w'\in\bfW_{\SO_{2(n'-k)}^{\epsilon_v\epsilon}}}
\epsilon_v R^{\Sp_{2n}}_{\bfT_v\times\bfT_w,s_n}\otimes R^{\rmO^{\epsilon}_{2n'}}_{\bfT_v\times\bfT_{w'},s_{n'}}\\
+\frac{1}{2}\frac{1}{|W_{n'}^\epsilon|}\frac{1}{|\bfW_{\Sp_{2(n-n')}}|}
\sum_{v\in W_{n'}^\epsilon}
\sum_{w\in\bfW_{\Sp_{2(n-n')}}}
\epsilon R^{\Sp_{2n}}_{\bfT_v\times\bfT_w,s_n}\otimes R^{\rmO^\epsilon_{2n'}}_{\bfT_v,s_{n'}}
\end{multline*}
\end{itemize}
where $W_k(s_k)$ denotes the stabilizer of $s_k$ in $W_k$.
\end{lem}
\begin{proof}
First suppose that $n'>n$.
For $v\in W_k$,
because $T_v^*=\bigcup_{l=0}^kT_{v,l}^*$,
we have
\begin{multline*}
\omega_{\Sp_{2n},\rmO^\epsilon_{2n'}}^\sharp
=\frac{1}{2}\sum_{k=0}^n\sum_{l=0}^k\frac{1}{|W_k|}\frac{1}{|\bfW_{\Sp_{2(n-k)}}|}
\frac{1}{|\bfW_{\SO_{2(n'-k)}}|}\sum_{v\in W_k}
\sum_{s\in T_{v,l}^*} \\
\sum_{w\in\bfW_{\Sp_{2(n-k)}}}\sum_{w'\in\bfW_{\SO_{2(n'-k)}^{\epsilon_v\epsilon}}}
\epsilon_v R^{\Sp_{2n}}_{\bfT_v^*\times\bfT_w^*,s_n}\otimes
R^{\rmO^{\epsilon}_{2n'}}_{\bfT_v^*\times\bfT_{w'}^*,s_{n'}}.
\end{multline*}
Now we interchange the order of the summations
\[
\frac{1}{|W_k|}\sum_{v\in W_k}\sum_{s\in T^*_{v,l}}
=\frac{1}{|W_k|}\sum_{s\in X_{k,l}}\sum_{v\in W_k(s)}
\quad\text{and}\quad
\sum_{k=0}^n\sum_{l=0}^k
=\sum_{l=0}^n\sum_{k=l}^n{},
\]
then we have
\begin{multline*}
\omega_{\Sp_{2n},\rmO^\epsilon_{2n'}}^\sharp
=\frac{1}{2}\sum_{l=0}^n\sum_{k=l}^n\frac{1}{|W_k|}\frac{1}{|\bfW_{\Sp_{2(n-k)}}|}\frac{1}{|\bfW_{\SO_{2(n'-k)}}|}
\sum_{s\in X_{k,l}}\sum_{v\in W_k(s)} \\
\sum_{w\in\bfW_{\Sp_{2(n-k)}}}\sum_{w'\in\bfW_{\SO_{2(n'-k)}^{\epsilon_v\epsilon}}}
\epsilon_v R^{\Sp_{2n}}_{\bfT_v^*\times\bfT_w^*,s_n}\otimes
R^{\rmO^{\epsilon}_{2n'}}_{\bfT_v^*\times\bfT_{w'}^*,s_{n'}}.
\end{multline*}
Because $|(s)|=\frac{|W_k|}{|W_k(s)|}$,
we have
\begin{multline*}
\omega_{\Sp_{2n},\rmO^\epsilon_{2n'}}^\sharp
=\frac{1}{2}\sum_{l=0}^n\sum_{k=l}^n\frac{1}{|\bfW_{\Sp_{2(n-k)}}|}\frac{1}{|\bfW_{\SO_{2(n'-k)}}|}
\sum_{(s)\subset X_{k,l}}\sum_{v\in W_k(s)} \\
\frac{1}{|W_k(s)|}\sum_{w\in\bfW_{\Sp_{2(n-k)}}}\sum_{w'\in\bfW_{\SO_{2(n'-k)}^{\epsilon_v\epsilon}}}
\epsilon_v R^{\Sp_{2n}}_{\bfT_v^*\times\bfT_w^*,s_n}\otimes
R^{\rmO^{\epsilon}_{2n'}}_{\bfT_v^*\times\bfT_{w'}^*,s_{n'}}.
\end{multline*}
For $k\geq l$,
each orbit in $X_{k,l}$ by the action of $W_k$ is of the form $(s_k)$ for some unique $(s)\subset X_{l,l}$,
so
\begin{multline*}
\omega_{\Sp_{2n},\rmO^\epsilon_{2n'}}^\sharp
=\frac{1}{2}\sum_{l=0}^n\sum_{(s)\subset X_{l,l}}\sum_{k=l}^n\frac{1}{|W_k(s_k)|}\frac{1}{|\bfW_{\Sp_{2(n-k)}}|}\frac{1}{|\bfW_{\SO_{2(n'-k)}}|}
\sum_{v\in W_k(s_k)} \\
\sum_{w\in\bfW_{\Sp_{2(n-k)}}}\sum_{w'\in\bfW_{\SO_{2(n'-k)}^{\epsilon_v\epsilon}}}
\epsilon_v R^{\Sp_{2n}}_{\bfT_v^*\times\bfT_w^*,s_n}\otimes
R^{\rmO^{\epsilon}_{2n'}}_{\bfT_v^*\times\bfT_{w'}^*,s_{n'}}.
\end{multline*}
Because now (\ref{0603}) is an orthogonal decomposition,
we must have
\begin{multline*}
\omega_{\Sp_{2n},\rmO^\epsilon_{2n'},s}^\sharp
=\frac{1}{2}\sum_{k=l}^{n}\frac{1}{|W_k(s_k)|}\frac{1}{|\bfW_{\Sp_{2(n-k)}}|}\frac{1}{|\bfW_{\SO_{2(n'-k)}}|}
\sum_{v\in W_k(s_k)}\\
\sum_{w\in\bfW_{\Sp_{2(n-k)}}}\sum_{w'\in\bfW_{\SO_{2(n'-k)}^{\epsilon_v\epsilon}}}
\epsilon_v R^{\Sp_{2n}}_{\bfT_v\times\bfT_w,s_n}\otimes R^{\rmO_{2n'}^\epsilon}_{\bfT_v\times\bfT_{w'},s_{n'}},
\end{multline*}
i.e., the lemma for case $n'>n$ is proved.

The proof for the other case is similar and omitted.
\end{proof}

\subsection{Lusztig correspondence of the Weil character}\label{0506}
Let $G$ be a member in a dual pair of a finite symplectic group and a finite even orthogonal group.
Let $s$ be a semisimple element in the connected component $(G^*)^0$ of dual group of $G$.
Define
\begin{itemize}
\item $G^{(0)}=G^{(0)}(s)=\prod_{\langle\lambda\rangle\subset\{\lambda_1,\ldots,\lambda_l\},\ \lambda\neq\pm 1}G_{[\lambda]}(s)$;

\item $G^{(1)}=G^{(1)}(s)=G_{[-1]}(s)$;

\item $G^{(2)}=G^{(2)}(s)=(G_{[1]}(s))^*$, the dual group of $G_{[1]}(s)$
\end{itemize}
where $G_{[\lambda]}(s)$ is given in Subsection~\ref{0501}.
Now $C_{G^*}(s)\simeq G^{(0)}\times G^{(1)}\times (G^{(2)})^*$,
hence we have a bijection
\[
\cale(C_{G^*}(s))_1\simeq\cale(G^{(0)}\times G^{(1)}\times G^{(2)})_1.
\]
Therefore we have a one-to-one correspondence
\begin{equation}\label{0505}
\Xi_s\colon\cale(G)_s\longrightarrow
\cale(G^{(0)}\times G^{(1)}\times G^{(2)})_1
\end{equation}
by Proposition~\ref{0201}.
The bijection $\Xi_s$ is also called a (modified) \emph{Lusztig correspondence}.
Now (\ref{0505}) can be extended linearly to an isometry (still denoted by $\Xi_s$) 
of inner product spaces:
\[
\Xi_s\colon\calv(G)_s\longrightarrow
\calv(G^{(0)})_1\otimes\calv(G^{(1)})_1\otimes\calv(G^{(2)})_1.
\]
Similar to (\ref{0524}), $s$ can be written as
\begin{equation}\label{0613}
s=s^{(0)}\times s^{(1)}\times s^{(2)}
\end{equation}
where $s^{(1)}$ (resp.~$s^{(2)}$) is the part whose eigenvalues are all equal to $-1$ (resp.~$1$),
and $s^{(0)}$ is the part whose eigenvalues do not contain $-1$ or $1$.

Now if $\rho$ is an irreducible character in $\cale(G)_s$,
$\Xi_s(\rho)$ can be written as $\rho^{(0)}\otimes\rho^{(1)}\otimes\rho^{(2)}$ where
$\rho^{(j)}\in\cale(G^{(j)})_1$ for $j=0,1,2$.
Moreover by Proposition~\ref{0202},
if $\rho$ is cuspidal, then each $\rho^{(j)}$ is cuspidal.
When $G$ is an even orthogonal group,
then both $G^{(1)}$ and $G^{(2)}$ are even orthogonal groups,
we will let $\{\rho_i\}$ denote the image under $\Xi_s^{-1}$ of the set
\begin{multline}\label{0508}
\{\rho^{(0)}\otimes\rho^{(1)}\otimes\rho^{(2)},\rho^{(0)}\otimes(\rho^{(1)}\cdot\sgn)\otimes\rho^{(2)}, \\
\rho^{(0)}\otimes\rho^{(1)}\otimes(\rho^{(2)}\cdot\sgn),
\rho^{(0)}\otimes(\rho^{(1)}\cdot\sgn)\otimes(\rho^{(2)}\cdot\sgn)\}.
\end{multline}
So we have $i=1,2,3,4$ if both $G^{(1)}$ and $G^{(2)}$ are nontrivial;
$i=1,2$ if exactly one of $G^{(1)}$ and $G^{(2)}$ is nontrivial;
and $i=1$ if both $G^{(1)}$ and $G^{(2)}$ are trivial.
Now suppose that $G^{(1)}$ and $G^{(2)}$ are non-trivial, so now
the Lusztig correspondence $\Xi_s$ is a bijection from $\{\rho_1,\rho_2,\rho_3,\rho_4\}$
to the set in (\ref{0508}).

The following lemma is extracted from the proof of \cite{amr} th\'eor\`eme 2.6.

\begin{lem}\label{0503}
Let $(\bfG,\bfG')=(\Sp_{2n},\rmO_{2n'}^\epsilon)$,
and $s\in X_{l,l}$ for some $l\leq\min(n,n')$.
Then
\begin{enumerate}
\item[(i)] $G^{(0)}(s_n)$ and $G'^{(0)}(s_{n'})$ are products of unitary groups or general linear groups,
and two groups are isomorphic;

\item[(ii)] $G^{(1)}(s_n)$ and $G'^{(1)}(s_{n'})$ are isomorphic even orthogonal groups;

\item[(iii)] $(G^{(2)}(s_n),G'^{(2)}(s_{n'}))$ is a dual pair of a symplectic group and an even orthogonal group.
\end{enumerate}
\end{lem}
\begin{proof}
Now $s_n$ is a semisimple element in $\SO_{2n+1}(q)$.
From Subsection~\ref{0501}, we know that $G^{(0)}(s_n)$ is a product of unitary groups or general linear groups,
and $G^{(1)}(s_n)=\rmO^{\epsilon'}_{2\nu_{-1}(s_n)}(q)$ for some $\epsilon'$ depending on $s$.
Moreover, $G^{(2)}(s_n)$ is the dual group of $\SO_{2\nu_1(s_n)+1}(q)$, hence
$G^{(2)}(s_n)=\Sp_{2\nu_1(s_n)}(q)$.

Similarly, $s_{n'}$ is a semisimple element in $\rmO_{2n'}^\epsilon(q)$,
and hence $G'^{(0)}(s_{n'})$ is a product of unitary groups,
and $G'^{(1)}(s_{n'})=\rmO^{\epsilon'}_{2\nu_{-1}(s_{n'})}(q)$.
Moreover, $G'^{(2)}(s_{n'})$ is the dual group of $\rmO_{2\nu_1(s_{n'})}^{\epsilon''}(q)$,
i.e., $G'^{(2)}(s_{n'})=\rmO_{2\nu_1(s_{n'})}^{\epsilon''}(q)$ for some $\epsilon',\epsilon''$ depending on $s$ and $\epsilon$.

Clearly $\nu_\lambda(s_n)=\nu_\lambda(s_{n'})=\nu_\lambda(s)$ when $\lambda\neq 1$,
so $G^{(0)}(s_n)\simeq G'^{(0)}(s_{n'})$ and $G^{(1)}(s_n)\simeq G'^{(1)}(s_{n'})$.
Moreover, $\nu_1(s_n)=n-l$, $\nu_1(s_{n'})=n'-l$,
and hence
\[
(G^{(2)}(s_n),G'^{(2)}(s_{n'}))=(\Sp_{2(n-l)}(q),\rmO_{2(n'-l)}^{\epsilon''}(q)).
\]
\end{proof}

For $s\in X_{l,l}$,
we will identify $G^{(0)}(s_{n})\simeq G'^{(0)}(s_{n'})$ and $G^{(1)}(s_{n})\simeq G'^{(1)}(s_{n'})$.
Let $\Xi_{(s_n,s_{n'})}$ denote $\Xi_{s_n}\otimes\Xi_{s_{n'}}$.
Then we have an isometry
\begin{multline*}
\Xi_{(s_n,s_{n'})}\colon\calv(\Sp_{2n}(q))_{s_n}\otimes\calv(\rmO^\epsilon_{2n'}(q))_{s_{n'}}\longrightarrow \\
\bigl[\calv(G^{(0)})_1\otimes\calv(G'^{(0)})_1\bigr]\otimes\bigl[\calv(G^{(1)})_1\otimes\calv(G'^{(1)})_1\bigr]
\otimes\bigl[\calv(G^{(2)})_1\otimes\calv(G'^{(2)})_1\bigr].
\end{multline*}
The following proposition is also from \cite{amr} th\'eor\`eme 2.6.

\begin{prop}\label{0604}
Let $(\bfG,\bfG')=(\Sp_{2n},\rmO^\epsilon_{2n'})$, $s\in X_{l,l}$ for some $l\leq\min(n,n')$.
Then
\[
\Xi_{(s_n,s_{n'})}(\omega_{\bfG,\bfG',s}^\sharp)
=R_{\bfG^{(0)},1}^\sharp\otimes R_{\bfG^{(1)},1}^\sharp \otimes\omega_{\bfG^{(2)},\bfG'^{(2)},1}^\sharp
\]
where $R_{\bfG^{(j)},1}^\sharp$ is given in Subsection~\ref{0309}.
\end{prop}
\begin{proof}
Suppose that $s\in X_{l,l}$.
From (\ref{0218}) and (\ref{0505}) we know that
\[
\Xi_{s_n}(R^{\bfG}_{\bfT_v^*\times\bfT_w^*,s_n})
=\epsilon_\bfG\epsilon_{\bfG^{(0)}\times\bfG^{(1)}\times\bfG^{(2)}}
R_{\bfT_v^*\times\bfT_w^*,1}^{\bfG^{(0)}\times\bfG^{(1)}\times\bfG^{(2)}}.
\]
Suppose that $n<n'$.
From Lemma~\ref{0504}, we have
\begin{multline}\label{0511}
\Xi_{(s_n,s_{n'})}(\omega_{\bfG,\bfG',s}^\sharp)
=\frac{1}{2}\sum_{k=l}^{n}\frac{1}{|W_k(s_k)|}\frac{1}{|\bfW_{\Sp_{2(n-k)}}|}\frac{1}{|\bfW_{\SO_{2(n'-k)}}|}
\sum_{v\in W_k(s_k)}\\
\sum_{w\in\bfW_{\Sp_{2(n-k)}}}\sum_{w'\in\bfW_{\SO_{2(n'-k)}^{\epsilon_v\epsilon}}}
\epsilon_v\epsilon_{\bfG}\epsilon_{\bfG^{(0)}\times\bfG^{(1)}\times\bfG^{(2)}}\epsilon_{\bfG'}\epsilon_{\bfG'^{(0)}\times\bfG'^{(1)}\times\bfG'^{(2)}}
R^{\bfG^{(0)}\times\bfG^{(1)}\times\bfG^{(2)}}_{\bfT_v^*\times\bfT_w^*,1}\otimes R^{\bfG'^{(0)}\times\bfG'^{(1)}\times\bfG'^{(2)}}_{\bfT_v^*\times\bfT_{w'}^*,1}.
\end{multline}

We know that $W_k(s_k)=\bfW_{\bfG^{(0)}}\times\bfW_{\bfG^{(1)}}\times W_{k-l}$
and write $v=(v_0,v_1,v_2)$ for $v_0\in\bfW_{\bfG^{(0)}}$, $v_1\in\bfW_{\bfG^{(1)}}$ and $v_2\in W_{k-l}$.
Then $\bfT^*_v=\bfT^*_{v_0}\times\bfT^*_{v_1}\times\bfT^*_{v_2}$ with
$\bfT^*_{v_0}\subset\bfG^{(0)}$, $\bfT^*_{v_1}\subset\bfG^{(1)}$ and $\bfT^*_{v_2}\times\bfT_w^*\subset\bfG^{(2)}$.
Moreover, it is easy to check that
$\epsilon_\bfG\epsilon_{\bfG^{(2)}}=(-1)^n\cdot(-1)^{n-l}=(-1)^l$;
$\epsilon_{\bfG^{(0)}}=\epsilon_{\bfG'^{(0)}}$ and $\epsilon_{\bfG^{(1)}}=\epsilon_{\bfG'^{(1)}}$
by Lemma~\ref{0503};
$\epsilon_{\bfG'}=(-1)^{n'}\epsilon_v\epsilon_{w'}$ and
$\epsilon_{\bfG'^{(2)}}=(-1)^{n'-l}\epsilon_{v_2}\epsilon_{w'}$ by the conditions that $\bfT_v^*\times\bfT_{w'}^*\subset\bfG'^*$
and $\bfT_{v_2}^*\times\bfT_{w'}^*\subset\bfG'^{(2)}$.
And hence
\begin{align*}
\epsilon_v\epsilon_\bfG\epsilon_{\bfG^{(0)}\times\bfG^{(1)}\times\bfG^{(2)}}\epsilon_{\bfG'}
\epsilon_{\bfG'^{(0)}\times\bfG'^{(1)}\times\bfG'^{(2)}}
&=\epsilon_v\epsilon_\bfG\epsilon_{\bfG^{(0)}}\epsilon_{\bfG^{(1)}}\epsilon_{\bfG^{(2)}}
\epsilon_{\bfG'}\epsilon_{\bfG'^{(0)}}\epsilon_{\bfG'^{(1)}}\epsilon_{\bfG'^{(2)}} \\
&=\epsilon_{v_2}.
\end{align*}
Therefore, (\ref{0511}) becomes
\begin{multline*}
\Xi_{(s_n,s_{n'})}(\omega_{\bfG,\bfG',s}^\sharp)
=\Biggl[\frac{1}{|\bfW_{\bfG^{(0)}}|}\sum_{v_1\in \bfW_{\bfG^{(0)}}}
R_{\bfT^*_{v_0},1}^{\bfG^{(0)}}\otimes R_{\bfT^*_{v_0},1}^{\bfG'^{(0)}}\Biggr]\otimes{}\\
\Biggl[\frac{1}{|\bfW_{\bfG^{(1)}}|}\sum_{v_1\in \bfW_{\bfG^{(1)}}}
R_{\bfT^*_{v_1},1}^{\bfG^{(1)}}\otimes R_{\bfT^*_{v_1},1}^{\bfG'^{(1)}}\Biggr]\otimes
\Biggl[\frac{1}{2}\sum_{k=l}^{n}\frac{1}{|W_{k-l}|}\frac{1}{|\bfW_{\Sp_{2(n-k)}}|}\frac{1}{|\bfW_{\SO_{2(n'-k)}}|}
\sum_{v_2\in W_{k-l}}\\
\sum_{w\in\bfW_{\Sp_{2(n-k)}}}\sum_{w'\in\bfW_{\SO_{2(n'-k)}^{\epsilon_v\epsilon}}}
\epsilon_{v_2} R^{\bfG^{(2)}}_{\bfT_{v_2}\times\bfT_w,1}\otimes R^{\bfG'^{(2)}}_{\bfT_{v_2}\times\bfT_{w'},1}\Biggr].
\end{multline*}
By (\ref{0305}) and (\ref{0306}),
we have
\[
\frac{1}{|\bfW_{\bfG^{(j)}}|}\sum_{v_1\in \bfW_{\bfG^{(j)}}}
R_{\bfT^*_{v_1},1}^{\bfG^{(j)}}\otimes R_{\bfT^*_{v_1},1}^{\bfG'^{(j)}}
=R_{\bfG^{(j),1}}^\sharp
\]
for $j=0,1$.
Apply Lemma~\ref{0504} with $s=1\in X_{0,0}$ and $k'=k-l$,
we have
\begin{multline*}
\frac{1}{2}\sum_{k=l}^{n}\frac{1}{|W_{k-l}|}\frac{1}{|\bfW_{\Sp_{2(n-k)}}|}\frac{1}{|\bfW_{\SO_{2(n'-k)}}|}
\sum_{v_2\in W_{k-l}}\\
\sum_{w\in\bfW_{\Sp_{2(n-k)}}}\sum_{w'\in\bfW_{\SO_{2(n'-k)}^{\epsilon_v\epsilon}}}
\epsilon_{v_2} R^{\bfG^{(2)}}_{\bfT_{v_2}\times\bfT_w,1}\otimes R^{\bfG'^{(2)}}_{\bfT_{v_2}\times\bfT_{w'},1}\\
=\frac{1}{2}\sum_{k'=0}^{n-l}\frac{1}{|W_{k'}|}\frac{1}{|\bfW_{\Sp_{2(n-l-k')}}|}\frac{1}{|\bfW_{\SO_{2(n'-l-k')}}|}
\sum_{v_2\in W_{k'}}\\
\sum_{w\in\bfW_{\Sp_{2(n-l-k')}}}\sum_{w'\in\bfW_{\SO_{2(n'-l-k')}^{\epsilon_v\epsilon}}}
\epsilon_{v_2} R^{\bfG^{(2)}}_{\bfT_{v_2}\times\bfT_w,1}\otimes R^{\bfG'^{(2)}}_{\bfT_{v_2}\times\bfT_{w'},1}
=\omega^\sharp_{\bfG^{(2)},\bfG'^{(2)},1}.
\end{multline*}
Note that now $\bfG^{(2)}=\Sp_{2(n-l)}$ and $\bfG'^{(2)}=\rmO^{\epsilon_v\epsilon}_{2(n'-l)}$.
Hence the proposition is proved for the case $n<n'$.

The proof for the case $n'\leq n$ is similar and omitted.
\end{proof}

\subsection{The main result I}
Now we have our first main result of this article:

\begin{thm}\label{0606}
Let $(\bfG,\bfG')=(\Sp_{2n},\rmO^\epsilon_{2n'})$,
and let $\rho\in\cale(G)_s$ and $\rho'\in\cale(G')_{s'}$ for some semisimple elements $s\in G^*$
and $s'\in (G'^*)^0$.
Write $\Xi_s(\rho)=\rho^{(0)}\otimes\rho^{(1)}\otimes\rho^{(2)}$ and\/
$\Xi_{s'}(\rho')=\rho'^{(0)}\otimes\rho'^{(1)}\otimes\rho'^{(2)}$,
and let $\{\rho'_i\}$ be given in (\ref{0508}).
Then one of the $\rho\otimes\rho'_i$'s occurs in $\omega_{\bfG,\bfG'}^\psi$ if and only if
the following conditions hold:
\begin{itemize}
\item $s^{(0)}=s'^{(0)}$ (up to conjugation), and $\rho^{(0)}=\rho'^{(0)}$;

\item $G^{(1)}\simeq G'^{(1)}$, and $\rho^{(1)}$ is equal to $\rho'^{(1)}$ or $\rho'^{(1)}\cdot\sgn$;

\item $\rho^{(2)}\otimes\rho'^{(2)}$ or $\rho^{(2)}\otimes(\rho'^{(2)}\cdot\sgn)$ occurs
$\omega_{\bfG^{(2)},\bfG'^{(2)},1}$.
\end{itemize}
\end{thm}
\begin{proof}
First suppose that $\rho\otimes\rho'_i$ occurs in $\omega^\psi_{\bfG,\bfG'}$ for some $i$.
Because the $\sgn$ character of $G'^{(1)}$ or $G'^{(2)}$ is of order $2$,
without loss of generality,
we may just assume that $\rho\otimes\rho'$ occurs in $\omega_{\bfG,\bfG'}^\psi$.
By (\ref{0603}), $\rho\otimes\rho'$ occurs in $\omega^\psi_{\bfG,\bfG',t}$ for some $t\in X_{l,l}$ 
and some $l\leq\min(n,n')$.
Therefore $s=t_n$ and $s'=t_{n'}$ up to conjugation,
and so $s^{(0)}=s'^{(0)}$ up to conjugation.
Write
\[
\Xi_{(s,s')}(\omega^\psi_{\bfG,\bfG',t})
=\Omega^{(0)}\otimes\Omega^{(1)}\otimes\Omega^{(2)}
\]
where $\Omega^{(j)}\in\calv(G^{(j)}\times G'^{(j)})$ for $j=0,1,2$.
We know that $\Omega^{(j)}$ is a character of $G^{(j)}\times G'^{(j)}$ and we have
$\rho^{(j)}\otimes\rho'^{(j)}$ occurs in $\Omega^{(j)}$ for each $j=0,1,2$.
\begin{itemize}
\item Now $G^{(0)}$ is isomorphic to $G'^{(0)}$ by Lemma~\ref{0503}.
We identify $G^{(0)}$ with $G'^{(0)}$ and by Proposition~\ref{0604},
we have $\Omega^{(0)\sharp}=R_{\bfG^{(0)},1}^\sharp$.
Because now $G^{(0)}$ is a product of unitary groups or general linear groups,
every character of $G^{(0)}\times G^{(0)}$ is uniform,
and hence $\Omega^{(0)}=R_{\bfG^{(0)},1}$.
This implies that $\rho^{(0)}=\rho'^{(0)}$.

\item By Lemma~\ref{0503},
$G^{(1)}=G'^{(1)}=\rmO_{2\nu_{-1}(t)}^{\epsilon'}(q)$ for some $\epsilon'$ depending on $t$.
By Proposition~\ref{0604},
we have $\Omega^{(1)\sharp}=R_{\bfG^{(1)},1}^\sharp$.
Then by Proposition~\ref{0311},
we see that either $\rho^{(1)}=\rho'^{(1)}$ or $\rho^{(1)}=\rho'^{(1)}\cdot\sgn$.

\item By Lemma~\ref{0503},
$(\bfG^{(2)},\bfG'^{(2)})$ is a dual pair of a symplectic group and an even orthogonal group.
By Proposition~\ref{0604}, we have $\Omega^{(2)\sharp}=\omega^\sharp_{\bfG^{(2)},\bfG'^{(2)},1}$.
Now both $\Omega^{(2)}$ and $\omega_{\bfG^{(2)},\bfG'^{(2)},1}$ are non-negative integral combination
of irreducible characters of $G^{(2)}\times G'^{(2)}$ whose uniform projections are equal.
Moreover, $\rho^{(2)}\otimes\rho'^{(2)}$ occurs in $\Omega^{(2)}$,
from the proof in \cite{pan-finite-unipotent} subsection 5.4,
we know that either $\rho^{(2)}\otimes\rho'^{(2)}$ or $\rho^{(2)}\otimes(\rho'^{(2)}\cdot\sgn)$
occurs in $\omega_{\bfG^{(2)},\bfG'^{(2)},1}$.
\end{itemize}

Conversely suppose that the following conditions hold:
\begin{itemize}
\item $s^{(0)}=s'^{(0)}$, and $\rho^{(0)}=\rho'^{(0)}$;

\item $G^{(1)}\simeq G'^{(1)}$, and $\rho^{(1)}$ is equal to $\rho'^{(1)}$ or $\rho'^{(1)}\cdot\sgn$;

\item $\rho^{(2)}\otimes\rho'^{(2)}$ or $\rho^{(2)}\otimes(\rho'^{(2)}\cdot\sgn)$ occurs in
$\omega_{\bfG^{(2)},\bfG'^{(2)},1}$.
\end{itemize}
The isomorphism $G^{(1)}\simeq G'^{(1)}$ implies $s^{(1)}=s'^{(1)}$.
Then we see that $s=t_n$ and $s'=t_{n'}$ up to conjugacy for some $t\in X_{l,l}$ for some $l\leq\min(n,n')$.
Now $\rho^{(0)}\otimes\rho'^{(0)}$ occurs in $R_{\bfG^{(1)},1}$,
$\rho^{(1)}\otimes\rho'^{(1)}$ or $\rho^{(1)}\otimes(\rho'^{(1)}\cdot\sgn)$ occurs in $R_{\bfG^{(1)},1}$,
and $\rho^{(2)}\otimes\rho'^{(2)}$ or $\rho^{(2)}\otimes(\rho'^{(2)}\cdot\sgn)$ occurs in
$\omega_{\bfG^{(2)},\bfG'^{(2)},1}$.
Hence we have some $\rho\otimes\rho'_i$ occurs in
\[
\Omega:=\Xi^{-1}_{(s,s')}(R_{\bfG^{(0)},1}\otimes R_{\bfG^{(1)},1}\otimes\omega_{\bfG^{(2)},\bfG'^{(3)},1}).
\]
Now both $\Omega$ and $\omega^\psi_{\bfG,\bfG',t}$ are non-negative integral combinations of irreducible
characters of $G\times G'$ whose uniform projection are the same by Proposition~\ref{0604}.
Moreover, some $\rho\otimes\rho'_i$ occurs in $\Omega$,
from the proof in \cite{pan-finite-unipotent} subsection 5.4,
we see that some $\rho\otimes\rho'_i$ occurs in $\omega^\psi_{\bfG,\bfG',t}$.
Therefore, some $\rho\otimes\rho'_i$ occurs in $\omega^\psi_{\bfG,\bfG'}$.
\end{proof}

Recall that for a dual pair $(\bfG,\bfG')$, $\rho\in\cale(G)$, $\rho'\in\cale(G')$,
we define
\begin{align*}
\Theta_{\bfG'}^\psi(\rho) &=\{\,\rho'\in\cale(G')\mid (\rho,\rho')\in\Theta_{\bfG,\bfG'}^\psi\,\}, \\
\Theta_{\bfG}^\psi(\rho') &=\{\,\rho\in\cale(G)\mid (\rho,\rho')\in\Theta_{\bfG,\bfG'}^\psi\,\}.
\end{align*}

\begin{cor}
Let $(\bfG,\bfG')=(\Sp_{2n},\rmO_{2n'}^\epsilon)$.
Suppose that $(\rho,\rho')\in\Theta^\psi_{\bfG,\bfG'}$, and $\rho\in\cale(G)_s$, $\rho'\in\cale(G')_{s'}$.
Then $\Theta^\psi_{\bfG'}(\rho)\subset\cale(G')_{s'}$ and $\Theta^\psi_{\bfG}(\rho')\subset\cale(G)_{s}$.
\end{cor}
\begin{proof}
Suppose that $\rho\in\cale(G)_s$, $\rho'\in\cale(G')_{s'}$ and
$(\rho,\rho')\in\Theta^\psi_{\bfG,\bfG'}$.
From the proof of Theorem~\ref{0606}, we know that there exists a semisimple element $t$
such that $s=t_n$ and $s'=t_{n'}$ up to a conjugation, i.e.,
the conjugacy class of $s'$ is uniquely determined by $s$ and vice versa.
\end{proof}


\section{Compatibility for Odd Orthogonal Groups}

In this section,
we consider the dual pair $(\Sp_{2n},\SO_{2n'+1})$ over a finite field of odd characteristic.
Most of the arguments in this section are parallel to those in the previous section and
will be sketchy.

\subsection{Decomposition of the Weil character}
The following proposition is from \cite{pan-odd} theorem 3.8:

\begin{prop}\label{0308}
Consider the dual pair $(\Sp_{2n},\SO_{2n'+1})$ over a finite field of odd characteristic.
We have the decomposition
\begin{multline*}
\omega_{\Sp_{2n},\SO_{2n'+1}}^\sharp\cdot(1\otimes\chi_{\SO_{2n'+1}})\\
=\sum_{k=0}^{\min(n,n')}\frac{1}{|W_k|}\frac{1}{|\bfW_{\Sp_{2(n-k)}}|}\frac{1}{|\bfW_{\SO_{2(n'-k)+1}}|}
\sum_{v\in W_k}\sum_{\theta\in\cale(T_v)}\\
\sum_{w\in\bfW_{\Sp_{2(n-k)}}}\sum_{w'\in\bfW_{\SO_{2(n'-k)+1}}}
\epsilon_w R^{\Sp_{2n}}_{\bfT_v\times\bfT_w,\theta\otimes\theta_w}\otimes R^{\SO_{2n'+1}}_{\bfT_v\times\bfT_{w'},\theta\otimes\theta_{w'}}
\end{multline*}
where $\theta_w$ is given in Subsection~\ref{0310} and $\chi_{\SO_{2n'+1}}$ is given in (\ref{0320}).
\end{prop}

For a semisimple element $s$ in a maximal torus $T^*$ of $\SO_{2k+1}(q)$,
and for $\epsilon=\pm$ and $k\leq n$,
let $s^\flat_{n,\epsilon}$ denote the semisimple element $(s,-1)\in\SO_{2k+1}(q)\times\SO_{2(n-k)}^\epsilon(q)\subseteq\SO_{2n+1}(q)$.
If $s$ in a maximal torus $T^*$ in $\Sp_{2k}(q)$ and $n'\geq k$,
then $s^\flat_{n'}=(s,-1)$ is a semisimple element in $\Sp_{2k}(q)\times\Sp_{2(n'-k)}(q)\subseteq\Sp_{2n'}(q)$.
Similar to the case for even orthogonal groups,
we can rewrite the decomposition in Proposition~\ref{0308} as
\begin{multline*}
\omega_{\Sp_{2n},\SO_{2n'+1}}^\sharp\cdot(1\otimes\chi_{\SO_{2n'+1}})\\
=\sum_{k=0}^{\min(n,n')}\frac{1}{|W_k|}\frac{1}{|\bfW_{\Sp_{2(n-k)}}|}\frac{1}{|\bfW_{\SO_{2(n'-k)+1}}|}
\sum_{v\in W_k}\sum_{s\in T_v^*}\\
\sum_{w\in\bfW_{\Sp_{2(n-k)}}}\sum_{w'\in\bfW_{\SO_{2(n'-k)+1}}}
\epsilon_w R^{\Sp_{2n}}_{\bfT_v^*\times\bfT_w^*,s_{n,\epsilon_w}^\flat}\otimes R^{\SO_{2n'+1}}_{\bfT_v^*\times\bfT_{w'}^*,s^\flat_{n'}}.
\end{multline*}
For $k\leq\min(n,n')$, $v\in W_k$, and $l\leq k$,
we define
\begin{align*}
T_{v,l}^{\flat*} &=\{\,s\in T_v^*\mid \nu_{-1}(s)=k-l\,\},\\
X_{k,l}^\flat &=\bigcup_{v\in W_k}T^{\flat*}_{v,l}
\end{align*}
where $\nu_{-1}(s)$ is given in Subsection~\ref{0501}.
Then $X_{k,l}^\flat$ is a subset of $\SO_{2k+1}(q)$ with an action by $W_k$.

For $s\in X_{l,l}^\flat$ with $l\leq\min(n,n')$,
note that the two elements $s^\flat_{n,+}$ and $s^\flat_{n,-}$ in $\SO_{2n+1}(q)$ are not conjugate.
For $\epsilon=\pm$,
let $\omega^\psi_{\Sp_{2n},\SO_{2n'+1},s,\epsilon}$ denote the orthogonal projection of
$\omega^\psi_{\Sp_{2n},\SO_{2n'+1}}\cdot(1\otimes\chi_{\SO_{2n'+1}})$ over
$\calv(\Sp_{2n}(q))_{s_{n,\epsilon}^\flat}\otimes\calv(\SO_{2n'+1}(q))_{s_{n'}^\flat}$.
Similar to (\ref{0603}),
we have the decomposition
\begin{equation}\label{0608}
\omega^\psi_{\Sp_{2n},\SO_{2n'+1}}\cdot(1\otimes\chi_{\SO_{2n'+1}})
=\sum_{l=0}^{\min(n,n')}\sum_{(s)\subset X^\flat_{l,l}}
(\omega^\psi_{\Sp_{2n},\SO_{2n'+1},s,+} +\omega^\psi_{\Sp_{2n},\SO_{2n'+1},s,-}).
\end{equation}

\begin{lem}\label{0701}
Suppose $s\in X^\flat_{l,l}$ for some $l\leq\min(n,n')$ and $\epsilon=\pm$.
Then
\begin{multline*}
\omega_{\Sp_{2n},\SO_{2n'+1},s,\epsilon}^\sharp
=\sum_{k=l}^{\min(n,n')}\frac{1}{|W_k(s^\flat_{k,\epsilon})|}\frac{1}{|\bfW_{\Sp_{2(n-k)}}|}\frac{1}{|\bfW_{\SO_{2(n'-k)+1}}|}
\sum_{v\in W_k(s^\flat_{k,\epsilon})}\\
\sum_{w\in\bfW^\epsilon_{\Sp_{2(n-k)}}}
\sum_{w'\in\bfW_{\SO_{2(n'-k)+1}}}
\epsilon R^{\Sp_{2n}}_{\bfT^*_v\times\bfT^*_w,s_{n,\epsilon}^\flat}\otimes R^{\SO_{2n'+1}}_{\bfT^*_v\times\bfT^*_{w'},s_{n'}^\flat}.
\end{multline*}
\end{lem}
\begin{proof}
The proof is similar to that of Lemma~\ref{0504}.
In fact, we only need to notice that for $k\geq l$,
each orbit in $X_{k,l}^\flat$ by the action of $W_k$ is of the form $(s_{k,\epsilon})$ for $\epsilon=\pm$
and some unique $(s)\subset X_{l,l}^\flat$.
Note that $\epsilon_w=\epsilon$ if $w\in\bfW^\epsilon_{\Sp_{2(n-k)}}$.
\end{proof}

\subsection{Lusztig correspondence of the Weil character}\label{0607}
Let $G$ be a member in a finite dual pair of a symplectic group and a special odd orthogonal group.
Let $s$ be a semisimple element in the dual group $G^*$ of $G$.
\begin{enumerate}
\item Suppose that $G$ is a symplectic group.
Then $G^*$ is a special odd orthogonal group.
We define
\begin{itemize}
\item $G^{(0)}=G^{(0)}(s)
=\prod_{\langle\lambda\rangle\subset\{\lambda_1,\ldots,\lambda_l\},\ \lambda\neq\pm 1}G_{[\lambda]}(s)$;

\item $G^{(1)}=G^{(1)}(s)=G_{[-1]}(s)$;

\item $G^{(2)}=G^{(2)}(s)=(G_{[1]}(s))^*$, the dual group of $G_{[1]}(s)$
\end{itemize}
where $G_{[\lambda]}(s)$ is given in Subsection~\ref{0501}.
Then we have
\[
C_{G^*}(s)\simeq G^{(0)}\times G^{(1)}\times (G^{(2)})^*,
\]
where $G^{(1)}(s)$ is an even orthogonal group and $G^{(2)}(s)$ is a symplectic group.
As in (\ref{0505}) we have a (modified) Lusztig correspondence
\begin{equation}\label{0601}
\Xi_s\colon\cale(G)_s\longrightarrow
\cale(G^{(0)}\times G^{(1)}\times G^{(2)})_1.
\end{equation}
So we will write $\Xi_s(\rho)=\rho^{(0)}\otimes\rho^{(1)}\otimes\rho^{(2)}$
where $\rho^{(j)}$ is in $\cale(G^{(j)})_1$ for $j=0,1,2$.
Let $\{\rho_i\}$ denote the image under $\Xi_s^{-1}$ of the set
\begin{equation}\label{0709}
\{\rho^{(0)}\otimes\rho^{(1)}\otimes\rho^{(2)},\rho^{(0)}\otimes(\rho^{(1)}\cdot\sgn)\otimes\rho^{(2)}\}.
\end{equation}
Hence $i=1,2$ if $G^{(1)}$ is nontrivial; and $i=1$ otherwise.

\item Suppose that $G$ is a special odd orthogonal group.
Then $G^*$ is a symplectic group.
We define
\begin{itemize}
\item $G^{(0)}=G^{(0)}(s)
=\prod_{\langle\lambda\rangle\subset\{\lambda_1,\ldots,\lambda_l\},\ \lambda\neq\pm 1}G_{[\lambda]}(s)$;

\item $G^{(1)}=G^{(1)}(s)=G_{[-1]}(s)$;

\item $G^{(2)}=G^{(2)}(s)=G_{[1]}(s)$,
\end{itemize}
Now
\[
C_{G^*}(s)\simeq G^{(0)}\times G^{(1)}\times G^{(2)}
\]
where $G^{(1)}$ and $G^{(2)}$ are both symplectic groups.
Let
\begin{equation}\label{0609}
\Xi_s\colon\cale(G)_s\longrightarrow
\cale(G^{(0)}\times G^{(1)}\times G^{(2)})_1
\end{equation}
be a (modified) Lusztig correspondence and write
$\Xi_s(\rho)=\rho^{(0)}\otimes\rho^{(1)}\otimes\rho^{(2)}$.
\end{enumerate}
We can extend $\Xi_s$ to be an isometry
\[
\Xi_s\colon\calv(G)_s\rightarrow\calv(G^{(0)})_1\otimes\calv(G^{(1)})_1\otimes\calv(G^{(2)})_1.
\]
As in (\ref{0613}), the semisimple element $s$ can be written as
\begin{equation}
s=s^{(0)}\times s^{(1)}\times s^{(2)}
\end{equation}
where $s^{(1)}$ (resp.~$s^{(2)}$) is the part whose eigenvalues are all equal to $-1$ (resp.~$1$),
and $s^{(0)}$ is the part whose eigenvalues do not contain $-1$ or $1$.

\begin{lem}\label{0605}
Let $(\bfG,\bfG')=(\Sp_{2n},\SO_{2n'+1})$, $\epsilon=\pm$,
and $s\in X^\flat_{l,l}$ with $l\leq\min(n,n')$.
Then
\begin{enumerate}
\item[(i)] $G^{(0)}(s_{n,\epsilon}^\flat)$ and $G'^{(0)}(s^\flat_{n'})$ are product of unitary groups or general linear groups,
and two groups are isomorphic;

\item[(ii)] $(G^{(1)}(s_{n,\epsilon}^\flat),G'^{(1)}(s^\flat_{n'}))$ is
a dual pair of an even orthogonal group and a symplectic group;

\item[(iii)] $G^{(2)}(s_{n,\epsilon}^\flat)$ and $G'^{(2)}(s^\flat_{n'})$ are isomorphic symplectic groups.
\end{enumerate}
\end{lem}
\begin{proof}
Now $s_{n,\epsilon}^\flat$ is a semisimple element in $G^*=\SO_{2n+1}(q)$.
From Subsection~\ref{0501}, we know that $G^{(0)}(s_{n,\epsilon}^\flat)$ is a product of unitary groups or general linear groups,
$G^{(1)}(s^\flat_{n,\epsilon})=\rmO_{2\nu_{-1}(s_{n,\epsilon}^\flat)}^\epsilon(q)$,
and $G^{(2)}(s_{n,\epsilon}^\flat)=\Sp_{2\nu_1(s^\flat_{n,\epsilon})}(q)$.
Similarly, $s_{n'}^\flat$ is a semisimple element in $G'^*=\Sp_{2n'}(q)$,
and hence $G'^{(0)}(s_{n'}^\flat)$ is a product of unitary groups or general linear groups,
$G'^{(1)}(s_{n'}^\flat)=\Sp_{2\nu_{1}(s_{n'}^\flat)}(q)$,
and $G'^{(2)}(s_{n'}^\flat)=\Sp_{2\nu_{-1}(s_{n'}^\flat)}(q)$.

Clearly $\nu_\lambda(s_{n,\epsilon}^\flat)=\nu_\lambda(s_{n'}^\flat)=\nu_\lambda(s)$ when $\lambda\neq -1$,
so $G^{(0)}(s_{n,\epsilon}^\flat)$ and $G'^{(0)}(s^\flat_{n'})$ are isomorphic;
$G^{(2)}(s_{n,\epsilon}^\flat)$ and $G'^{(2)}(s^\flat_{n'})$ are also isomorphic.
Finally, $\nu_{-1}(s_{n,\epsilon}^\flat)=n-l$, $\nu_{-1}(s_{n'}^\flat)=n'-l$,
and hence
$(G^{(1)}(s_{n,\epsilon}^\flat),G'^{(1)}(s^\flat_{n'}))
=(\rmO_{2(n-l)}^\epsilon(q),\Sp_{2(n'-l)}(q))$.
\end{proof}

Let $\Xi_{(s^\flat_{n,\epsilon},s^\flat_{n'})}$ denote $\Xi_{s^\flat_{n,\epsilon}}\otimes\Xi_{s^\flat_{n'}}$.
Then we have an isometry
\begin{multline*}
\Xi_{(s^\flat_{n,\epsilon},s^\flat_{n'})}\colon\calv(\Sp_{2n}(q))_{s^\flat_{n,\epsilon}}\otimes\calv(\SO_{2n'+1}(q))_{s^\flat_{n'}}\longrightarrow \\
\bigl[\calv(G^{(0)})_1\otimes\calv(G'^{(0)})_1\bigr]\otimes\bigl[\calv(G^{(1)})_1\otimes\calv(G'^{(1)})_1\bigr]
\otimes\bigl[\calv(G^{(2)})_1\otimes\calv(G'^{(2)})_1\bigr].
\end{multline*}
Then again we should have the analogue of Proposition~\ref{0604} for the dual pair of a symplectic group
and a special odd orthogonal group:

\begin{prop}\label{0602}
Let $(\bfG,\bfG')=(\Sp_{2n},\SO_{2n'+1})$, $\epsilon=\pm$,
$s\in X^\flat_{l,l}$ with $l\leq\min(n,n')$.
Then
\[
\Xi_{(s_{n,\epsilon}^\flat,s_{n'}^\flat)}(\omega_{\bfG,\bfG',s,\epsilon}^\sharp)
=R_{\bfG^{(0)},1}^\sharp\otimes\omega_{\bfG^{(1)},\bfG'^{(1)},1}^\sharp
\otimes R_{\bfG^{(2)},1}^\sharp.
\]
\end{prop}
\begin{proof}
The proof is similar to that of Proposition~\ref{0604}.
Note that
\begin{align*}
\Xi_{s^\flat_{n,\epsilon}}(R^{\bfG}_{\bfT_v^*\times\bfT_w^*,s^\flat_{n,\epsilon}})
&=\epsilon_\bfG\epsilon_{\bfG^{(0)}\times\bfG^{(1)}\times\bfG^{(2)}}
R_{\bfT_v^*\times\bfT_w^*,1}^{\bfG^{(0)}\times\bfG^{(1)}\times\bfG^{(2)}},\\
\Xi_{s^\flat_{n'}}(R^{\bfG'}_{\bfT_v^*\times\bfT_{w'}^*,s^\flat_{n'}})
&=\epsilon_{\bfG'}\epsilon_{\bfG'^{(0)}\times\bfG'^{(1)}\times\bfG'^{(2)}}
R_{\bfT_v^*\times\bfT_{w'}^*,1}^{\bfG'^{(0)}\times\bfG'^{(1)}\times\bfG'^{(2)}}.
\end{align*}
We know that $W_k(s_k)=\bfW_{\bfG^{(0)}}\times\bfW_{\bfG^{(1)}}\times W_{k-l}$
and write $v=(v_0,v_1,v_2)$ for $v_0\in\bfW_{\bfG^{(0)}}$, $v_1\in\bfW_{\bfG^{(1)}}$ and $v_2\in W_{k-l}$.
Then by Lemma~\ref{0701} and the same argument in the proof of Proposition~\ref{0604},
we have
\begin{multline*}
\Xi_{(s^\flat_{n,\epsilon},s^\flat_{n'})}(\omega_{\bfG,\bfG',s,\epsilon}^\sharp)
=\Biggl[\frac{1}{|\bfW_{\bfG^{(0)}}|}\sum_{v_0\in \bfW_{\bfG^{(0)}}}
R_{\bfT^*_{v_0},1}^{\bfG^{(0)}}\otimes R_{\bfT^*_{v_0},1}^{\bfG'^{(0)}}\Biggr]\otimes{}\\
\Biggl[\sum_{k=l}^{n}\frac{1}{|W_{k-l}|}\frac{1}{|\bfW_{\Sp_{2(n-k)}}|}\frac{1}{|\bfW_{\SO_{2(n'-k)+1}}|}
\sum_{v_1\in W_{k-l}}\\
\sum_{w\in\bfW^\epsilon_{\Sp_{2(n-k)}}}\sum_{w'\in\bfW_{\SO_{2(n'-k)+1}}}
\epsilon_{v_1} R^{\bfG^{(1)}}_{\bfT_{v_1}\times\bfT_w,1}\otimes R^{\bfG'^{(1)}}_{\bfT_{v_1}\times\bfT_{w'},1}\Biggr]\otimes \\
\Biggl[\frac{1}{|\bfW_{\bfG^{(2)}}|}\sum_{v_2\in \bfW_{\bfG^{(2)}}}
R_{\bfT^*_{v_2},1}^{\bfG^{(2)}}\otimes R_{\bfT^*_{v_2},1}^{\bfG'^{(2)}}\Biggr].
\end{multline*}
For the middle part, we apply Lemma~\ref{0504} with $s=1\in X_{0,0}$ and $k'=k-l$,
$|\bfW_{\Sp_{2(n-k)}}|=2|\bfW_{\SO^\epsilon_{2(n-k)}}|$,
and identify $\bfW^\epsilon_{\Sp_{2(n-k)}}=\bfW_{\SO^\epsilon_{2(n-k)}}$ and
$\bfW_{\SO_{2(n'-k)+1}}=\bfW_{\Sp_{2(n'-k)}}$,
we conclude that
\begin{multline*}
\sum_{k=l}^{n}\frac{1}{|W_{k-l}|}\frac{1}{|\bfW_{\Sp_{2(n-k)}}|}\frac{1}{|\bfW_{\SO_{2(n'-k)+1}}|}
\sum_{v_1\in W_{k-l}}\\
\sum_{w\in\bfW^\epsilon_{\Sp_{2(n-k)}}}\sum_{w'\in\bfW_{\SO_{2(n'-k)+1}}}
\epsilon_{v_1} R^{\bfG^{(1)}}_{\bfT_{v_1}\times\bfT_w,1}\otimes R^{\bfG'^{(1)}}_{\bfT_{v_1}\times\bfT_{w'},1}\\
=\frac{1}{2}\sum_{k'=0}^{n-l}\frac{1}{|W_{k'}|}\frac{1}{|\bfW_{\SO^\epsilon_{2(n-l-k')}}|}\frac{1}{|\bfW_{\Sp_{2(n'-l-k')}}|}
\sum_{v_1\in W_{k'}}\\
\sum_{w\in\bfW_{\SO^{\epsilon_{v_1}\epsilon}_{2(n-l-k')}}}\sum_{w'\in\bfW_{\Sp_{2(n'-l-k')}}}
\epsilon_{v_1} R^{\bfG^{(1)}}_{\bfT_{v_1}\times\bfT_w,1}\otimes R^{\bfG'^{(1)}}_{\bfT_{v_1}\times\bfT_{w'},1}
=\omega^\sharp_{\bfG^{(1)},\bfG'^{(1)},1}.
\end{multline*}
Note that now $\bfG^{(1)}=\rmO^\epsilon_{2(n-l)}$ and $\bfG'^{(1)}=\Sp_{2(n'-l)}$.
Then the proposition is proved.
\end{proof}

\subsection{The main result II}\label{0712}
Now we have the compatibility of the Lusztig correspondence and the Howe correspondence for a dual pair
of a symplectic group and a special odd orthogonal group.

\begin{thm}\label{0502}
Let $(\bfG,\bfG')=(\Sp_{2n},\SO_{2n'+1})$,
and let $\rho\in\cale(G)_s$ and $\rho'\in\cale(G')_{s'}$ for some semisimple elements $s\in G^*$
and $s'\in G'^*$.
Write $\Xi_s(\rho)=\rho^{(0)}\otimes\rho^{(1)}\otimes\rho^{(2)}$,
$\Xi_{s'}(\rho')=\rho'^{(0)}\otimes\rho'^{(1)}\otimes\rho'^{(2)}$,
and let $\{\rho_i\}$ be defined as in (\ref{0709}).
Then one of the $\rho_i\otimes\rho'$'s occurs in $\omega^\psi_{\bfG,\bfG'}$
if and only if the following conditions hold:
\begin{itemize}
\item $s^{(0)}=-s'^{(0)}$ (up to conjugation), and $\rho^{(0)}=\rho'^{(0)}$;

\item $G^{(2)}(s)\simeq G'^{(1)}(s')$, $\rho^{(2)}=\rho'^{(1)}$;

\item $\rho^{(1)}\otimes\rho'^{(2)}$ or $(\rho^{(1)}\cdot\sgn)\otimes\rho'^{(2)}$ occurs in
$\omega_{\bfG^{(1)}(s),\bfG'^{(2)}(s'),1}$.
\end{itemize}
\end{thm}
\begin{proof}
Let $\rho\in\cale(G)_s$ and $\rho'\in\cale(G')_{s'}$.
It is clear that $\rho'\chi_{\bfG'}\in\cale(G')_{-s'}$,
$\bfG'^{(0)}(-s')=\bfG'^{(0)}(s')$, $\bfG'^{(1)}(-s')=\bfG'^{(2)}(s')$,
$\bfG'^{(2)}(-s')=\bfG'^{(1)}(s')$, and
\begin{equation}\label{0710}
\Xi_{-s'}(\rho'\chi_{\bfG'})
=\rho'^{(0)}\otimes\rho'^{(2)}\otimes\rho'^{(1)}
\end{equation}
(\cf.~Subsection~\ref{0310}).

First suppose that $\rho_i\otimes\rho'$ occurs in $\omega^\psi_{\bfG,\bfG'}$ for some $i$.
Without loss of the generality, we may just assume that $\rho\otimes\rho'$ occurs in $\omega^\psi_{\bfG,\bfG'}$.
So $\rho\otimes(\rho'\chi_{\bfG'})$ occurs in $\omega^\psi_{\bfG,\bfG'}\cdot(1\otimes\chi_{\bfG'})$.
By (\ref{0608}),
we know that $\rho\otimes(\rho'\chi_{\bfG'})$ occurs in $\omega^\psi_{\bfG,\bfG',t,\epsilon}$
for some $\epsilon=\pm$ and some $t\in X^\flat_{l,l}$ for some integer $l\leq\min(n,n')$.
Then $s=t^\flat_{n,\epsilon}$ and $-s'=t^\flat_{n'}$ up to conjugation.
Write
\[
\Xi_{(t^\flat_{n,\epsilon},t^\flat_{n'})}(\omega^\psi_{\bfG,\bfG',t,\epsilon})
=\Omega^{(0)}\otimes\Omega^{(1)}\otimes\Omega^{(2)}
\]
where $\Omega^{(j)}\in\calv(G^{(j)}(s)\times G'^{(j)}(-s'))$ for $j=0,1,2$.
Then $\Omega^{(j)}$ is a character of $G^{(j)}(s)\times G'^{(j)}(-s')$
and $\rho^{(j)}\otimes\rho'^{(j)}$ occurs in $\Omega^{(j)}$.
\begin{itemize}
\item Now $s^{(0)}=-s'^{(0)}$ up to conjugation,
$G^{(0)}(s)$ is a product of unitary groups or general linear groups and isomorphic to
$G'^{(0)}(-s')\simeq\ G'^{(0)}(s')$ by Lemma~\ref{0605}.
We identity $G^{(0)}(s)$ and $G'^{(0)}(s')$ and by Proposition~\ref{0602},
we have $\Omega^{(0)\sharp}=R_{G^{(0)},1}^\sharp$.
Because now $G^{(0)}$ is a product of unitary groups or general linear groups,
every irreducible character of $G^{(0)}\times G^{(0)}$ is uniform,
and hence $\Omega^{(0)}=R_{G^{(0)}(s),1}$.
This implies that $\rho^{(0)}=\rho'^{(0)}$.

\item By Lemma~\ref{0605},
$G^{(2)}(s)=G'^{(2)}(-s')=G'^{(1)}(s')=\Sp_{2\nu_1(s)}(q)$.
By Proposition~\ref{0602},
we have $\Omega^{(2)\sharp}=R_{G^{(2)}(s),1}^\sharp$.
Then by Proposition~\ref{0311},
we have $\rho^{(2)}=\rho'^{(1)}$.

\item By Lemma~\ref{0605},
$(\bfG^{(1)}(s),\bfG'^{(1)}(-s'))=(\bfG^{(1)}(s),\bfG'^{(2)}(s'))$ is a dual pair
of an even orthogonal group and a symplectic group,
and by Proposition~\ref{0602} we have $\Omega^{(1)\sharp}=\omega^\sharp_{\bfG^{(1)}(s),\bfG'^{(2)}(s'),1}$.
Now $\Omega^{(1)}$ is a non-negative integral combination of irreducible unipotent characters of
$G^{(1)}(s)\times G'^{(2)}(s')$ whose uniform projection is equal to the uniform projection of
$\omega_{\bfG^{(1)}(s),\bfG'^{(2)}(s'),1}$.
By our assumption $\rho^{(1)}\otimes\rho'^{(2)}$ occurs in $\Omega^{(2)}$,
from the proof in \cite{pan-finite-unipotent} subsection~5.4,
we know that either $\rho^{(1)}\otimes\rho'^{(2)}$ or $(\rho^{(1)}\cdot\sgn)\otimes\rho'^{(2)}$
occurs in $\omega_{\bfG^{(1)}(s),\bfG'^{(2)}(s'),1}$.
\end{itemize}

Conversely suppose that
\begin{itemize}
\item $s^{(0)}=-s'^{(0)}$ and $\rho^{(0)}=\rho'^{(0)}$,

\item $G^{(2)}(s)\simeq G'^{(1)}(s')$ and $\rho^{(2)}=\rho'^{(1)}$,

\item $\rho^{(1)}\otimes\rho'^{(2)}$ or $(\rho^{(1)}\cdot\sgn)\otimes\rho'^{(2)}$ occurs in
$\omega_{\bfG^{(1)}(s),\bfG'^{(2)}(s'),1}$.
\end{itemize}
Then we see that $s=t^\flat_{n,\epsilon}$ and $-s'=t^\flat_{n'}$, up to conjugation,
for some $t\in X^\flat_{l,l}$ and some $l\leq\min(n,n')$.
Now $\rho^{(0)}\otimes\rho'^{(0)}$ occurs in $R_{\bfG^{(0)}(s),1}$,
$\rho^{(2)}\otimes\rho'^{(1)}$ occurs in $R_{\bfG^{(2)}(s),1}$,
and $\rho^{(1)}\otimes\rho'^{(2)}$ or $(\rho^{(1)}\cdot\sgn)\otimes\rho'^{(2)}$ occurs in
$\omega_{\bfG^{(2)}(s),\bfG'^{(2)}(s'),1}$.
Hence we have some $\rho_i\otimes(\rho'\chi_{\bfG'})$ occurs in the inverse image
\[
\Omega:=\Xi^{-1}_{(s,-s')}(R_{\bfG^{(0)}(s),1}\otimes R_{\bfG^{(2)}(s),1}
\otimes\omega_{\bfG^{(1)}(s),\bfG'^{(2)}(s'),1}).
\]
By Proposition~\ref{0602}, we know that $\Omega$ and $\omega^\psi_{\bfG,\bfG',t,\epsilon}$ have the same uniform projection.
From the proof in \cite{pan-finite-unipotent} subsection 5.4,
we see that $\rho_i\otimes(\rho'\chi_{\bfG'})$ occurs in $\omega^\psi_{\bfG,\bfG',t,\epsilon}$ for some $i$.
Therefore, by (\ref{0608}),
some $\rho_i\otimes(\rho'\chi_{\bfG'})$ occurs in $\omega^\psi_{\bfG,\bfG'}\cdot(1\otimes\chi_{\bfG'})$,
i.e., some $\rho_i\otimes\rho'$ occurs in $\omega^\psi_{\bfG,\bfG'}$.
\end{proof}

Let $\bfG=\Sp_{2n}$ and $s\in G^*$.
We can write $\bfG^{(1)}(s)=\rmO^{\epsilon^{(1)}}_{2n^{(1)}}$ for some $\epsilon^{(1)}$ and $n^{(1)}$.
It is known that there exists a unique class of semisimple $s^\dag\in G^*$ such that
$\bfG^{(0)}(s^\dag)\simeq \bfG^{(0)}(s)$,
$\bfG^{(1)}(s^\dag)\simeq \rmO^{-\epsilon^{(1)}}_{2n^{(1)}}$,
$\bfG^{(2)}(s^\dag)\simeq \bfG^{(2)}(s)$.
Clearly $(s^\dag)^\dag=s$.

\begin{cor}\label{0711}
Let $(\bfG,\bfG')=(\Sp_{2n},\SO_{2n'+1})$.
Suppose that $(\rho,\rho')\in\Theta^\psi_{\bfG,\bfG'}$, and $\rho\in\cale(G)_s$, $\rho'\in\cale(G')_{s'}$.
Then
\begin{enumerate}
\item[(i)] $\Theta^\psi_{\bfG'}(\rho)\subset\cale(G')_{s'}$,

\item[(ii)] $\Theta^\psi_{\bfG}(\rho')\subset\cale(G)_{s}\cup \cale(G)_{s^\dag}$.
\end{enumerate}
\end{cor}
\begin{proof}
Suppose that $\rho\in\cale(G)_s$ and $\rho'\in\cale(G')_{s'}$, and $(\rho,\rho')\in\Theta^\psi_{\bfG,\bfG'}$.
From the proof of Theorem~\ref{0502}, we know that there exists a semisimple element $t\in X^\flat_{l,l}$
for some $l\leq\min(n,n')$ such that $s=t^\flat_{n,\epsilon}$ for some $\epsilon$ and $-s'=t^\flat_{n'}$
up to conjugation.

Now if $s$ is given, then $s'$ is unique determined by $s$.
On the other hand, if $s'$ is given, then there are at most two possible $s$, namely
$s_1=t^\flat_{n,+}$ and $s_2=t^\flat_{n,-}$.
Note that from the definition, we have $s_2=s_1^\dag$ and $s_1=s_2^\dag$.
\end{proof}

The full centralizer of $\Sp_{2n}$ in $\Sp_{2n(2n'+1)}$ is in fact $\rmO_{2n+1}$ but not $\SO_{2n+1}$.
This is the reason that (ii) of the corollary is not really satisfactory.
The correspondence for the dual pair $(\Sp_{2n},\rmO_{2n'+1})$ will be discussed in Subsection~\ref{0809}.

\begin{exam}
We know that
\[
\omega^\psi_{\Sp_{2n},\SO_1}=\chi_{\frac{q^n+1}{2}}\otimes{\bf1}_{\SO_1}+\chi_{\frac{q^n-1}{2}}\otimes{\bf1}_{\SO_1}
\]
where $\chi_{\frac{q^n+1}{2}}\in\cale(\Sp_{2n}(q))_{s}$ is an irreducible character of degree $\frac{q^n+1}{2}$
such that $C_{G^*}(s)=\rmO^+_{2n}(q)$,
and $\chi_{\frac{q^n-1}{2}}\in\cale(\Sp_{2n}(q))_{s^\dag}$ is an irreducible character of degree $\frac{q^n-1}{2}$
such that $C_{G^*}(s^\dag)=\rmO^-_{2n}(q)$.
\end{exam}

\section{Lusztig correspondence and Howe correspondence}\label{0806}

\subsection{Lusztig parametrization for classical groups}

For a symplectic group or an orthogonal group $\bfG$,
let $\call_1$ denote a bijection from $\cals_\bfG\rightarrow\cale(G)_1$ denoted by $\Lambda\mapsto\rho_\Lambda$
given in Proposition~\ref{0232}.
Of course, such a parametrization of unipotent characters also exists for a general linear group
or a unitary group.
Combining $\call_1$ (for $\bfG^{(0)}\times\bfG^{(1)}\times\bfG^{(2)}$)
and the inverse of $\Xi_s$ in (\ref{0505}) or (\ref{0609}):
\[
\cals_{\bfG^{(0)}}\times\cals_{\bfG^{(1)}}\times\cals_{\bfG^{(2)}}
\rightarrow \cale(G^{(0)}\times G^{(1)}\times G^{(2)})_1
\rightarrow \cale(G)_s
\]
we obtain a bijection
\begin{align}
\begin{split}
\call_s\colon\cals_{\bfG^{(0)}}\times\cals_{\bfG^{(1)}}\times\cals_{\bfG^{(2)}} &\rightarrow\cale(G)_s \\
(x,\Lambda_1,\Lambda_2) &\mapsto  \rho_{x,\Lambda_1,\Lambda_2}.
\end{split}
\end{align}
Such a bijection $\call_s$ is also called a \emph{modified Lusztig correspondence}.

\subsection{Howe correspondence for $(\Sp_{2n},\rmO^\epsilon_{2n'})$}

Now we consider the dual pair $(\bfG,\bfG')=(\Sp_{2n},\rmO^\epsilon_{2n'})$.
Theorem~\ref{0606} can be reformulated as follows:

\begin{prop}\label{0802}
Let $(\bfG,\bfG')=(\Sp_{2n},\rmO_{2n'}^\epsilon)$ where $\epsilon=+$ or $-$,
and $s\in G^*$, $s'\in(G'^*)^0$ semisimple.
Let $\call_s$ and $\call_{s'}$ be any modified Lusztig correspondences for $\bfG$ and $\bfG'$ respectively.
Then one of
\[
(\rho_{x,\Lambda_1,\Lambda_2},\rho_{x',\Lambda_1',\Lambda_2'}),\quad
(\rho_{x,\Lambda_1,\Lambda_2},\rho_{x',\Lambda_1'^\rmt,\Lambda_2'}),\quad
(\rho_{x,\Lambda_1,\Lambda_2},\rho_{x',\Lambda_1',\Lambda_2'^\rmt}),\quad
(\rho_{x,\Lambda_1,\Lambda_2},\rho_{x',\Lambda_1'^\rmt,\Lambda_2'^\rmt})
\]
occurs in $\Theta_{\bfG,\bfG'}^\psi$ if and only if
\begin{itemize}
\item $s^{(0)}=s'^{(0)}$ (up to conjugation) and $x=x'$,

\item $\Lambda_1=\Lambda'_1$ or $\Lambda_1=\Lambda_1'^\rmt$, and

\item $(\Lambda_2,\Lambda'_2)$ or $(\Lambda_2,\Lambda_2'^\rmt)$ is in $\calb_{\bfG^{(2)},\bfG'^{(2)}}$.
\end{itemize}
\end{prop}

In general, the modified Lusztig correspondences $\call_s,\call_{s'}$ are not uniquely determined.
However, in \cite{pan-ambiguity} we show that the correspondences $\call_s,\call_{s'}$
can be properly chosen such that they are compatible with the parabolic induction and also compatible with
the theta correspondence for $(\bfG,\bfG')$.
The following is \cite{pan-ambiguity} theorem 9.12 which removes the ambiguity in Proposition~\ref{0802}:

\begin{prop}\label{0808}
Let $(\bfG,\bfG')=(\Sp_{2n},\rmO_{2n'}^\epsilon)$ where $\epsilon=+$ or $-$,
and $s\in G^*$, $s'\in(G'^*)^0$ semisimple.
Let $\call_s,\call_{s'}$ be the modified Lusztig correspondences for $\bfG,\bfG'$ respectively
given in \cite{pan-ambiguity}.
Then $(\rho_{x,\Lambda_1,\Lambda_2},\rho_{x',\Lambda_1',\Lambda_2'})\in \Theta_{\bfG,\bfG'}^\psi$
if and only if
\begin{itemize}
\item $s^{(0)}=s'^{(0)}$ (up to conjugation) and $x=x'$,

\item $\Lambda_1=\Lambda'_1$, and

\item $(\Lambda_2,\Lambda'_2)\in\calb_{\bfG^{(2)},\bfG^{(2)}}$.
\end{itemize}
\end{prop}

\subsection{Howe correspondence for $(\Sp_{2n},\SO_{2n'+1})$}\label{0807}
In this subsection, we consider the dual pair $(\bfG,\bfG')=(\Sp_{2n},\SO_{2n'+1})$.
Now Theorem~\ref{0502} can be reformulated as follows:

\begin{prop}\label{0801}
Let $(\bfG,\bfG')=(\Sp_{2n},\SO_{2n'+1})$, and $s\in G^*$, $s'\in(G'^*)^0$ semisimple.
Let $\call_s$ and $\call_{s'}$ be any modified Lusztig correspondences for $\bfG$ and $\bfG'$ respectively.
Then one of
\[
(\rho_{x,\Lambda_1,\Lambda_2},\rho_{x',\Lambda_1',\Lambda_2'}),\quad
(\rho_{x,\Lambda_1^\rmt,\Lambda_2},\rho_{x',\Lambda_1',\Lambda_2'})
\]
occurs in $\Theta_{\bfG,\bfG'}^\psi$ if and only if
\begin{itemize}
\item $s^{(0)}=-s'^{(0)}$ (up to conjugation) and $x=x'$,

\item $\Lambda_2=\Lambda'_1$, and

\item $(\Lambda_1,\Lambda'_2)$ or $(\Lambda_1^\rmt,\Lambda'_2)$ is in $\calb_{\bfG^{(1)},\bfG'^{(2)}}$.
\end{itemize}
\end{prop}

Again, for the dual pair $(\bfG,\bfG')=(\Sp_{2n},\SO_{2n'+1})$, Lusztig correspondences $\call_s,\call_{s'}$
can be properly chosen such that they are compatible with the parabolic induction and also compatible with
the theta correspondence.
The following result is \cite{pan-ambiguity} theorem~9.9:

\begin{prop}\label{0803}
Let $(\bfG,\bfG')=(\Sp_{2n},\SO_{2n'+1})$, and $s\in G^*$, $s'\in G'^*$ semisimple.
Let $\call_s,\call_{s'}$ be the modified Lusztig correspondences for $\bfG$ and $\bfG'$
respectively given in \cite{pan-ambiguity}.
Then $(\rho_{x,\Lambda_1,\Lambda_2},\rho_{x',\Lambda_1',\Lambda_2'})$ occurs in $\Theta_{\bfG,\bfG'}^\psi$
if and only if
\begin{itemize}
\item $s^{(0)}=-s'^{(0)}$ (up to conjugation) and $x=x'$,

\item $\Lambda_2=\Lambda'_1$, and

\item $(\Lambda_1,\Lambda'_2)\in\calb_{\bfG^{(1)},\bfG^{(2)}}$.
\end{itemize}
\end{prop}

\subsection{Howe correspondence for $(\Sp_{2n},\rmO_{2n'+1})$}\label{0809}
In this subsection, let $\bfG'=\rmO_{2n'+1}$ and $\bfG'^0=\SO_{2n'+1}$.
It is known that  $\bfG'\simeq\bfG'^0\times\{\pm1\}$,
and so we have
\[
\cale(G')=\{\,\rho'\otimes{\bf1},\rho'\otimes\sgn\mid\rho'\in\cale(G'^0)\,\}.
\]
Now for $s'\in(G'^0)^*$ and $\epsilon'=\pm$, we define
\begin{align*}
\cale(G')_{s',\epsilon'}=\begin{cases}
\{\,\rho'\otimes{\bf1}\mid\rho'\in\cale(G'^0)_{s'}\,\}, & \text{if $\epsilon'=+$};\\
\{\,\rho'\otimes\sgn\mid\rho'\in\cale(G'^0)_{s'}\,\}, & \text{if $\epsilon'=-$}.
\end{cases}
\end{align*}
Each $\cale(G')_{s',\epsilon'}$ will be regarded as a \emph{Lusztig series} of $G'$.
It is clear that for $\rho'\in\cale(G')_{s',\epsilon'}$ we have $\rho'|_{G'^0}\in\cale(G'^0)_{s'}$.
Moreover, we have
\[
\cale(G')=\bigcup_{(s')\subset(G'^0)^*,\ \epsilon'=\pm}\cale(G')_{s',\epsilon'}.
\]
For $\rho'\in\cale(\rmO_{2n'+1}(q))$,
because $-1$ is in the center of $\rmO_{2n'+1}(q)$,
we know that $\rho'(-1)=\epsilon_{\rho'}\rho'(1)$ where $\epsilon_{\rho'}=\pm1$.
It is clear that if $\rho'\in\cale(G')_{s',\epsilon'}$, then $\epsilon_{\rho'}=\epsilon'$.

Similarly for $\rho\in\cale(\Sp_{2n}(q))$,
we have $\rho(-1)=\epsilon_\rho\rho(1)$ where $\epsilon_\rho=\pm1$.

\begin{lem}\label{0810}
Let $s\in\SO_{2n+1}(q)$ be semisimple,
$\rho\in\cale(\Sp_{2n}(q))_s$ and $\rho''\in\cale(\Sp_{2n}(q))_{s^\dag}$
where $s^\dag$ is given as in Subsection~\ref{0712}.
Write $\rho=\rho_{x,\Lambda_1,\Lambda_2}$ and $\rho''=\rho_{x'',\Lambda''_1,\Lambda_2''}$
under any modified Lusztig correspondences $\call_s,\call_{s^\dag}$.
Suppose that $x=x''$ and $\Lambda_2=\Lambda_2''$.
Then $\epsilon_\rho=-\epsilon_{\rho''}$.
\end{lem}
\begin{proof}
By the similar argument in \cite{waldspurger} \S 4.4 and \S 4.11, we can see that
\[
\epsilon_\rho\epsilon_{\rho''}=(-1)^{\frac{{\rm def}(\Lambda_1)}{2}+\frac{{\rm def}(\Lambda_1'')}{2}}.
\]
From the definition of $s^\dag$ we know that $\bfG^{(1)}(s)=\rmO^{\epsilon^{(1)}}_{2n^{(1)}}$
and $\bfG^{(1)}(s^\dag)=\rmO^{-\epsilon^{(1)}}_{2n^{(1)}}$ for some $\epsilon^{(1)}$ and $n^{(1)}$.
From the definition of $\cals_\bfG$ in Subsection~\ref{0234},
we see that one of ${\rm def}(\Lambda_1),{\rm def}(\Lambda''_1)$ is $\equiv 0\pmod 4$ and the other is
$\equiv 2\pmod 4$.
Therefore $\frac{{\rm def}(\Lambda_1)}{2}+\frac{{\rm def}(\Lambda_1'')}{2}$ must be odd.
Note that if we choose another $\call_s$ or $\call_{s^\dag}$,
from \cite{pan-ambiguity}, we know that $\Lambda_1$ might be changed to $\Lambda_1^\rmt$,
or $\Lambda''_1$ might be changed to $\Lambda_1''^\rmt$,
however the parity of $\frac{{\rm def}(\Lambda_1)}{2}+\frac{{\rm def}(\Lambda_1'')}{2}$ is unchanged.
\end{proof}

Now we have the following refinement of Corollary~\ref{0711}:
\begin{prop}\label{0811}
Let $(\bfG,\bfG')=(\Sp_{2n},\rmO_{2n'+1})$.
Suppose that $(\rho,\rho')\in\Theta^\psi_{\bfG,\bfG'}$, 
and $\rho\in\cale(G)_s$, $\rho'\in\cale(G')_{s',\epsilon'}$.
Then $\Theta^\psi_{\bfG'}(\rho)\subset\cale(G')_{s',\epsilon'}$,
$\Theta^\psi_{\bfG}(\rho')\subset\cale(G)_s$.
\end{prop}
\begin{proof}
Let $V$ be a $2n$-dimensional symplectic space, and let $V'$ be a $(2n'+1)$-dimensional quadratic space
so that $G=\Sp(V)$ and $G'=\rmO(V')$.
Then we have a homomorphism
\[
\iota\colon\Sp(V)\times\rmO(V')\rightarrow\Sp(V\otimes V').
\]
Note that $\iota(-1_{V},1_{V'})=\iota(1_V,-1_{V'})=-1_{V\otimes V'}$
where $-1_{V}\colon V\rightarrow V$ given by $v\mapsto -v$ for $v\in V$, and so on.
Now suppose that $\rho\in\cale(G)_s$, $\rho'\in\cale(G')_{s',\epsilon'}$,
and $(\rho,\rho')\in\Theta^\psi_{\bfG,\bfG'}$.
Then we have
\[
\epsilon_\rho\rho(1_V)\rho'(1_{V'})
=\rho(-1_V)\rho'(1_{V'})
=\rho(1_V)\rho'(-1_{V'})
=\epsilon_{\rho'}\rho(1_V)\rho'(1_{V'}),
\]
i.e., $\epsilon_\rho=\epsilon_{\rho'}$.

Let $\rho'^0=\rho'|_{G'^0}$.
Then we have $\rho'^0\in\cale(G'^0)_{s'}$.
Suppose that $\rho''\in\Theta^\psi_\bfG(\rho')$, i.e., $(\rho'',\rho')\in\Theta^\psi_{\bfG,\bfG'}$.
Then we have $(\rho'',\rho'^0)\in\Theta^\psi_{\bfG,\bfG'^0}$.
Then by Corollary~\ref{0711}, we know that $\rho''\in\cale(G)_s\cup\cale(G)_{s^\dag}$.
But we have $\epsilon_{\rho''}=\epsilon_{\rho'}=\epsilon_\rho$.
Then we have $\rho''\in\cale(G)_s$ by Lemma~\ref{0810}.
Hence we conclude that $\Theta^\psi_\bfG(\rho')\subset\cale(G)_s$.

Next suppose that $\rho'''\in\Theta^\psi_{\bfG'}(\rho)$.
Then we have $\epsilon_{\rho'''}=\epsilon_\rho=\epsilon_{\rho'}=\epsilon'$,
i.e., $\rho'''$ is in $\cale(G')_{s',\epsilon'}$.
Therefore, we have $\Theta^\psi_{\bfG'}(\rho)\subset\cale(G')_{s',\epsilon'}$.
\end{proof}

\begin{prop}
Let $(\bfG,\bfG')=(\Sp_{2n},\rmO_{2n'+1})$,
and let $\rho\in\cale(G)_s$ and $\rho'\in\cale(G')_{s',\epsilon'}$ for some semisimple elements $s\in G^*$
and $s'\in (G'^0)^*$, and some $\epsilon'=\pm$.
Write $\Xi_s(\rho)=\rho^{(0)}\otimes\rho^{(1)}\otimes\rho^{(2)}$,
$\Xi_{s'}(\rho'|_{G'^0})=\rho'^{(0)}\otimes\rho'^{(1)}\otimes\rho'^{(2)}$
where $\Xi_s,\Xi_{s'}$ are the modified Lusztig correspondences for $\bfG$ and $\bfG'^0$
respectively given in \cite{pan-ambiguity}.
Then $(\rho,\rho')$ occurs in $\Theta_{\bfG,\bfG'}^\psi$ if and only if
\begin{itemize}
\item $s^{(0)}=-s'^{(0)}$ (up to conjugation), and $\rho^{(0)}=\rho'^{(0)}$;

\item $\bfG^{(2)}\simeq\bfG'^{(1)}$ and $\rho^{(2)}=\rho'^{(1)}$;

\item $(\rho^{(1)},\rho'^{(2)})$ occurs in $\Theta_{\bfG^{(1)},\bfG'^{(2)},1}$;

\item $\epsilon_\rho=\epsilon'$.
\end{itemize}
\end{prop}
\begin{proof}
Now $\bfG'^0=\SO_{2n'+1}$ and let $\rho'^0=\rho'|_{G'^0}$.
First suppose that $(\rho,\rho')$ occurs in $\Theta^\psi_{\bfG,\bfG'}$.
Then $(\rho,\rho'^0)$ occurs in $\Theta^\psi_{\bfG,\bfG'^0}$.
By Proposition~\ref{0803}, we have $s^{(0)}=-s'^{(0)}$ (up to conjugation), and $\rho^{(0)}=\rho'^{(0)}$;
$\bfG^{(2)}\simeq\bfG'^{(1)}$ and $\rho^{(2)}=\rho'^{(1)}$;
$(\rho^{(1)},\rho'^{(2)})$ occurs in $\Theta_{\bfG^{(1)},\bfG'^{(2)},1}$.
Moreover, we have $\epsilon_\rho=\epsilon'$ from the proof of Proposition~\ref{0811}.

Conversely, suppose the conditions in the proposition are all satisfied.
By Proposition~\ref{0803}, we see that $(\rho,\rho'^0)$ occurs in $\Theta^\psi_{\bfG,\bfG'^0}$.
Then from the proof of Proposition~\ref{0811}, we have $(\rho,\rho'^0\otimes\epsilon_\rho)$ occurs in
$\Theta^\psi_{\bfG,\bfG'}$.
Finally, the condition $\epsilon_\rho=\epsilon'$ implies that $\rho'^0\otimes\epsilon_\rho\simeq\rho'$
and hence $(\rho,\rho')$ occurs in $\Theta^\psi_{\bfG,\bfG'}$.
\end{proof}

\bibliography{refer}

\providecommand{\bysame}{\leavevmode\hbox to3em{\hrulefill}\thinspace}
\providecommand{\MR}{\relax\ifhmode\unskip\space\fi MR }
\providecommand{\MRhref}[2]{%
  \href{http://www.ams.org/mathscinet-getitem?mr=#1}{#2}
}
\providecommand{\href}[2]{#2}
\begin{thebibliography}{AMR96}

\bibitem[AM93]{adams-moy}
J.~Adams and A.~Moy, \emph{Unipotent representations and reductive dual pairs
  over finite fields}, Trans. Amer. Math. Soc. \textbf{340} (1993), 309--321.

\bibitem[AMR96]{amr}
A.-M. Aubert, J.~Michel, and R.~Rouquier, \emph{Correspondance de {H}owe pour
  les groupes r\'eductifs sur les corps finis}, Duke Math. J. \textbf{83}
  (1996), 353--397.

\bibitem[Car85]{carter-finite}
R.~Carter, \emph{Finite groups of {L}ie type, conjugacy classes and complex
  characters}, John Wiley \& Sons, England, 1985.

\bibitem[DL76]{dl}
P.~Deligne and G.~Lusztig, \emph{Representations of reductive groups over
  finite fields}, Ann. of Math. \textbf{103} (1976), 103--161.

\bibitem[DM90]{DM-Lusztig}
F.~Digne and J.~Michel, \emph{On {L}usztig's parametrization of characters of
  finite groups of {L}ie type}, Ast\'erisque \textbf{181--182} (1990),
  113--156.

\bibitem[DM91]{DM}
\bysame, \emph{Representations of finite groups of {L}ie type}, Cambridge
  University Press, Cambridge, 1991.

\bibitem[FS90]{FS}
P.~Fong and B.~Srinivasan, \emph{Brauer trees in classical groups}, J. Algebra
  \textbf{131} (1990), 179--225.

\bibitem[G{\'e}r77]{gerardin}
P.~G{\'e}rardin, \emph{Weil representations associated to finite fields}, J.
  Algebra \textbf{46} (1977), 54--101.

\bibitem[GM20]{GM-guide}
M.~Geck and G.~Malle, \emph{The character theory of finite groups of {L}ie
  type---a guided tour}, Cambridge Studies in Advanced Mathematics, no. 187,
  Cambridge University Press, Cambridge, 2020.

\bibitem[GP00]{Geck-Pfeiffer}
M.~Geck and G.~Pfeiffer, \emph{Characters of finite {C}oxeter groups and
  {I}wahori-{H}ecke algebras}, Oxford University Press, New York, 2000.

\bibitem[Lus77]{lg}
G.~Lusztig, \emph{Irreducible representations of finite classical groups},
  Invent. Math. \textbf{43} (1977), 125--175.

\bibitem[Lus78]{lg-CBMS}
\bysame, \emph{Representations of finite {C}hevalley groups}, C.B.M.S.~Regional
  Conference Series in Math., no.~39, American Mathematical Society,
  Providence, 1978.

\bibitem[Lus81]{lg-symplectic}
\bysame, \emph{Unipotent characters of the symplectic and odd orthogonal groups
  over a finite field}, Invent. Math. \textbf{64} (1981), 263--296.

\bibitem[Lus82]{lg-orthogonal}
\bysame, \emph{Unipotent characters of the even orthogonal groups over a finite
  field}, Trans. Amer. Math. Soc. \textbf{272} (1982), 733--751.

\bibitem[Lus84]{lusztig-book}
\bysame, \emph{Characters of reductive groups over a finite field}, Ann. Math.
  Stud., no. 107, Princeton University Press, Princeton, 1984.

\bibitem[Pan16]{pan-odd}
S.-Y. Pan, \emph{Weil representation of finite symplectic groups and finite odd
  orthogonal groups}, J. Algebra \textbf{453} (2016), 291--324.

\bibitem[Pan19a]{pan-finite-unipotent}
\bysame, \emph{Howe correspondence of unipotent characters for a finite
  symplectic/even-orthogonal dual pair}, arXiv:1901.00623 (2019).

\bibitem[Pan19b]{pan-chain01}
\bysame, \emph{Supercuspidal representations and preservation principle of
  theta correspondence}, J. Reine Angew. Math. \textbf{750} (2019), 1--52.

\bibitem[Pan21]{pan-uniform}
\bysame, \emph{Uniform projection of the {W}eil character}, Adv. Math.
  \textbf{381} (2021), 107624, 80 pp.

\bibitem[Pan22]{pan-ambiguity}
\bysame, \emph{On the unicity and the ambiguity of {L}usztig parametrizations
  for finite classical groups}, arXiv:2206.04900 (2022).

\bibitem[Ser77]{sr-finite}
J.-P. Serre, \emph{Linear representations of finite groups}, Springer-Verlag,
  Berlin-Heidelberg-New York, 1977.

\bibitem[Sri79]{srinivasan}
B.~Srinivasan, \emph{Weil representations of finite classical groups}, Invent.
  Math. \textbf{51} (1979), 143--153.

\bibitem[Wal04]{waldspurger}
J.-L. Waldspurger, \emph{Une conjecture de {L}usztig pour les groupes
  classiques}, M\'em. Soc. Math. Fr. (N.S.) (2004), no.~96, vi+166 pp.

\end{thebibliography}
\bibliographystyle{amsalpha}

\end{document}